# Quantile Processes for Semi and Nonparametric Regression

**Shih-Kang Chao**[*]

*Department of Statistics, Purdue University, West Lafayette, IN 47906.*
*e-mail:* skchao74@purdue.edu

**Stanislav Volgushev**[†]

*Department of Statistical Sciences, University of Toronto, Sidney Smith Hall, 100 St George St, Toronto, ON M5S 3G3.*
*e-mail:* volgushe@utstat.toronto.edu

**Guang Cheng**[‡]

*Department of Statistics, Purdue University, West Lafayette, IN 47906.*
*e-mail:* chengg@purdue.edu

**Abstract:** A collection of quantile curves provides a complete picture of conditional distributions. A properly centered and scaled version of estimated curves at various quantile levels gives rise to the so-called quantile regression process (QRP). In this paper, we establish weak convergence of QRP in a general series approximation framework, which includes linear models with increasing dimension, non-parametric models and partial linear models. An interesting consequence is obtained in the last class of models, where parametric and non-parametric estimators are shown to be asymptotically independent. Applications of our general process convergence results include the construction of non-crossing quantile curves and the estimation of conditional distribution functions. As a result of independent interest, we obtain a series of Bahadur representations with exponential bounds for tail probabilities of all remainder terms. Bounds of this kind are potentially useful in analyzing statistical inference procedures under divide-and-conquer setup.

**MSC 2010 subject classifications:** Primary 62F12, 62G20, 62G08.
**Keywords and phrases:** Bahadur representation, Quantile regression process, Semi/nonparametric model, Series estimation.

## Contents

[*]Partially supported by Office of Naval Research (ONR N00014-15-1-2331).
[†]During that time the second author was supported by the Sonderforschungsbereich "Statistical modelling of nonlinear dynamic processes" (SFB 823), Teilprojekt (C1), of the Deutsche Forschungsgemeinschaft.
[‡]Research Sponsored by NSF CAREER Award DMS-1151692, DMS-1418042, and Office of Naval Research (ONR N00014-15-1-2331).

## 1. Introduction

Quantile regression is widely applied in various scientific fields such as economics (Koenker and Hallock, 2001), biology (Briollais and Durrieu, 2014) and ecology (Cade and Noon, 2003). By focusing on *a collection* of conditional quantiles instead of a single conditional mean, quantile regression allows to describe the impact of predictors on the entire conditional distribution of the response. A properly scaled and centered version of these estimated curves form an underlying (conditional) *quantile regression process* (QRP, see Section 2 for a formal definition). The weak convergence of QRP is useful in developing statistical inference procedures, such as hypothesis testing on Hadamard differentiable $M$ and $L$ estimators (Fernholz, 1983, Chapter 7), testing on conditional distributions (Bassett and Koenker, 1982) and Wilcoxon test (van der Vaart and Wellner, 1996, Example 3.9.18). Applications in econometrics include detection of treatment effect on the conditional distribution after an intervention (Koenker and Xiao, 2002; Qu and Yoon, 2015) and testing Gini indices (Barrett and Donald, 2009). Please see Remark 4.2 for more details.

The asymptotic behavior of QRP depends on the model imposed for quantile regression. Existing literature on QRP is either concerned with models of *fixed* dimension (Koenker and Xiao, 2002; Angrist et al., 2006), or with a linearly interpolated version based on kernel smoothing (Qu and Yoon, 2015). However, this excludes many important cases such as linear models with growing dimension and partial linear models. For such models, establishing weak convergence of QRP becomes non-trivial since classical Donsker theorems (e.g. those in van der Vaart and Wellner (1996)) may not be directly applied. An additional challenge for partial linear models comes from the fact that their parametric and non-parametric components converge at different rates.

In this paper we consider a general model which is of the following (approximate) form

$$Q(x;\tau) \approx \mathbf{Z}(x)^\top \boldsymbol{\gamma}_n(\tau), \tag{1.1}$$

where $Q(x;\tau)$ denotes the $\tau$-th quantile of the distribution of $Y$ conditional on $X = x \in \mathbb{R}^d$ and $\mathbf{Z}(x) \in \mathbb{R}^m$ is a transformation vector of $x$. As noted by Belloni et al. (2016), the above framework incorporates a variety of estimation procedures such as parametric (Koenker and Bassett, 1978), non-parametric (He and Shi, 1994) and semi-parametric (He and Shi, 1996) ones. For example, $\mathbf{Z}(x) = x$ corresponds to a linear model (with potentially increasing dimension), while $\mathbf{Z}(x)$ can be chosen as powers, trigonometrics or local polynomials in the non-parametric basis expansion (where $m$ diverges at a proper rate). Partially linear and additive models are also covered by (1.1). Therefore, our weak convergence results are developed in a very broad context.

Models that can be expressed in the form (1.1) were previously studied by Belloni et al. (2016) in a very general setting, which we also consider here. In the following, we provide a detailed description of the main contributions in the present paper, and compare them with Belloni et al. (2016).

1. **Partially linear models**: A key result in the present paper is obtained for partially linear models

$$Q(X;\tau) = V^\top \boldsymbol{\alpha}(\tau) + h(W;\tau), \tag{1.2}$$

where $X = (V^\top, W)^\top \in \mathbb{R}^{k+k'}$, $\boldsymbol{\alpha}(\tau)$ is an unknown Euclidean vector and $h(W;\tau)$ is an unknown smooth function. Here, $k$ and $k'$ are both fixed. In the spirit of (1.1), we can estimate $(\boldsymbol{\alpha}(\tau), h(\cdot;\tau))$ based on the following series approximation

$$h(W;\tau) \approx \widetilde{\mathbf{Z}}(W)^\top \boldsymbol{\beta}_n^\dagger(\tau),$$

where $\widetilde{\mathbf{Z}}(W)$ is a transformation vector of $W$. We provide *joint* asymptotic results for the parametric and non-parametric part in partially linear models [see Section 3 and Section 5.3], with a $n^{1/2}$ scaling for the parametric part and a scaling with slower rate for the non-parametric part [see Theorem 3.1 and Theorem 5.4].

To the best of our knowledge, this is the first time that the joint asymptotic as processes in $\tau$ for quantile regression is established – in fact, even the pointwise result is new. In particular, we prove that the "joint asymptotics phenomenon" discovered by Cheng and Shang (2015) even holds for non-smooth loss functions with multivariate nonparametric covariates.

This joint asymptotic result does not follow directly from the results of Belloni et al. (2016), because of the specific centering sequence (defined in their equation (2.2)) they consider and the matrix $J^{-1}(u)$ in their Theorem 2 where $J(u)$ is non-diagonal with increasing dimension. To derive our Theorem 3.1, it is necessary to choose an appropriate centering sequence (see Remark 3.3), apply our new Bahadur representations in Section 5, and provide a detailed analysis of the matrix $J^{-1}(u)$.

2. **Centering and tail bounds on remainder terms in Bahadur representations**: Theorem 2 in Belloni et al. (2016) provides a Bahadur representation for the estimated coefficients $\widehat{\boldsymbol{\gamma}}_n$ centering at $\boldsymbol{\beta}_n$ with an $O_P$ term on the remainders, where $\boldsymbol{\beta}_n$ minimizes a QR series approximation problem [see equation (2.2) in their paper]. In Section 5 we provide similar expansions, with a main difference that we allow for more general centering sequences that satisfy certain approximation conditions. This is important in Section 3 of our paper where we consider partially linear models. Moreover, we provide explicit *exponential* tail bounds on corresponding remainders which is a somewhat stronger result compared to the $O_P$ bounds in Belloni et al. (2016). This result is, for instance, utilized in Volgushev et al. (2017). The findings in that paper cannot be obtained from the $O_P$ bounds in Belloni et al. (2016).

3. **Approximation by a sequence of Gaussian processes v.s. convergence to a fixed limiting process**: All results in Belloni et al. (2016) which are uniform in the quantile index $\tau$ are stated in terms of approximating the quantile regression process and weighted versions thereof by a *sequence* of Gaussian processes which depend on $n$ [see their Theorem 5, Theorem 11, Theorem 12]. In contrast, we show that there exists a single Gaussian process which is the (weak) limit of the leading term in the Bahadur representation. Showing convergence to this weak limit requires proving asymptotic tightness of the leading term, which is a major challenge in our proof [see Section A.4] and does not follow from the approximation by a series of Gaussian processes as in Belloni et al. (2016). This is also a key ingredient in our Section 4 where we utilize the functional delta method together with compact differentiability of the rearrangement operator [established in Chernozhukov et al. (2010)]. Note that the application of the delta method requires *convergence to a fixed limit*, which does not follow directly from the results in Belloni et al. (2016). On the other hand, Belloni et al. (2016) provide approximations which are uniform in $x$ and $\tau$ while we only consider results that are pointwise in $x$.

4. **New bounds for local basis functions**: Last but not the least, in Sections 2.2 and 5.2, we provide results for models with "local basis structure" (for instance $B$-splines). For such basis functions, we show that the conditions on model dimension can be relaxed from $m^4 = o(n^{1-\varepsilon})$ [required by Theorem 12 of Belloni et al. (2016)] to $m^2(\log n)^6 = o(n)$ in the case of $B$-splines [see the discussion below Assumption (B1)].

Given the discussions above, we would also like to point out that Belloni et al. (2016) discuss other aspects such as bootstrap approximations which are not covered in our paper. In summary, both Belloni et al. (2016) and the present paper consider the same model setup, but focus on different aspects of the resulting theory, and none of the two papers is more general than the other.

The rest of this paper is organized as follows. Section 2 presents the weak convergence of QRP under general series approximation framework. Section 3 discusses the QRP in quantile partial linear models. As an application of our weak convergence theory, Section 4 considers various functionals of the quantile regression process. A detailed discussion on our novel Bahadur representations is given in Section 5, and all proofs are deferred to the appendix.

**Notation.** Denote $\{(X_i, Y_i)\}_{i=1}^n$ i.i.d. samples in $\mathcal{X} \times \mathbb{R}$ where $\mathcal{X} \subset \mathbb{R}^d$. Here, the distribution of $(X_i, Y_i)$ and the dimension $d$ can depend on $n$, i.e. triangular arrays. For brevity, let $\mathbf{Z} = \mathbf{Z}(X)$ and $\mathbf{Z}_i = \mathbf{Z}(X_i)$. Define the empirical measure of $(Y_i, \mathbf{Z}_i)$ by $\mathbb{P}_n$, and the true underlying measure by $P$ with the corresponding expectation as $\mathbb{E}$. Note that the measure $P$ depends on $n$ for triangular array cases, but this dependence is omitted in the notation. Denote by $\|\mathbf{b}\|$ the $L^2$-norm of a vector $\mathbf{b}$. $\lambda_{\min}(A)$ and $\lambda_{\max}(A)$ are the smallest and largest eigenvalue of a matrix $A$. $\mathbf{0}_k$ denotes a $k$-dimensional 0 vector, and $I_k$ be the $k$-dimensional identity matrix for $k \in \mathbb{N}$. Define

$$\rho_\tau(u) := (\tau - \mathbf{1}(u \leq 0))u,$$

where $\mathbf{1}(\cdot)$ is the indicator function. $\mathcal{C}^\eta(\mathcal{X})$ denotes the class of $\eta$-continuously differentiable functions on a set $\mathcal{X}$. $\mathcal{C}(0,1)$ denotes the class of continuous functions defined on $(0,1)$. Define

$$\psi(Y_i, \mathbf{Z}_i; \mathbf{b}, \tau) := \mathbf{Z}_i(\mathbf{1}\{Y_i \leq \mathbf{Z}_i^\top \mathbf{b}\} - \tau), \quad \mu(\mathbf{b}, \tau) := \mathbb{E}[\psi(Y_i, \mathbf{Z}_i; \mathbf{b}, \tau)] = \mathbb{E}[\mathbf{Z}_i\{F_{Y|X}(\mathbf{Z}_i^\top \mathbf{b}|X) - \tau\}],$$

and for a vector $\boldsymbol{\gamma}_n(\tau) \in \mathbb{R}^m$, we define the following quantities

$$g_n := g_n(\boldsymbol{\gamma}_n) := \sup_{\tau \in \mathcal{T}} \|\mu(\boldsymbol{\gamma}_n(\tau), \tau)\| = \sup_{\tau \in \mathcal{T}} \left\|\mathbb{E}[\mathbf{Z}_i\{F_{Y|X}(\mathbf{Z}_i^\top \boldsymbol{\gamma}_n(\tau)|X) - \tau\}]\right\| \tag{1.3}$$

Let $\mathcal{S}^{m-1} := \{\mathbf{u} \in \mathbb{R}^m : \|\mathbf{u}\| = 1\}$ denote the unit sphere in $\mathbb{R}^m$. For a set $\mathcal{I} \subset \{1, ..., m\}$, define

$$\mathbb{R}_\mathcal{I}^m := \{\mathbf{u} = (u_1, ..., u_m)^\top \in \mathbb{R}^m : u_j \neq 0 \text{ if and only if } j \in \mathcal{I}\}$$
$$\mathcal{S}_\mathcal{I}^{m-1} := \{\mathbf{u} = (u_1, ..., u_m)^\top \in \mathcal{S}^{m-1} : u_j \neq 0 \text{ if and only if } j \in \mathcal{I}\}$$

Finally, consider the class of functions

$$\Lambda_c^\eta(\mathcal{X}, \mathcal{T}) :=$$
$$\left\{ f_\tau \in \mathcal{C}^{\lfloor \eta \rfloor}(\mathcal{X}) : \tau \in \mathcal{T}, \sup_{|\boldsymbol{j}| \leq \lfloor \eta \rfloor} \sup_{x, \tau \in \mathcal{T}} |D^{\boldsymbol{j}} f_\tau(x)| \leq c, \sup_{|\boldsymbol{j}| = \lfloor \eta \rfloor} \sup_{x \neq y, \tau \in \mathcal{T}} \frac{|D^{\boldsymbol{j}} f_\tau(x) - D^{\boldsymbol{j}} f_\tau(y)|}{\|x - y\|^{\eta - \lfloor \eta \rfloor}} \leq c \right\}, \tag{1.4}$$

where $\lfloor \eta \rfloor$ denotes the integer part of a real number $\eta$, and $|\boldsymbol{j}| = j_1 + ... + j_d$ for $d$-tuple $\boldsymbol{j} = (j_1, ..., j_d)$. For simplicity, we sometimes write $\sup_\tau(\inf_\tau)$ and $\sup_x(\inf_x)$ instead of $\sup_{\tau \in \mathcal{T}}(\inf_{\tau \in \mathcal{T}})$ and $\sup_{x \in \mathcal{X}}(\inf_{x \in \mathcal{X}})$ throughout the paper.

## 2. Weak Convergence Results

In this section, we first present our weak convergence results of QRP in a general series approximation framework that covers linear models with increasing dimension, nonparametric models and partial linear models. Furthermore, we demonstrate that the use of polynomial splines with local support, such as $B$-splines, significantly weakens the sufficient conditions required in the above general framework.

### 2.1. General Series Estimator

Consider a general series estimator $\widehat{Q}(x;\tau) := \widehat{\gamma}(\tau)^\top \mathbf{Z}(x)$, where for each fixed $\tau$

$$\widehat{\gamma}(\tau) := \underset{\gamma \in \mathbb{R}^m}{\operatorname{argmin}} \sum_{i=1}^n \rho_\tau(Y_i - \gamma^\top \mathbf{Z}_i), \tag{2.1}$$

and $m$ is allowed to grow as $n \to \infty$, and assume the following conditions:

(A1) Assume that $\|\mathbf{Z}_i\| \leq \xi_m = O(n^b)$ almost surely with $b > 0$, and that $1/M \leq \lambda_{\min}(\mathbb{E}[\mathbf{ZZ}^\top]) \leq \lambda_{\max}(\mathbb{E}[\mathbf{ZZ}^\top]) \leq M$ holds uniformly in $n$ for some fixed constant $M > 0$.
(A2) The conditional distribution $F_{Y|X}(y|x)$ is twice differentiable w.r.t. $y$. Denote the corresponding derivatives by $f_{Y|X}(y|x)$ and $f'_{Y|X}(y|x)$. Assume that $\bar{f} := \sup_{y,x} |f_{Y|X}(y|x)| < \infty$ and $\overline{f'} := \sup_{y,x} |f'_{Y|X}(y|x)| < \infty$ uniformly in $n$.
(A3) Assume that uniformly in $n$, there exists a constant $f_{\min} > 0$ such that

$$\inf_{\tau \in \mathcal{T}} \inf_x f_{Y|X}(Q(x;\tau)|x) \geq f_{\min}.$$

In the above assumptions, uniformity in $n$ is necessary as we consider triangular arrays. Assumptions (A2) and (A3) are fairly standard in the quantile regression literature. Hence, we only make a few comments on Assumption (A1). In linear models where $\mathbf{Z}(X) = X$ and $m = d$, it holds that $\xi_m \lesssim \sqrt{m}$ if each component of $X$ is bounded almost surely. If $B$-splines $\widetilde{\mathbf{B}}(x)$ defined in Section 4.3 of Schumaker (1981) are adopted, then one needs to use its re-scaled version $\mathbf{B}(x) = m^{1/2}\widetilde{\mathbf{B}}(x)$ as $\mathbf{Z}(x)$ such that (A1) holds (cf. Lemma 6.2 of Zhou et al. (1998)). In this case, we have $\xi_m \asymp \sqrt{m}$. In addition, Assumptions (A1) and (A3) imply that for any sequence of $\mathbb{R}^m$-valued (non-random) functions $\gamma_n(\tau)$ satisfying $\sup_{\tau \in \mathcal{T}} \sup_x |\gamma_n(\tau)^\top \mathbf{Z}(x) - Q(x;\tau)| = o(1)$, the smallest eigenvalues of the matrices

$$\widetilde{J}_m(\tau) := \mathbb{E}[\mathbf{ZZ}^\top f_{Y|X}(\gamma_n(\tau)^\top \mathbf{Z}|X)], \quad J_m(\tau) := \mathbb{E}[\mathbf{ZZ}^\top f_{Y|X}(Q(X;\tau)|X)]$$

are bounded away from zero uniformly in $\tau$ for all $n$.

Define for any $\mathbf{u} \in \mathbb{R}^m$,

$$\chi_{\gamma_n}(\mathbf{u}, \mathbf{Z}) := \sup_{\tau \in \mathcal{T}} \left| \mathbf{u}^\top J_m(\tau)^{-1} \mathbb{E}\Big[\mathbf{Z}_i\big(\mathbf{1}\{Y_i \leq Q(X_i;\tau)\} - \mathbf{1}\{Y_i \leq \mathbf{Z}_i^\top \gamma_n(\tau)\}\big)\Big]\right|.$$

We are now ready to state our weak convergence result for QRP based on the general series estimators.

**Theorem 2.1.** *Suppose (A1)-(A3) hold and $m^3 \xi_m^2 (\log n)^3 = o(n)$. Let $\gamma_n(\cdot) : \mathcal{T} \to \mathbb{R}^m$ be a sequence of functions such that $g_n := g_n(\gamma_n(\tau)) = o(n^{-1/2})$ (see (1.3)), $c_n = c_n(\gamma_n) := \sup_{x,\tau \in \mathcal{T}} |Q(x;\tau) - \mathbf{Z}(x)^\top \gamma_n(\tau)|$ and $mc_n \log n = o(1)$. Then for any $\mathbf{u}_n \in \mathbb{R}^m$ satisfying $\chi_{\gamma_n}(\mathbf{u}_n, \mathbf{Z}) = o(\|\mathbf{u}_n\| n^{-1/2})$ and $\widehat{\gamma}(\tau)$ defined in (2.1),*

$$\mathbf{u}_n^\top(\widehat{\gamma}(\tau) - \gamma_n(\tau)) = -\frac{1}{n}\mathbf{u}_n^\top J_m(\tau)^{-1} \sum_{i=1}^n \mathbf{Z}_i(\mathbf{1}\{Y_i \leq Q(X_i;\tau)\} - \tau) + o_P\left(\frac{\|\mathbf{u}_n\|}{\sqrt{n}}\right) \tag{2.2}$$

*where the remainder term is uniform in $\tau \in \mathcal{T}$. In addition, if the following limit*

$$H(\tau_1, \tau_2; \mathbf{u}_n) := \lim_{n \to \infty} \|\mathbf{u}_n\|^{-2} \mathbf{u}_n^\top J_m^{-1}(\tau_1) \mathbb{E}[\mathbf{ZZ}^\top] J_m^{-1}(\tau_2) \mathbf{u}_n (\tau_1 \wedge \tau_2 - \tau_1 \tau_2) \tag{2.3}$$

*exists for any $\tau_1, \tau_2 \in \mathcal{T}$, then*

$$\frac{\sqrt{n}}{\|\mathbf{u}_n\|} \left( \mathbf{u}_n^\top \widehat{\boldsymbol{\gamma}}(\cdot) - \mathbf{u}_n^\top \boldsymbol{\gamma}_n(\cdot) \right) \rightsquigarrow \mathbb{G}(\cdot) \ in \ \ell^\infty(\mathcal{T}), \tag{2.4}$$

*where $\mathbb{G}(\cdot)$ is a centered Gaussian process with the covariance function $H$ defined as (2.3). In particular, there exists a version of $\mathbb{G}$ with almost surely continuous sample paths.*

The proof of Theorem 2.1 is given in Section A.1. Theorem 2.1 holds under very general conditions. For transformations $\mathbf{Z}$ that have a specific local structure, the assumptions on $m, \xi_m$ can be relaxed considerably. Details are provided in Section 2.2.

In the end, we illustrate Theorem 2.1 in linear quantile regression models with *increasing* dimension, in which $g_n$, $c_n$ and $\chi_{\boldsymbol{\gamma}_n}(\mathbf{u}, \mathbf{Z})$ are trivially zero. As far as we are aware, this is the first quantile process result for linear models with increasing dimension.

**Corollary 2.2.** *(Linear models with increasing dimension) Suppose (A1)-(A3) hold with $\mathbf{Z}(X) = X$ and $Q(x;\tau) = x^\top \boldsymbol{\gamma}_n(\tau)$ for any $x$ and $\tau \in \mathcal{T}$. Assume that $m^3 \xi_m^2 (\log n)^3 = o(n)$. In addition, if $\mathbf{u}_n \in \mathbb{R}^m$ is such that the following limit*

$$H_1(\tau_1, \tau_2; \mathbf{u}_n) := \lim_{n \to \infty} \|\mathbf{u}_n\|^{-2} \mathbf{u}_n^\top J_m^{-1}(\tau_1) \mathbb{E}[XX^\top] J_m^{-1}(\tau_2) \mathbf{u}_n (\tau_1 \wedge \tau_2 - \tau_1 \tau_2) \tag{2.5}$$

*exists for any $\tau_1, \tau_2 \in \mathcal{T}$, then (2.4) holds with the covariance function $H_1$ defined in (2.5). Moreover, by setting $\mathbf{u}_n = x_0$, we have for any fixed $x_0$*

$$\frac{\sqrt{n}}{\|x_0\|} \{\widehat{Q}(x_0; \cdot) - Q(x_0; \cdot)\} \rightsquigarrow \mathbb{G}(x_0; \cdot) \ in \ \ell^\infty(\mathcal{T}),$$

*where $\mathbb{G}(x_0; \cdot)$ is a centered Gaussian process with covariance function $H_1(\tau_1, \tau_2; x_0)$. In particular, there exists a version of $\mathbb{G}$ with almost surely continuous sample paths.*

### 2.2. Local Basis Series Estimator

In this section, we assume that $\mathbf{Z}(\cdot)$ corresponds to a basis expansion with "local" support. Our main motivation for considering this setting is that it allows to considerably weaken assumptions on $m, \xi_m$ made in the previous section. To distinguish such basis functions from the general setting in the previous section, we shall use the notation $\mathbf{B}$ instead of $\mathbf{Z}$. Let $\widehat{\boldsymbol{\beta}}(\tau)$ be defined as

$$\widehat{\boldsymbol{\beta}}(\tau) := \underset{\mathbf{b} \in \mathbb{R}^m}{\operatorname{argmin}} \sum_i \rho_\tau \{Y_i - \mathbf{b}^\top \mathbf{B}_i\}, \tag{2.6}$$

where $\mathbf{B}_i = \mathbf{B}(X_i)$. The notion of "local support" is made precise in the following sense.

(L) For each $x$, the basis vector $\mathbf{B}(x)$ has zeroes in all but at most $r$ *consecutive* entries, where $r$ is fixed. Moreover, $\sup_{x,\tau} \mathbb{E}[|\mathbf{B}(x)^\top \widetilde{J}_m(\tau)^{-1} \mathbf{B}(X)|] = O(1)$.

The above assumption holds for certain choices of basis functions, e.g., univariate $B$-splines.

**Example 2.3.** Let $\mathcal{X} = [0, 1]$, assume that (A2)-(A3) hold and that the density of $X$ over $\mathcal{X}$ is uniformly bounded away from zero and infinity. Consider the space of polynomial splines of order $q$ with $k$ uniformly spaced knots $0 = t_1 < ... < t_k = 1$ in the interval $[0, 1]$. The space of such splines can be represented through linear combinations of the basis functions $B_1, ..., B_{k-q-1}$ with each basis function $B_j$ having support contained in the interval $[t_j, t_{j+q+1}]$. Let $\mathbf{B}(x) := (B_1(x), ..., B_{k-q-1}(x))^\top$. Then the first part of assumption (L) holds with $r = q$. The condition $\sup_{x,\tau} \mathbb{E}[|\mathbf{B}(x)^\top \widetilde{J}_m^{-1}(\tau) \mathbf{B}(X)|] = O(1)$ is verified in the Appendix, see Section A.2.

Condition (L) ensures that the matrix $\widetilde{J}_m(\tau)$ has a band structure, which is useful for bounding the off-diagonal entries of $\widetilde{J}_m^{-1}(\tau)$. See Lemma 6.3 in Zhou et al. (1998) for additional details.

Throughout this section, consider the specific centering

$$\boldsymbol{\beta}_n(\tau) := \operatorname*{argmin}_{\mathbf{b} \in \mathbb{R}^m} \mathbb{E}\big[(\mathbf{B}^\top \mathbf{b} - Q(X;\tau))^2 f_{Y|X}(Q(X;\tau)|X)\big], \tag{2.7}$$

where $\mathbf{B} = \mathbf{B}(X)$. For basis functions satisfying condition (L), assumptions in Theorem 2.1 in the previous section can be replaced by the following weaker version.

(B1) Assume that $\xi_m^4 (\log n)^6 = o(n)$ and letting $\widetilde{c}_n := \sup_{x,\tau} |\boldsymbol{\beta}_n(\tau)^\top \mathbf{B}(x) - Q(x;\tau)|$ with $\widetilde{c}_n^2 = o(n^{-1/2})$, where $\|\mathbf{B}(X_i)\| \leq \xi_m$ almost surely.

Note that the condition $\xi_m^4 (\log n)^6 = o(n)$ in (B1) is less restrictive than $m^3 \xi_m^2 (\log n)^3 = o(n)$ required in Theorem 2.1 under many situations. For instance, in the setting of Example 2.3 where $\xi_m \asymp \sqrt{m}$, we only require $m^2 (\log n)^6 = o(n)$, which is weaker than $m^4 (\log n)^3 = o(n)$ in Theorem 2.1. This improvement is made possible based on the local structure of the spline basis.

In the setting of Example 2.3, bounds on $\widetilde{c}_n$ can be obtained provided that the function $x \mapsto Q(x;\tau)$ is smooth for all $\tau \in \mathcal{T}$. For instance, assuming that $Q(\cdot;\cdot) \in \Lambda_c^\eta(\mathcal{X}, \mathcal{T})$ with $\mathcal{X} = [0,1]$ and integer $\eta$, Remark B.1 shows that $\widetilde{c}_n = O(m^{-\lfloor \eta \rfloor})$. Thus the condition $\widetilde{c}_n^2 = o(n^{-1/2})$ holds provided that $m^{-2\lfloor \eta \rfloor} = o(n^{1/2})$. Since for splines we have $\xi_m \sim m^{1/2}$, this is compatible with the restrictions imposed in assumption (B1) provided that $\eta \geq 1$.

**Theorem 2.4.** *(Nonparametric models with local basis functions) Assume that conditions (A1)-(A3) hold with $\mathbf{Z} = \mathbf{B}$, (L) holds for $\mathbf{B}$ and (B1) for $\boldsymbol{\beta}_n(\tau)$. Assume that the set $\mathcal{I}$ consists of at most $L$ consecutive integers, where $L \geq 1$ is fixed. Then for any $\mathbf{u}_n \in \mathbb{R}_\mathcal{I}^m$, (2.2) holds with $\widehat{\boldsymbol{\gamma}}(\tau)$, $\boldsymbol{\gamma}_n(\tau)$ and $\mathbf{Z}$ being replaced by $\widehat{\boldsymbol{\beta}}(\tau)$, $\boldsymbol{\beta}_n(\tau)$ and $\mathbf{B}$. In addition, if the following limit*

$$\widetilde{H}(\tau_1, \tau_2; \mathbf{u}_n) := \lim_{n \to \infty} \|\mathbf{u}_n\|^{-2} \mathbf{u}_n^\top J_m^{-1}(\tau_1) \mathbb{E}[\mathbf{B}\mathbf{B}^\top] J_m^{-1}(\tau_2) \mathbf{u}_n (\tau_1 \wedge \tau_2 - \tau_1 \tau_2) \tag{2.8}$$

*exists for any $\tau_1, \tau_2 \in \mathcal{T}$, then (2.4) holds with the same replacement as above, and the limit $\mathbb{G}$ is a centered Gaussian process with covariance function $\widetilde{H}$ defined as (2.8). Moreover, for any $x_0$, let $\widehat{Q}(x_0;\tau) := \mathbf{B}(x_0)^\top \widehat{\boldsymbol{\beta}}(\tau)$ and assume that $\widetilde{c}_n = o(\|\mathbf{B}(x_0)\| n^{-1/2})$. Then*

$$\frac{\sqrt{n}}{\|\mathbf{B}(x_0)\|} \{\widehat{Q}(x_0;\cdot) - Q(x_0;\cdot)\} \rightsquigarrow \mathbb{G}(x_0;\cdot) \text{ in } \ell^\infty(\mathcal{T}), \tag{2.9}$$

*where $\mathbb{G}(x_0;\cdot)$ is a centered Gaussian process with covariance function $\widetilde{H}(\tau_1, \tau_2; \mathbf{B}(x_0))$. In particular, there exists a version of $\mathbb{G}$ with almost surely continuous sample paths.*

The proof of Theorem 2.4 is given in Section A.2.

**Remark 2.5.** The proof of Theorem 2.4 and the related Bahadur representation result in Section 5.2 crucially rely on the fact that the elements of $\widetilde{J}_m(\tau)^{-1}$ decay exponentially fast in their distance from the diagonal, i.e. a bound of the form $|(\widetilde{J}_m(\tau)^{-1})_{i,j}| \leq C\gamma^{|i-j|}$ for some $\gamma < 1$. Assumption (L) provides one way to guarantee such a result. We conjecture that similar results can be obtained for more classes of basis functions as long as the entries of $\widetilde{J}_m(\tau)^{-1}$ decay exponentially fast in their distance from suitable subsets of indices in $(j, j') \in \{1, ..., m\}^2$. This kind of result can be obtained for matrices $\widetilde{J}_m(\tau)$ with specific sparsity patterns, see for instance Demko et al. (1984). In particular, we conjecture that such arguments can be applied for tensor product $B$-splines, see Example 1 in Section 5 of Demko et al. (1984). A detailed investigation of this interesting topic is left to future research.

We conclude this section by discussing a special case where the limit in (2.9) can be characterized more explicitly.

**Remark 2.6.** The covariance function $\widetilde{H}$ can be explicitly characterized under $\mathbf{u}_n = \mathbf{B}(x)$ and univariate B-splines $\mathbf{B}(x)$ on $x \in [0,1]$, with an order $r$ and equidistant knots $0 = t_1 < ... < t_k = 1$. Assume additional to (A3) that

$$\sup_{t \in \mathcal{X}, \tau \in \mathcal{T}} \left| \partial_x |_{x=t} f_{Y|X}\big(Q(x;\tau)|x\big) \right| < C, \text{ where } C > 0 \text{ is a constant}, \tag{2.10}$$

and the density $f_X(x)$ for $X$ is bounded above, then under $\widetilde{c}_n = o(\|\mathbf{B}(x_0)\|n^{-1/2})$, (2.9) in Theorem 2.4 can be rewritten as

$$\sqrt{\frac{n}{\mathbf{B}(x_0)^\top \mathbb{E}[\mathbf{BB}^\top]^{-1} \mathbf{B}(x_0)}} \Big( \mathbf{B}(x_0)^\top \widehat{\boldsymbol{\beta}}(\cdot) - Q(x_0; \cdot) \Big) \rightsquigarrow \mathbb{G}(\cdot; x_0) \text{ in } \ell^\infty(\mathcal{T}), \tag{2.11}$$

where the Gaussian process $\mathbb{G}(\cdot; x_0)$ is defined by the following covariance function

$$\widetilde{H}(\tau_1, \tau_2; x_0) = \frac{\tau_1 \wedge \tau_2 - \tau_1 \tau_2}{f_{Y|X}(Q(x_0; \tau_1)|x_0) f_{Y|X}(Q(x_0; \tau_2)|x_0)}.$$

Although we only show the univariate case here, the same arguments are expected to hold for tensor-product B-spline based on the same reasoning. See Section A.2 for a proof of this remark.

## 3. Joint Weak Convergence for Partial Linear Models

In this section, we consider partial linear models of the form

$$Q(X;\tau) = V^\top \boldsymbol{\alpha}(\tau) + h(W;\tau), \tag{3.1}$$

where $X = (V^\top, W^\top)^\top \in \mathbb{R}^{k+k'}$ and $k, k' \in \mathbb{N}$ are fixed. An interesting *joint* weak convergence result is obtained for $(\widehat{\boldsymbol{\alpha}}(\tau), \widehat{h}(w_0;\tau))$ at any fixed $w_0$. More precisely, $\widehat{\boldsymbol{\alpha}}(\tau)$ and $\widehat{h}(w;\tau)$ (after proper scaling and centering) are proved to be asymptotically independent at any fixed $\tau \in \mathcal{T}$. Therefore, the "joint asymptotics phenomenon" first discovered in Cheng and Shang (2015) persists even for *non-smooth* quantile loss functions. Such a theoretical result is practically useful for joint inference on $\boldsymbol{\alpha}(\tau)$ and $h(W;\tau)$; see Cheng and Shang (2015).

Expanding $w \mapsto h(w;\tau)$ in terms of basis vectors $w \mapsto \widetilde{\mathbf{Z}}(w)$, we can approximate (3.1) through the series expansion $\mathbf{Z}(x)^\top \boldsymbol{\gamma}_n^\dagger(\tau)$ by setting $\mathbf{Z}(x) = (v^\top, \widetilde{\mathbf{Z}}(w)^\top)^\top$. In this section, $\widetilde{\mathbf{Z}} : \mathbb{R}^{k'} \to \mathbb{R}^m$ is regarded as a general basis expansion that does not need to satisfy the local support assumptions in the previous section. Estimation is performed in the following form

$$\widehat{\boldsymbol{\gamma}}^\dagger(\tau) = (\widehat{\boldsymbol{\alpha}}(\tau)^\top, \widehat{\boldsymbol{\beta}}^\dagger(\tau)^\top)^\top := \operatorname*{argmin}_{\mathbf{a} \in \mathbb{R}^k, \mathbf{b} \in \mathbb{R}^m} \sum_i \rho_\tau \big(Y_i - \mathbf{a}^\top V_i - \mathbf{b}^\top \widetilde{\mathbf{Z}}(W_i)\big). \tag{3.2}$$

For a theoretical analysis of $\widehat{\boldsymbol{\gamma}}^\dagger$, define population coefficients $\boldsymbol{\gamma}_n^\dagger(\tau) := (\boldsymbol{\alpha}(\tau)^\top, \boldsymbol{\beta}_n^\dagger(\tau)^\top)^\top$, where

$$\boldsymbol{\beta}_n^\dagger(\tau) := \operatorname*{argmin}_{\boldsymbol{\beta} \in \mathbb{R}^m} \mathbb{E}[f_{Y|X}(Q(X;\tau)|X)(h(W;\tau) - \boldsymbol{\beta}^\top \widetilde{\mathbf{Z}}(W))^2] \tag{3.3}$$

similar to (2.7); see Remark 3.3 for additional explanations.

To state our main result, we need to define a class of functions

$$\mathcal{U}_\tau := \Big\{ w \mapsto g(w) \Big| g \text{ measurable and } \mathbb{E}[g^2(W) f_{Y|X}(Q(X;\tau)|X)] < \infty, w \in \mathbb{R}^{k'} \Big\}.$$

For $V \in \mathbb{R}^k$, define for $j = 1, ..., k$,

$$h_{VW,j}(\cdot;\tau) := \operatorname*{argmin}_{g \in \mathcal{U}_\tau} \mathbb{E}[(V_j - g(W))^2 f_{Y|X}(Q(X;\tau)|X)] \tag{3.4}$$

where $V_j$ denotes the $j$-th entry of the random vector $V$. By the definition of $h_{VW,j}$, we have for all $\tau \in \mathcal{T}$ and $g \in \mathcal{U}_\tau$,
$$\mathbb{E}[(V - h_{VW}(W;\tau))g(W)f_{Y|X}(Q(X;\tau)|X)] = \mathbf{0}_k. \tag{3.5}$$

The matrix $A$ is defined as coeffcient matrix of the best series approximation of $h_{VW}(W;\tau)$:
$$A(\tau) := \underset{A}{\operatorname{argmin}} \, \mathbb{E}[f_{Y|X}(Q(X;\tau)|X)\|h_{VW}(W;\tau) - A\widetilde{\mathbf{Z}}(W)\|^2]. \tag{3.6}$$

The following two assumptions are needed in our main results.

(C1) Define $c_n^\dagger := \sup_{\tau,w}|\widetilde{\mathbf{Z}}(w)^\top \boldsymbol{\beta}_n^\dagger(\tau) - h(w;\tau)|$ and assume that
$$\xi_m c_n^\dagger = o(1); \tag{3.7}$$
$$\sup_{\tau \in \mathcal{T}} \mathbb{E}[f_{Y|X}(Q(X;\tau)|X)\|h_{VW}(W;\tau) - A(\tau)\widetilde{\mathbf{Z}}(W))\|^2] = O(\lambda_n^2) \text{ with } \xi_m \lambda_n^2 = o(1); \tag{3.8}$$

(C2) We have $\max_{j \leq k} |V_j| < C$ almost surely for some constant $C > 0$.

Bounds on $c_n^\dagger$ can be obtained under various assumptions on the basis expansion and smoothness of the function $w \mapsto h(w;\tau)$. Assume for instance that $\mathcal{W} = [0,1]^{k'}$, that $h(\cdot;\cdot) \in \Lambda_c^\eta(\mathcal{W},\mathcal{T})$ and that $\widetilde{\mathbf{Z}}$ corresponds to a tensor product B-spline basis of order $q$ on $\mathcal{W}$ with $m^{1/k'}$ equidistant knots in each coordinate. Assuming that $(V,W)$ has a density $f_{V,W}$ such that $0 < \inf_{v,w} f_{V,W}(v,w) \leq \sup_{v,w} f_{V,W}(v,w) < \infty$ and $q > \eta$, we show in Remark B.1 that $c_n^\dagger = O(m^{-\lfloor\eta\rfloor/k'})$.

Assumption (3.8) essentially states that $h_{VW}$ can be approximated by a series estimator sufficiently well. This assumption is necessary to ensure that $\boldsymbol{\alpha}(\tau)$ is estimable at a parametric rate without under-smoothing when estimating $h(\cdot;\tau)$. In general, (3.8) is a non-trivial high-level assumption. It can be verified under smoothness conditions on the joint density of $(X,Y)$ by applying arguments similar to those in Appendix S.1 of Cheng et al. (2014).

In addition to (C1)-(C2), we need the following condition.

(B1') Assume that
$$\left(\frac{m\xi_m^{2/3}\log n}{n}\right)^{3/4} + c_n^{\dagger 2}\xi_m = o(n^{-1/2}).$$

Moreover, assume that $c_n^\dagger \lambda_n = o(n^{-1/2})$ and $mc_n^\dagger \log n = o(1)$.

We now are ready to state the main result of this section.

**Theorem 3.1.** *Let Conditions (A1)-(A3) hold with $\mathbf{Z} = (V^\top, \widetilde{\mathbf{Z}}(W)^\top)^\top$, (B1') and (C1)-(C2) hold for $\boldsymbol{\beta}_n^\dagger(\tau)$ defined in (3.3). For any sequence $\mathbf{w}_n \in \mathbb{R}^m$ with $\mathbb{E}\big[|\mathbf{w}_n^\top M_2(\tau_2)^{-1}\widetilde{\mathbf{Z}}(W)|\big] = o(\|\mathbf{w}_n\|)$ where $M_2(\tau) := \mathbb{E}\big[\widetilde{\mathbf{Z}}(W)\widetilde{\mathbf{Z}}(W)^\top f_{Y|X}(Q(X;\tau)|X)\big]$, if*
$$\Gamma_{22}(\tau_1,\tau_2) = \lim_{n\to\infty} \|\mathbf{w}_n\|^{-2}\mathbf{w}_n^\top M_2(\tau_1)^{-1}\mathbb{E}[\widetilde{\mathbf{Z}}(W)\widetilde{\mathbf{Z}}(W)^\top]M_2(\tau_2)^{-1}\mathbf{w}_n \tag{3.9}$$

*exists, then*
$$\begin{pmatrix} \sqrt{n}(\widehat{\boldsymbol{\alpha}}(\cdot) - \boldsymbol{\alpha}(\cdot)) \\ \frac{\sqrt{n}}{\|\mathbf{w}_n\|}\mathbf{w}_n^\top(\widehat{\boldsymbol{\beta}}^\dagger(\cdot) - \boldsymbol{\beta}_n^\dagger(\cdot)) \end{pmatrix} \rightsquigarrow (\mathbb{G}_1(\cdot),...,\mathbb{G}_k(\cdot),\mathbb{G}_h(\cdot))^\top \text{ in } (\ell^\infty(\mathcal{T}))^{k+1}, \tag{3.10}$$

*and the multivariate process $(\mathbb{G}_1(\cdot),...,\mathbb{G}_k(\cdot),\mathbb{G}_h(\cdot))$ has the covariance function*
$$\Gamma(\tau_1,\tau_2;\mathbf{w}_n) = (\tau_1 \wedge \tau_2 - \tau_1\tau_2)\begin{pmatrix} \Gamma_{11}(\tau_1,\tau_2) & \mathbf{0}_k \\ \mathbf{0}_k^\top & \Gamma_{22}(\tau_1,\tau_2) \end{pmatrix} \tag{3.11}$$

*with*

$$\Gamma_{11}(\tau_1, \tau_2) = M_{1,h}(\tau_1)^{-1} \mathbb{E}\big[(V - h_{VW}(W; \tau_1))(V - h_{VW}(W; \tau_2))^\top\big] M_{1,h}(\tau_2)^{-1} \quad (3.12)$$

where $M_{1,h}(\tau) = \mathbb{E}\big[(V - h_{VW}(W; \tau))(V - h_{VW}(W; \tau))^\top f_{Y|X}(Q(X; \tau)|X)\big]$. In addition, at any fixed $w_0 \in \mathbb{R}^{k'}$, let $\mathbf{w}_n = \widetilde{\mathbf{Z}}(w_0)$ satisfy the above conditions, $\widehat{h}(w_0; \tau) = \widetilde{\mathbf{Z}}(w_0)^\top \widehat{\boldsymbol{\beta}}^\dagger(\tau)$, $c_n^\dagger = o(\|\widetilde{\mathbf{Z}}(w_0)\| n^{-1/2})$, then

$$\begin{pmatrix} \sqrt{n}\{\widehat{\boldsymbol{\alpha}}(\cdot) - \boldsymbol{\alpha}(\cdot)\} \\ \frac{\sqrt{n}}{\|\widetilde{\mathbf{Z}}(w_0)\|}\{\widehat{h}(w_0; \cdot) - h(w_0; \cdot)\} \end{pmatrix} \rightsquigarrow \big(\mathbb{G}_1(\cdot), ..., \mathbb{G}_k(\cdot), \mathbb{G}_h(w_0; \cdot)\big)^\top \text{ in } (\ell^\infty(\mathcal{T}))^{k+1}, \quad (3.13)$$

where $(\mathbb{G}_1(\cdot), ..., \mathbb{G}_k(\cdot), \mathbb{G}_h(w_0; \cdot))$ are centered Gaussian processes with joint covariance function $\Gamma_{w_0}(\tau_1, \tau_2)$ of the form (3.11) where $\Gamma_{22}(\tau_1, \tau_2)$ is defined through the limit in (3.9) with $\mathbf{w}_n$ replaced by $\widetilde{\mathbf{Z}}(w_0)$. In particular, there exists a version of $\mathbb{G}_h(w_0; \cdot)$ with almost surely continuous sample paths.

The proof of Theorem 3.1 is presented in Section A.3. The invertibility of the matrices $M_{1,h}(\tau)$ and $M_2(\tau)$ is discussed in Remark 5.5. In general, $\widehat{\boldsymbol{\alpha}}(\tau)$ is not semiparametric efficient, as its covariance matrix $\tau(1-\tau)\Gamma_{11}$ does not achieve the efficiency bound given in Section 5 of Lee (2003).

The joint asymptotic process convergence result (in $\ell^\infty(\mathcal{T})$) presented in Theorem 3.1 is new in the quantile regression literature. The block structure of covariance function $\Gamma$ defined in (3.11) implies that $\widehat{\boldsymbol{\alpha}}(\tau)$ and $\widehat{h}(w_0; \tau)$ are asymptotically independent for any fixed $\tau$. This effect was recently discovered by Cheng and Shang (2015) in the case of mean regression, named as *joint asymptotics phenomenon*.

**Remark 3.2.** We point out that $\mathbb{E}\big[|\mathbf{w}_n^\top M_2(\tau_2)^{-1} \widetilde{\mathbf{Z}}(W)|\big] = o(\|\mathbf{w}_n\|)$ is a crucial sufficient condition for asymptotic independence between the parametric and nonparametric parts. We conjecture that this condition is also necessary. This condition holds, for example, for $\mathbf{w}_n = \widetilde{\mathbf{Z}}(w_0)$ or $\mathbf{w}_n = \partial_{w_j} \widetilde{\mathbf{Z}}(w_0)$ at a fixed $w_0$, $j = 1, ..., k'$, where $\widetilde{\mathbf{Z}}(w)$ is a vector of $B$-spline basis. However, this condition may not hold for other estimators. Consider for instance the case $\mathcal{W} = [0, 1]$, $B$-splines of order zero $\widetilde{\mathbf{Z}}$ and the vector $\mathbf{w}_n = \int_0^\lambda \widetilde{\mathbf{Z}}(w) dw$ for some $\lambda > 0$. In this case $\|\mathbf{w}_n\| \asymp 1$, and one can show that $\mathbb{E}\big[|\mathbf{w}_n^\top M_2(\tau_2)^{-1} \widetilde{\mathbf{Z}}(W)|\big] \asymp 1$ instead. A more detailed investigation of related questions is left to future research.

**Remark 3.3.** A seemingly more natural choice for the centering vector, which was also considered in Belloni et al. (2016), is

$$\boldsymbol{\gamma}_n^*(\tau) = (\boldsymbol{\alpha}_n^*(\tau), \boldsymbol{\beta}_n^*(\tau)) := \arg\min_{(\mathbf{a}, \mathbf{b})} \mathbb{E}[\rho_\tau(Y - \mathbf{a}^\top V - \mathbf{b}^\top \widetilde{\mathbf{Z}}(W))], \quad (3.14)$$

which gives $g_n(\boldsymbol{\gamma}_n^*(\tau)) = 0$. However, a major drawback of centering with $\boldsymbol{\gamma}_n^*(\tau)$ is that it is impossible to find a good bound for the difference $\boldsymbol{\alpha}_n^*(\tau) - \boldsymbol{\alpha}(\tau)$ without imposing restrictive assumptions. However, such a bound is needed to show that the bias of $\boldsymbol{\alpha}_n^*(\tau)$ is of order $o(n^{-1/2})$ which is required establish (3.10) in Theorem 3.1.

## 4. Applications of Weak Convergence Results

In this section, we consider applications of the process convergence results to the estimation of conditional distribution functions and non-crossing quantile curves via rearrangement operators. For the former estimation, define the functional (see Dette and Volgushev (2008), Chernozhukov et al. (2010) or Volgushev (2013) for similar ideas)

$$\Phi : \begin{cases} \ell^\infty((\tau_L, \tau_U)) & \to \ \ell^\infty(\mathbb{R}) \\ \Phi(f)(y) & := \ \tau_L + \int_{\tau_L}^{\tau_U} \mathbf{1}\{f(\tau) < y\} d\tau. \end{cases}$$

A simple calculation shows that

$$\Phi(Q(x;\cdot))(y) = \begin{cases} \tau_L & \text{if } F_{Y|X}(y|x) < \tau_L \\ F_{Y|X}(y|x) & \text{if } \tau_L \leq F_{Y|X}(y|x) \leq \tau_U \\ \tau_U & \text{if } F_{Y|X}(y|x) > \tau_U. \end{cases}$$

The latter identity motivates the following estimator of the conditional distribution function

$$\widehat{F}_{Y|X}(y|x) := \tau_L + \int_{\tau_L}^{\tau_U} \mathbf{1}\{\widehat{Q}(x;\tau) < y\} d\tau,$$

where $\widehat{Q}(x;\tau)$ denotes the estimator of the conditional quantile function in any of the three settings discussed in Sections 2 and 3. By following the arguments in Chernozhukov et al. (2010), one can easily show that under suitable assumptions the functional $\Phi$ is compactly differentiable (see Section A.5 for more details). Hence, the general process convergence results in Sections 2 and 3 allow to easily establish the asymptotic properties of $\widehat{F}_{Y|X}$ - see Corollary 4.1 at the end of this section.

The second functional of interest is the monotone rearrangement operator, defined as follows

$$\Psi : \begin{cases} \ell^\infty((\tau_L, \tau_U)) \to \ell^\infty((\tau_L, \tau_U)) \\ \Psi(f)(\tau) := \inf\{y : \Phi(Q(x;\cdot))(y) \geq \tau\}. \end{cases}$$

The main motivation for considering $\Psi$ is that the function $\tau \mapsto \Psi(f)(\tau)$ is by construction non-decreasing. Thus for any initial estimator $\widehat{Q}(x;\cdot)$, its rearranged version $\Psi(\widehat{Q}(x;\cdot))(\tau)$ is an estimator of the conditional quantile function which avoids the issue of quantile crossing. For more detailed discussions of rearrangement operators and their use in avoiding quantile crossing we refer the interested reader to Dette and Volgushev (2008) and Chernozhukov et al. (2010).

**Corollary 4.1** (Convergence of $\widehat{F}(y|x)$ and $\Psi(\widehat{Q}(x;\tau))$). *For any fixed $x_0$ and an initial estimator $\widehat{Q}(x_0, \cdot)$, we have for any compact sets $[\tau_U, \tau_L] \subset \mathcal{T}, \mathcal{Y} \subset \mathcal{Y}_{0,\mathcal{T}} := \{y : F_{Y|X}(y|x_0) \in \mathcal{T}\}$*

$$a_n\{\widehat{F}_{Y|X}(\cdot|x_0) - F_{Y|X}(\cdot|x_0)\} \rightsquigarrow -f_{Y|X}(\cdot|x_0)\mathbb{G}(x_0; F_{Y|X}(\cdot|x_0)) \text{ in } \ell^\infty(\mathcal{Y}),$$
$$a_n\{\Psi(\widehat{Q}(x_0;\cdot))(\cdot) - Q(x_0;\cdot)\} \rightsquigarrow \mathbb{G}(x_0;\cdot) \text{ in } \ell^\infty((\tau_U, \tau_L)),$$

*where $\widehat{Q}(x_0, \cdot)$, the normalization $a_n$, and the process $\mathbb{G}(x_0;\cdot)$ are stated as follows*

1. *(Linear model with increasing dimension) Suppose $\mathbf{Z}(X) = X$, $\widehat{Q}(x_0, \cdot) = \widehat{\boldsymbol{\gamma}}(\cdot)^\top x_0$ and the conditions in Corollary 2.2 hold. In this case, we have $a_n = \sqrt{n}/\|x_0\|$. $\mathbb{G}(x_0;\cdot)$ is a centered Gaussian process with covariance function $H_1(\tau_1, \tau_2; x_0)$ defined in (2.5).*

2. *(Nonparametric model) Suppose $\widehat{Q}(x_0, \cdot) = \boldsymbol{\beta}(\cdot)^\top \mathbf{B}(x_0)$ and the conditions in Theorems 2.4 hold. In this case, we have $a_n = \sqrt{n}/\|\mathbf{B}(x_0)\|$. $\mathbb{G}(x_0;\cdot)$ is a centered Gaussian process with covariance function $\widetilde{H}(\tau_1, \tau_2; \mathbf{B}(x_0))$ defined in (2.8).*

3. *(Partial linear model) Suppose $x_0^\top = (v_0^\top, w_0^\top)$, $\widehat{Q}(x_0, \cdot) = \widehat{\boldsymbol{\gamma}}^\dagger(\tau)^\top (v^\top, \widetilde{\mathbf{Z}}(w_0)^\top)^\top$ and the conditions in Theorem 3.1 hold. In this case, we have $a_n = \sqrt{n}/\|\widetilde{\mathbf{Z}}(w_0)\|$. $\mathbb{G}(x_0;\cdot)$ is a centered Gaussian process with covariance function $\Gamma_{22}(\tau_1, \tau_2; \widetilde{\mathbf{Z}}(w_0))$ defined in (3.9).*

The proof of Corollary 4.1 is a direct consequence of the functional delta method, combined with the process convergence results established in Section 2 and Section 3 and Hadamard differentiability results of certain functionals established in Chernozhukov et al. (2010). Details can be found in Section A.5.

**Remark 4.2** (More Statistical Applications of Corollary 4.1). In practice, many quantities of interest can be written as Hadamard differentiable functionals of distribution functions such as some $M$ and $L$ estimators ([Fernholz, 1983](#), Chapter 7), conditional distributions ([Bassett and Koenker, 1982](#)), the Wilcoxon test statistics ([van der Vaart and Wellner, 1996](#), Example 3.9.18) and Gini indices ([Barrett and Donald, 2009](#)). Based on the chain rule of Hadamard derivative, Corollary [4.1](#) can be applied to prove the asymptotic normality of these statistical estimators. Moreover, detection of treatment effect on the conditional distribution after an intervention ([Koenker and Xiao, 2002](#); [Qu and Yoon, 2015](#)) is often based on Kolmogorov-Smirnov or Cramér-von Mises statistics, whose asymptotic distribution can also be found by applying Corollary [4.1](#).

## 5. Bahadur representations

In this section, we provide Bahadur representations for the estimators discussed in Sections [2](#) and [3](#). In Sections [5.1](#) and [5.2](#), we state Bahadur representations for general series estimators and a more specific choice of local basis function, respectively. In particular, the latter representation is developed with an improved remainder term. Section [5.3](#) contains a special case of the general theorem in Section [5.1](#) that is particularly tailored to partial linear models. The remainders in these representations are shown to have exponential tail probabilities (uniformly over $\mathcal{T}$).

### 5.1. A Fundamental Bahadur Representation

Our first result gives a Bahadur representation for $\widehat{\gamma}(\tau) - \gamma_n(\tau)$ for centering functions $\gamma_n$ satisfying certain conditions. Recall the definition of $\widehat{\gamma}(\tau)$ in [(2.1)](#). This kind of representation for quantile regression with an increasing number of covariates has previously been established in Theorem 2 of [Belloni et al. (2016)](#). Compared to their results, the Bahadur representation given below has several advantages. First, we allow for a more general centering. This is helpful for the analysis of partial linear models (see Sections [3](#) and [S.1.2](#)). Second, we provide exponential tail bounds on remainder terms, which is much more explicit and sharper than those in [Belloni et al. (2016)](#).

**Theorem 5.1.** *Suppose Conditions [(A1)-(A3)](#) hold and that additionally $m\xi_m^2 \log n = o(n)$. Then, for any $\gamma_n(\cdot)$ satisfying $g_n(\gamma_n) = o(\xi_m^{-1})$ and $c_n(\gamma_n) = o(1)$, we have*

$$\widehat{\gamma}(\tau) - \gamma_n(\tau) = -\frac{1}{n} J_m(\tau)^{-1} \sum_{i=1}^{n} \psi(Y_i, \mathbf{Z}_i; \gamma_n(\tau), \tau) + r_{n,1}(\tau) + r_{n,2}(\tau) + r_{n,3}(\tau) + r_{n,4}(\tau).$$

*The remainder terms $r_{n,j}$'s can be bounded as follows:*

$$\sup_{\tau \in \mathcal{T}} \|r_{n,1}(\tau)\| \leq \frac{1}{\inf_{\tau \in \mathcal{T}} \lambda_{\min}(J_m(\tau))} \frac{m\xi_m}{n} \quad a.s. \tag{5.1}$$

*Moreover, we have for any $\kappa_n \ll n/\xi_m^2$, sufficiently large $n$, and a constant $C$ independent of $n$*

$$P\Big(\sup_{\tau \in \mathcal{T}} \|r_{n,j}(\tau)\| \leq C \Re_j(\kappa_n)\Big) \geq 1 - 2e^{-\kappa_n}, \quad j = 2, 3, 4,$$

*where*

$$\Re_2(\kappa_n) := \xi_m \Big(\Big(\frac{m}{n} \log n\Big)^{1/2} + \Big(\frac{\kappa_n}{n}\Big)^{1/2} + g_n\Big)^2, \tag{5.2}$$

$$\Re_3(\kappa_n) := \Big(\Big(\frac{m \log n}{n}\Big)^{1/2} + \Big(\frac{\kappa_n}{n}\Big)^{1/2} + g_n\Big)^{1/2} \Big(\Big(\frac{m\xi_m \log n}{n}\Big)^{1/2} + \Big(\frac{\xi_m \kappa_n}{n}\Big)^{1/2}\Big), \tag{5.3}$$

$$\Re_4(\kappa_n) := c_n \Big(\Big(\frac{m}{n} \log n\Big)^{1/2} + \Big(\frac{\kappa_n}{n}\Big)^{1/2}\Big) + g_n. \tag{5.4}$$

The proof for Theorem [5.1](#) can be found in Section [S.1.1](#).

## 5.2. Bahadur Representation for Local Basis Series Estimator

In this section, we focus on basis expansions $\mathbf{B}$ satisfying (L) and derive a Bahadur representation for linear functionals of the form $\mathbf{u}_n^\top(\widehat{\boldsymbol{\beta}}(\tau) - \boldsymbol{\beta}_n(\tau))$, where the vector $\mathbf{u}_n$ can have at most a finite number of consecutive non-zero entries. Such linear functionals are of interest since the estimator of the quantile function itself as well as estimators of derivatives can be represented in exactly this form - see Remark 5.3 for additional details. The advantage of concentrating on vectors with this particular structure is that we can substantially improve the rates of remainder terms compared to the general setting in Theorem 5.1.

**Theorem 5.2.** *Suppose Conditions (A1)-(A3) and (L) hold with $\mathbf{Z}(x) = \mathbf{B}(x)$. Assume additionally that $m\xi_m^2(\log n)^2 = o(n)$ and that $\widetilde{c}_n = o(1)$ and that $\mathcal{I} \subset \{1, ..., m\}$ consists of at most $L$ consecutive integers. Then, for $\boldsymbol{\beta}_n(\tau)$ defined as (2.7) and $\mathbf{u}_n \in \mathcal{S}_\mathcal{I}^{m-1}$ we have*

$$\mathbf{u}_n^\top(\widehat{\boldsymbol{\beta}}(\tau) - \boldsymbol{\beta}_n(\tau)) = -\mathbf{u}_n^\top \widetilde{J}_m(\tau)^{-1} n^{-1} \sum_{i=1}^n \mathbf{B}_i (\mathbf{1}\{Y_i \leq \boldsymbol{\beta}_n(\tau)^\top \mathbf{B}_i\} - \tau) + \sum_{k=1}^4 r_{n,k}(\tau, \mathbf{u}_n), \tag{5.5}$$

*where the remainder terms $r_{n,j}$'s can be bounded as follows:*

$$\sup_{\mathbf{u}_n \in \mathcal{S}_\mathcal{I}^{m-1}} \sup_{\tau \in \mathcal{T}} |r_{n,1}(\tau, \mathbf{u}_n)| \lesssim \frac{\xi_m \log n}{n} \quad a.s. \tag{5.6}$$

$$\sup_{\mathbf{u}_n \in \mathcal{S}_\mathcal{I}^{m-1}} \sup_{\tau \in \mathcal{T}} |r_{n,4}(\tau, \mathbf{u}_n)| \leq \frac{1}{n} + \frac{1}{2}\overline{f'}\widetilde{c}_n^2 \sup_{\mathbf{u}_n \in \mathcal{S}_\mathcal{I}^{m-1}} \widetilde{\mathcal{E}}(\mathbf{u}_n, \mathbf{B}) \quad a.s. \tag{5.7}$$

*where $\widetilde{\mathcal{E}}(\mathbf{u}_n, \mathbf{B}) := \sup_\tau \mathbb{E}|\mathbf{u}_n \widetilde{J}_m(\tau)^{-1} \mathbf{B}|$. Moreover, we have for any $\kappa_n \ll n/\xi_m^2$, all sufficiently large $n$, and a constant $C$ independent of $n$*

$$P\Big(\sup_{\mathbf{u}_n \in \mathcal{S}_\mathcal{I}^{m-1}} \sup_{\tau \in \mathcal{T}} |r_{n,j}(\tau, \mathbf{u}_n)| \leq C\widetilde{\mathfrak{R}}_j(\kappa_n)\Big) \geq 1 - n^2 e^{-\kappa_n}, \quad j = 2, 3$$

*where*

$$\widetilde{\mathfrak{R}}_2(\kappa_n) := C \sup_{\mathbf{u}_n \in \mathcal{S}_\mathcal{I}^{m-1}} \widetilde{\mathcal{E}}(\mathbf{u}_n, \mathbf{B}) \Big(\frac{\xi_m(\log n + \kappa_n^{1/2})}{n^{1/2}} + \widetilde{c}_n^2\Big)^2, \tag{5.8}$$

$$\widetilde{\mathfrak{R}}_3(\kappa_n) := C\Big(\widetilde{c}_n \frac{\kappa_n^{1/2} \vee \log n}{n^{1/2}} + \frac{\xi_m^{1/2}(\kappa_n^{1/2} \vee \log n)^{3/2}}{n^{3/4}}\Big). \tag{5.9}$$

Theorem 5.2 is proved in Section S.1.2. We note that by Hölder's inequality and assumptions (A1)-(A3), we have a simple bound for

$$\sup_{\mathbf{u}_n \in \mathcal{S}_\mathcal{I}^{m-1}} \widetilde{\mathcal{E}}(\mathbf{u}_n, \mathbf{B}) \leq \sup_{\mathbf{u}_n \in \mathcal{S}^{m-1}} \sup_\tau \big(\mathbf{u}_n^\top J_m^{-1}(\tau) \mathbb{E}[\mathbf{B}\mathbf{B}^\top] J_m^{-1}(\tau) \mathbf{u}_n\big)^{1/2} = O(1).$$

**Remark 5.3.** Theorem 5.2 enables us to study several quantities associated with the quantile function $Q(x; \tau)$. For instance, consider the spline setting of Example 2.3. Setting $\mathbf{u}_n = \mathbf{B}(x)/\|\mathbf{B}(x)\|$ in the Theorem 5.2 yields a representation for $\widehat{Q}(x; \tau)$, while setting $\mathbf{u}_n = \mathbf{B}'(x)/\|\mathbf{B}'(x)\|$ yields a representation for the estimator of the derivative $\partial_x Q(x; \tau)$. Uniformity in $x$ follows once we observe that for different values of $x$, the support of the vector $\mathbf{B}(x)$ is always consecutive so that there is at most $n^l$, $l > 0$, number of different sets $\mathcal{I}$ that we need to consider.

## 5.3. Bahadur Representation for Partial Linear Models

In this section, we provide a joint Bahadur representation for the parametric and non-parametric part of this model. Recall the partial linear model $Q(X; \tau) = h(W; \tau) + \boldsymbol{\alpha}(\tau)^\top V$.

**Theorem 5.4.** *Let conditions (A1)-(A3), (C1)-(C2) hold with $\mathbf{Z} = (V^\top, \widetilde{\mathbf{Z}}(W)^\top)^\top$ and assume $m\xi_m^2 (\log n)^2 = o(n)$. Then*

$$\begin{pmatrix} \widehat{\boldsymbol{\alpha}}(\tau) - \boldsymbol{\alpha}(\tau) \\ \widehat{\boldsymbol{\beta}}^\dagger(\tau) - \boldsymbol{\beta}_n^\dagger(\tau) \end{pmatrix} = -J_m(\tau)^{-1} n^{-1} \sum_{i=1}^n \mathbf{Z}_i (\mathbf{1}\{Y_i \leq \{\boldsymbol{\gamma}_n(\tau)^\dagger\}^\top \mathbf{Z}_i\} - \tau) + \sum_{j=1}^4 r_{n,j}(\tau),$$

*where the remainder terms $r_{n,j}$'s satisfy the bounds stated in Theorem 5.1 with $g_n = \xi_m c_n^{\dagger 2}$. Additionally, the matrix $J_m^{-1}(\tau)$ can be represented as*

$$J_m^{-1}(\tau) = \begin{pmatrix} M_1(\tau)^{-1} & -M_1(\tau)^{-1} A(\tau) \\ -A(\tau)^\top M_1(\tau)^{-1} & M_2(\tau)^{-1} + A(\tau)^\top M_1(\tau)^{-1} A(\tau) \end{pmatrix} \quad (5.10)$$

*where*

$$M_1(\tau) := \mathbb{E}[f_{Y|X}(Q(X;\tau)|X)(V - A(\tau)\widetilde{\mathbf{Z}}(W))(V - A(\tau)\widetilde{\mathbf{Z}}(W))^\top],$$
$$M_2(\tau) = \mathbb{E}\big[\widetilde{\mathbf{Z}}(W)\widetilde{\mathbf{Z}}(W)^\top f_{Y|X}(Q(X;\tau)|X)\big],$$

*and $A(\tau)$ is defined in (3.6).*

See Section S.1.3 for the proof of Theorem 5.4.

**Remark 5.5.** We discuss the positive definiteness of $M_1(\tau)$ and $M_2(\tau)$. Following Condition (A1) with $\mathbf{Z} = (V^\top, \widetilde{\mathbf{Z}}(W)^\top)^\top$, we have

$$1/M \leq \inf_{\tau \in \mathcal{T}} \lambda_{\min}(M_1(\tau)) \leq \sup_{\tau \in \mathcal{T}} \lambda_{\max}(M_1(\tau)) \leq M; \quad (5.11)$$

$$1/M \leq \inf_{\tau \in \mathcal{T}} \lambda_{\min}(M_2(\tau)) \leq \sup_{\tau \in \mathcal{T}} \lambda_{\max}(M_2(\tau)) \leq M, \quad (5.12)$$

for all $n$. To see this, observe that $M_1(\tau) = [I_k| - A(\tau)] J_m(\tau) [I_k| - A(\tau)]^\top$ where $I_k$ is the $k$-dimensional identity matrix, $[A|B]$ denotes the block matrix with $A$ in the left block and $B$ in the right block, and

$$J_m(\tau) = \begin{pmatrix} M_1(\tau) + A(\tau) M_2(\tau) A(\tau)^\top & A(\tau) M_2(\tau) \\ M_2(\tau) A(\tau)^\top & M_2(\tau) \end{pmatrix},$$

whose form follows from the definition and the condition (3.5) (see the proof for Theorem 5.4 for more details). Thus, for an arbitrary nonzero vector $\mathbf{a} \in \mathbb{R}^k$, by Condition (A1),

$$0 < 1/M \leq \mathbf{a}^\top M_1(\tau) \mathbf{a} = \mathbf{a}^\top [I_k| - A(\tau)] J_m(\tau) [I_k| - A(\tau)]^\top \mathbf{a} \leq M < \infty$$

by the strictly positive definiteness of $J_m(\tau)$ for some $M > 0$.

The strictly positive definiteness of $M_2(\tau)$ follows directly from the observation that

$$0 < 1/M \leq \mathbf{b}^\top M_2(\tau) \mathbf{b} = (\mathbf{0}_k^\top, \mathbf{b}^\top) J_m(\tau) (\mathbf{0}_k^\top, \mathbf{b}^\top)^\top \leq M < \infty,$$

for all nonzero $\mathbf{b} \in \mathbb{R}^m$ and some $M > 0$.

# APPENDIX

This appendix gives technical details of the results shown in the main text. Appendix A contains all the proofs for weak convergence results in Theorems 2.1, 2.4 and 3.1. Appendix B discusses basis approximation errors with full technical details.

**Additional Notations.** Define for a function $x \mapsto f(x)$ that $\mathbb{G}_n(f) := n^{1/2} \int f(x)(d\mathbb{P}_n(x) - dP(x))$ and $\|f\|_{L_p(P)} = (\int |f(x)|^p dP(x))^{1/p}$ for $0 < p < \infty$. For a class of functions $\mathcal{G}$, let $\|\mathbb{P}_n - P\|_\mathcal{G} := \sup_{f \in \mathcal{G}} |\mathbb{P}_n f - P f|$. For any $\epsilon > 0$, the covering number $N(\epsilon, \mathcal{G}, L_p)$ is the minimal number of balls of radius $\epsilon$ (under $L_p$-norm) that is needed to cover $\mathcal{G}$. The bracketing number $N_{[\,]}(\epsilon, \mathcal{G}, L_p)$ is the minimal number of $\epsilon$-brackets that is needed to cover $\mathcal{G}$. An $\epsilon$-bracket refers to a pair of functions within an $\epsilon$ distance: $\|u - l\|_{L_p} < \epsilon$. Throughout the proofs, $C, C_1, C_2$ etc. will denote constants which do not depend on $n$ but may have different values in different lines.

# APPENDIX A: Proofs for Process Convergence

## A.1. Proof of Theorem 2.1

### A.1.1. Proof of (2.2)

Under conditions (A1)-(A3) and those in Theorem 2.1 it follows from Theorem 5.1 applied with $\kappa_n = c \log n$ for a suitable constant $c$ (note that the conditions $g_n = o(\xi_m^{-1})$ and $c_n = o(1)$ in Theorem 5.1 follow under the assumptions of Theorem 2.1) that

$$\mathbf{u}_n^\top (\widehat{\boldsymbol{\gamma}}(\tau) - \boldsymbol{\gamma}(\tau)) + \frac{1}{n} \mathbf{u}_n^\top J_m(\tau)^{-1} \sum_{i=1}^n \mathbf{Z}_i \big(\mathbf{1}\{Y_i \leq Q(X_i; \tau)\} - \tau\big) = I(\tau) + o_P(\|\mathbf{u}_n\| n^{-1/2}),$$

where the remainder term is uniform in $\mathcal{T}$ and

$$I(\tau) := -n^{-1} \mathbf{u}_n^\top J_m(\tau)^{-1} \sum_{i=1}^n \mathbf{Z}_i \big(\mathbf{1}\{Y_i \leq \mathbf{Z}_i^\top \boldsymbol{\gamma}_n(\tau)\} - \mathbf{1}\{Y_i \leq Q(X_i; \tau)\}\big).$$

Under the assumption $\chi_{\boldsymbol{\gamma}_n}(\mathbf{u}_n, \mathbf{Z}) = o(\|\mathbf{u}_n\| n^{-1/2})$, we have $\sup_{\tau \in \mathcal{T}} |\mathbb{E}[I(\tau)]| = o(\|\mathbf{u}_n\| n^{-1/2})$ and moreover

$$\sup_{\tau \in \mathcal{T}} |I(\tau) - \mathbb{E}[I(\tau)]| \leq \|\mathbf{u}_n\| [\inf_{\tau \in \mathcal{T}} \lambda_{\min}(J_m(\tau))]^{-1} \|\mathbb{P}_n - P\|_{\mathcal{G}_5},$$

where the class of functions $\mathcal{G}_5$ is defined as

$$\mathcal{G}_5(\mathbf{Z}, \boldsymbol{\gamma}_n) := \big\{(X, Y) \mapsto \mathbf{a}^\top \mathbf{Z}(X) \mathbf{1}\{\|\mathbf{Z}(X)\| \leq \xi_m\} \big(\mathbf{1}\{Y \leq \mathbf{Z}(X)^\top \boldsymbol{\gamma}_n(\tau)\} - \mathbf{1}\{Y \leq Q(X; \tau)\}\big) \\ \big| \tau \in \mathcal{T}, \mathbf{a} \in \mathcal{S}^{m-1}\big\}.$$

It remains to bound $\|\mathbb{P}_n - P\|_{\mathcal{G}_5}$. For any $f \in \mathcal{G}_5$ and a sufficiently large $C$, we obtain

$$|f| \leq |\mathbf{a}^\top \mathbf{Z}| \leq \xi_m,$$
$$\|f\|_{L_2(P)}^2 \leq 2\bar{f} c_n \lambda_{\max}(\mathbb{E}[\mathbf{Z}\mathbf{Z}^\top]) \leq C c_n.$$

By Lemma 21 of Belloni et al. (2016), the VC index of $\mathcal{G}_5$ is of the order $O(m)$. Therefore, we obtain from (S.2.2)

$$\mathbb{E}\|\mathbb{P}_n - P\|_{\mathcal{G}_5} \leq \widetilde{C} \Big[\Big(\frac{m c_n}{n} \log \frac{\xi_m}{\sqrt{c_n}}\Big)^{1/2} + \frac{m \xi_m}{n} \log \frac{\xi_m}{\sqrt{c_n}}\Big]. \tag{A.1}$$

For any $\kappa_n > 0$, let
$$r'_{N,3}(\kappa_n) = C\left[\left(\frac{mc_n}{n}\log\frac{\xi_m}{\sqrt{c_n}}\right)^{1/2} + \frac{m\xi_m}{n}\log\frac{\xi_m}{\sqrt{c_n}} + \left(\frac{c_n}{n}\kappa_n\right)^{1/2} + \frac{\xi_m}{n}\kappa_n\right]$$
for a sufficiently large constant $C > 0$. We obtain from (S.2.3) combined with (A.1)
$$P\left\{\sup_{\tau\in\mathcal{T}}|I(\tau)| \geq \|\mathbf{u}_n\|r'_{N,3}(\kappa_n) + \sup_{\tau\in\mathcal{T}}\left|\mathbb{E}[I(\tau)]\right|\right\} \leq e^{-\kappa_n}.$$
Finally, note that under condition $mc_n \log n = o(1)$ and
$$m^3\xi_m^2(\log n)^3 = o(n)$$
we have that $r'_{N,3}(\log n) = o(n^{-1/2})$. This completes the proof of (2.2). $\square$

### A.1.2. Proof of (2.4)

Throughout this subsection assume without loss of generality that $\|\mathbf{u}_n\| = 1$. It suffices to prove finite dimensional convergence and asymptotic equicontinuity. Asymptotic equicontinuity follows from (A.33). The existence of a version of the limit with continuous sample paths is a consequence of Theorem 1.5.7 and Addendum 1.5.8 in van der Vaart and Wellner (1996).

So, we only need to focus on finite dimensional convergence.

Let
$$\mathbb{G}_n(\tau) := \frac{1}{\sqrt{n}}\mathbf{u}_n^\top J_m(\tau)^{-1}\sum_{i=1}^n \mathbf{Z}_i(\mathbf{1}\{Y_i \leq Q(X_i;\tau)\} - \tau),$$
and $\mathbb{G}$ be the Gaussian process defined in (2.4). From Cramér-Wold theorem, the goal is to show for arbitrary set of $\{\tau_1,...,\tau_L\}$ and $\{\lambda_1,...,\lambda_L\} \in \mathbb{R}^L$, we have
$$\sum_{l=1}^L \lambda_l \mathbb{G}_n(\tau_l) \xrightarrow{d} \sum_{l=1}^L \lambda_l \mathbb{G}(\tau_l).$$

Let the triangular array $V_{n,i}(\tau) := n^{-1/2}\mathbf{u}_n^\top J_m(\tau)^{-1}\mathbf{Z}_i(\mathbf{1}\{Y_i \leq Q(X_i;\tau)\} - \tau)$. Then for all $\tau \in \mathcal{T}$, we have $\mathbb{E}[V_{n,i}(\tau)] = 0$, $|V_{n,i}| \leq n^{-1/2}\xi_m$ and $\text{var}(V_{n,i}(\tau)) = n^{-1}\mathbf{u}_n^\top J_m(\tau)^{-1}\mathbb{E}[\mathbf{Z}_i\mathbf{Z}_i^\top]J_m(\tau)^{-1}\mathbf{u}_n\tau(1-\tau) < \infty$ by Conditions (A1)-(A3). We can express $\mathbb{G}_n(\tau) = \sum_{i=1}^n V_{n,i}(\tau)$ and $\sum_{l=1}^L \lambda_l \mathbb{G}_n(\tau_l) = \sum_{i=1}^n \sum_{l=1}^L \lambda_l V_{n,i}(\tau_l)$. Observe that $\text{var}(\sum_{i=1}^n \sum_{l=1}^L \lambda_l V_{n,i}(\tau_l)) =: \sigma_{n,L}^2$ where
$$\sigma_{n,L}^2 = \sum_{l,l'=1}^L \lambda_l\lambda_{l'}\mathbf{u}_n^\top J_m(\tau_l)^{-1}\mathbb{E}[\mathbf{Z}_i\mathbf{Z}_i^\top]J_m(\tau_{l'})^{-1}\mathbf{u}_n(\tau_l \wedge \tau_{l'} - \tau_l\tau_{l'}).$$

If $0 = \lim_{n\to\infty} \sigma_{n,L}^2 = \sum_{l,l'=1}^L \lambda_l\lambda_{l'} H(\tau_l,\tau_{l'};\mathbf{u}_n) = \text{var}(\sum_{l=1}^L \lambda_l\mathbb{G}(\tau_l))$, then by Markov's inequality $\sum_{i=1}^n \sum_{l=1}^L \lambda_l V_{n,i}(\tau_l) \to 0$ in probability, which coincides with the distribution of $\sum_{l=1}^L \lambda_l\mathbb{G}(\tau_l)$, which is a single point mass at 0. Next, consider the case $\sigma_{n,L}^2 \to \sigma_L^2 > 0$. For sufficiently large $n$ and arbitrary $v > 0$, Markov's inequality implies
$$\sigma_{n,L}^{-2}\sum_{i=1}^n \mathbb{E}\left[\left(\sum_{l=1}^L \lambda_l V_{n,i}(\tau_l)\right)^2 \mathbf{1}\left(\sum_{l=1}^L \lambda_l V_{n,i}(\tau_l) > v\right)\right]$$
$$\lesssim \xi_m^2 n^{-1}\sigma_{n,L}^{-2}\sum_{i=1}^n \mathbb{E}\left[\mathbf{1}\left(\sum_{l=1}^L \lambda_l V_{n,i}(\tau_l) > v\right)\right]$$
$$\lesssim \xi_m^2 n^{-1}\sigma_{n,L}^{-2}v^{-2}\sum_{l,l'=1}^L \lambda_l\lambda_{l'}\mathbf{u}_n^\top J_m(\tau_l)^{-1}\mathbb{E}[\mathbf{Z}_i\mathbf{Z}_i^\top]J_m(\tau_{l'})^{-1}\mathbf{u}_n(\tau_l \wedge \tau_{l'} - \tau_l\tau_{l'})$$
$$= o(1)$$

since $\xi_m^2 n^{-1} = o(1)$ by the assumption $m\xi_m^2 \log n = o(n)$. Hence the Lindeberg condition is verified. The finite dimensional convergence follows from Cramér-Wold devise. This completes the proof. □

### A.2. Proofs of Theorem 2.4, Example 2.3 and Remark 2.6

We begin by introducing some notations and useful preliminary results. For a vector $\mathbf{u} = (u_1, ..., u_m)^\top$ and a set $\mathcal{I} \subset \{1, ..., m\}$, let $\mathbf{u}^{(\mathcal{I})} \in \mathbb{R}^m$ denote the vector that has entries $u_i$ for $i \in \mathcal{I}$ and zero otherwise. For a vector $\mathbf{a} \in \mathbb{R}^m$, let $k_\mathbf{a}$ denote the position of the first non-zero entry of $\mathbf{a}$ with $\|\mathbf{a}\|_0$ non-zero consecutive entries

$$\mathcal{I}(\mathbf{a}, D) := \{i : |i - k_\mathbf{a}| \leq \|\mathbf{a}\|_0 + D\}, \tag{A.2}$$

$$\mathcal{I}'(\mathbf{a}, D) := \{1 \leq j \leq m : \exists i \in \mathcal{I}(\mathbf{a}, D) \text{ such that } |j - i| \leq \|\mathbf{a}\|_0\}, \tag{A.3}$$

**Lemma A.1.** *Under (L), for an arbitrary vector $\mathbf{a} \in \mathbb{R}^m$ with at most $\|\mathbf{a}\|_0$ non-zero consecutive entries we have for a constant $\gamma \in (0, 1)$ independent of $n, \tau$*

$$|(\mathbf{a}^\top \widetilde{J}_m^{-1}(\tau))_j| \leq C_1 \|\mathbf{a}\|_\infty \sum_{q=k_\mathbf{a}}^{k_\mathbf{a} + \|\mathbf{a}\|_0} \gamma^{|q-j|}. \tag{A.4}$$

$$\|\mathbf{a}^\top \widetilde{J}_m^{-1}(\tau) - (\mathbf{a}^\top \widetilde{J}_m^{-1}(\tau))^{(\mathcal{I}(\mathbf{a},D))}\| \lesssim \|\mathbf{a}\|_\infty \|\mathbf{a}\|_0 \gamma^D \tag{A.5}$$

$$\|\mathbf{a}^\top J_m^{-1}(\tau) - (\mathbf{a}^\top J_m^{-1}(\tau))^{(\mathcal{I}(\mathbf{a},D))}\| \lesssim \|\mathbf{a}\|_\infty \|\mathbf{a}\|_0 \gamma^D \tag{A.6}$$

**Proof for Lemma A.1.** Under (L) the matrix $\mathbf{Z}(x)\mathbf{Z}(x)^\top$ has no non-zero entries that are further than $r$ away from the diagonal. Thus $\widetilde{J}_m^{-1}$ is a band matrix with band width no larger that $2r$. Apply similar arguments as in the proof of Lemma 6.3 in Zhou et al. (1998) to find that under (L) the entries of $\widetilde{J}_m^{-1}(\tau)$ satisfy

$$\sup_{\tau, m} |(\widetilde{J}_m^{-1}(\tau))_{j,k}| \leq C_1 \gamma^{|j-k|}$$

for some $\gamma \in (0, 1)$ and a constant $C_1$ where both $\gamma$ and $C_1$ do not depend on $n, \tau$. It follows that

$$|(\mathbf{a}^\top \widetilde{J}_m^{-1}(\tau))_j| \leq C_1 \|\mathbf{a}\|_\infty \sum_{q=k_\mathbf{a}}^{k_\mathbf{a} + \|\mathbf{a}\|_0} \gamma^{|q-j|}.$$

and thus (A.4) is established. For a proof of (A.5) note that by (A.4) we have

$$|(\mathbf{a}^\top \widetilde{J}_m^{-1}(\tau))_j| \leq C_1 \|\mathbf{a}\|_\infty \sum_{q=k_\mathbf{a}}^{k_\mathbf{a} + \|\mathbf{a}\|_0} \gamma^{|q-j|} \leq C_1 \|\mathbf{a}\|_\infty \|\mathbf{a}\|_0 \gamma^{|k_\mathbf{a} - j| - \|\mathbf{a}\|_0}.$$

By the definition of $\mathcal{I}(\mathbf{a}, D)$ we find

$$\|\mathbf{a}^\top \widetilde{J}_m^{-1}(\tau) - (\mathbf{a}^\top \widetilde{J}_m^{-1}(\tau))^{(\mathcal{I}(\mathbf{a},D))}\| \leq C \|\mathbf{a}\|_\infty \|\mathbf{a}\|_0 \gamma^D$$

for a constant $C$ independent of $n$. The proof of (A.6) is similar to the proof of (A.5). □

**Proof for Theorem 2.4.** By Theorem 5.2 and Condition (B1), we first obtain

$$\mathbf{u}_n^\top (\widehat{\boldsymbol{\beta}}(\tau) - \boldsymbol{\beta}(\tau)) = -\frac{1}{n} \mathbf{u}_n^\top \widetilde{J}_m(\tau)^{-1} \sum_{i=1}^n \mathbf{B}_i (\mathbf{1}\{Y_i \leq \mathbf{B}(x)^\top \boldsymbol{\beta}_n(\tau)\} - \tau) + o_P(\|\mathbf{u}_n\| n^{-1/2}). \tag{A.7}$$

Let $\widetilde{U}_{1,n}(\tau) := n^{-1} \mathbf{u}_n^\top \widetilde{J}_m(\tau)^{-1} \sum_{i=1}^n \mathbf{B}_i (\mathbf{1}\{Y_i \leq \mathbf{B}(x)^\top \boldsymbol{\beta}_n(\tau)\} - \tau)$. We claim that

$$\mathbf{u}_n^\top (\widetilde{U}_{1,n}(\tau) - U_n(\tau)) = o_P(\|\mathbf{u}_n\| n^{-1/2}), \tag{A.8}$$

where $\boldsymbol{U}_n(\tau) := n^{-1}\mathbf{u}_n^\top J_m(\tau)^{-1}\sum_{i=1}^n \mathbf{B}_i(\mathbf{1}\{Y_i \leq Q(X_i;\tau)\} - \tau)$. Given (A.8), the process convergence of $\mathbf{u}_n^\top(\widehat{\boldsymbol{\beta}}(\tau) - \boldsymbol{\beta}(\tau))$ and continuity of the sample paths of the limiting process follows from process convergence of $\mathbf{u}_n^\top \boldsymbol{U}_n(\tau)$, which can be shown via exactly the same steps as in Section A.1.2 by replacing $\mathbf{Z}$ by $\mathbf{B}$ given assumptions (A1)-(A3).

To show (A.8), we proceed in two steps. Given $\mathbf{u}_n \in \mathcal{S}_\mathcal{I}^{m-1}$, let $\boldsymbol{U}_{1,n}(\tau) := n^{-1}\mathbf{u}_n^\top \widetilde{J}_m(\tau)^{-1}\sum_{i=1}^n \mathbf{B}_i(\mathbf{1}\{Y_i \leq Q(X_i;\tau)\} - \tau)$.

**Step 1:** $\sup_{\tau\in\mathcal{T}}|\mathbf{u}_n^\top(\widetilde{\boldsymbol{U}}_{1,n}(\tau) - \boldsymbol{U}_{1,n}(\tau))| = o_P(n^{-1/2})$, for all $\mathbf{u}_n \in \mathcal{S}_\mathcal{I}^{m-1}$.

Let $\widetilde{I}_0(\tau) := \mathbb{E}[\mathbf{u}_n^\top(\widetilde{\boldsymbol{U}}_{1,n}(\tau) - \boldsymbol{U}_{1,n}(\tau))]$ and observe the decomposition

$$\mathbf{u}_n^\top(\widetilde{\boldsymbol{U}}_{1,n}(\tau) - \boldsymbol{U}_{1,n}(\tau)) - \mathbb{E}[\mathbf{u}_n^\top(\widetilde{\boldsymbol{U}}_{1,n}(\tau) - \boldsymbol{U}_{1,n}(\tau))] + \mathbb{E}[\mathbf{u}_n^\top(\widetilde{\boldsymbol{U}}_{1,n}(\tau) - \boldsymbol{U}_{1,n}(\tau))]$$
$$= \left(\mathbf{u}_n\widetilde{J}_m^{-1}(\tau) - (\mathbf{u}_n\widetilde{J}_m^{-1}(\tau))^{(\mathcal{I}(\mathbf{u}_n,D))}\right)(\mathbb{P}_n - P)\{\mathbf{B}_i(\mathbf{1}\{Y_i \leq \mathbf{B}(x)^\top\boldsymbol{\beta}_n(\tau)\} - \mathbf{1}\{Y_i \leq Q(X_i;\tau)\})\}$$
$$+ (\mathbb{P}_n - P)\{(\mathbf{u}_n\widetilde{J}_m^{-1}(\tau))^{(\mathcal{I}(\mathbf{u}_n,D))}\mathbf{B}_i(\mathbf{1}\{Y_i \leq \mathbf{B}(x)^\top\boldsymbol{\beta}_n(\tau)\} - \mathbf{1}\{Y_i \leq Q(X_i;\tau)\})\} + \widetilde{I}_0(\tau)$$
$$=: \widetilde{I}_1(\tau) + \widetilde{I}_2(\tau) + \widetilde{I}_0(\tau).$$

For $\sup_{\tau\in\mathcal{T}}|\widetilde{I}_0(\tau)|$, by the construction of $\boldsymbol{\beta}_n(\tau)$ in (2.7),

$$\sup_{\tau\in\mathcal{T}}|\widetilde{I}_0(\tau)| \leq \sup_{\tau\in\mathcal{T}}\left|\mathbf{u}_n^\top \widetilde{J}_m^{-1}(\tau)\mathbb{E}\left[\mathbf{B}_i(\mathbf{1}\{Y_i \leq Q(X_i;\tau)\} - \mathbf{1}\{Y_i \leq \mathbf{B}_i^\top\boldsymbol{\beta}_n(\tau)\})\right]\right|$$
$$\leq \widetilde{c}_n^2 \bar{f}' \sup_{\mathbf{u}\in\mathcal{S}^{m-1}}\mathbb{E}|\mathbf{u}^\top \widetilde{J}_m^{-1}(\tau)\mathbf{B}| \leq \widetilde{c}_n^2 \bar{f}' \sup_{\mathbf{u}\in\mathcal{S}^{m-1}}\left(\mathbf{u}^\top \widetilde{J}_m^{-1}(\tau)\mathbb{E}[\mathbf{B}\mathbf{B}^\top]\widetilde{J}_m^{-1}(\tau)\mathbf{u}\right)^{1/2} = o(n^{-1/2}),$$

where the final rate follows from assumptions (A2) and $\widetilde{c}_n^2 = o(n^{-1/2})$ in (B1).

By (A.5) in Lemma A.1, let $D = c\log n$ for large enough $c > 0$, we have almost surely

$$\sup_{\tau\in\mathcal{T}}|\widetilde{I}_1(\tau)| \leq \sup_{\tau\in\mathcal{T}}\|\mathbf{u}_n^\top \widetilde{J}_m^{-1}(\tau) - (\mathbf{u}_n^\top \widetilde{J}_m^{-1}(\tau))^{(\mathcal{I}(\mathbf{u}_n,D))}\|\xi_m \leq \|\mathbf{u}_n\|_\infty \|\mathbf{u}_n\|_0 n^{c\log\gamma}\xi_m = o(n^{-1/2}).$$

For bounding $\sup_{\tau\in\mathcal{T}}|\widetilde{I}_2(\tau)|$, observe that

$$\frac{1}{n}\sum_{i=1}^n (\mathbf{u}_n^\top \widetilde{J}_m^{-1}(\tau))^{(\mathcal{I}(\mathbf{u}_n,D))}\mathbf{B}_i(\mathbf{1}\{Y_i \leq \mathbf{B}_i^\top\boldsymbol{\beta}_n(\tau)\} - \mathbf{1}\{Y_i \leq Q(X_i;\tau)\})$$
$$= \frac{1}{n}\sum_{i=1}^n \mathbf{u}_n^\top \widetilde{J}_m^{-1}(\tau)\mathbf{B}_i^{(\mathcal{I}(\mathbf{u}_n,D))}(\mathbf{1}\{Y_i \leq \mathbf{B}_i^\top\boldsymbol{\beta}_n(\tau)\} - \mathbf{1}\{Y_i \leq Q(X_i;\tau)\})$$
$$= \frac{1}{n}\sum_{\{i:\mathrm{supp}(\mathbf{B}_i)\cap\mathcal{I}(\mathbf{u}_n,D)\neq\emptyset\}} \mathbf{u}_n^\top \widetilde{J}_m^{-1}(\tau)\mathbf{B}_i^{(\mathcal{I}(\mathbf{u}_n,D))}(\mathbf{1}\{Y_i \leq \mathbf{B}_i^\top\boldsymbol{\beta}_n(\tau)\} - \mathbf{1}\{Y_i \leq Q(X_i;\tau)\}) \quad (A.9)$$
$$= \frac{1}{n}\sum_{\{i:\mathbf{B}_i^{(\mathcal{I}'(\mathbf{u}_n,D)^c)}=\mathbf{0}\}} \mathbf{u}_n^\top \widetilde{J}_m^{-1}(\tau)\mathbf{B}_i^{(\mathcal{I}(\mathbf{u}_n,D))}(\mathbf{1}\{Y_i \leq \mathbf{B}_i^\top\boldsymbol{\beta}_n(\tau)\} - \mathbf{1}\{Y_i \leq Q(X_i;\tau)\})$$
$$= \frac{1}{n}\sum_{i=1}^n \mathbf{u}_n^\top \widetilde{J}_m^{-1}(\tau)\mathbf{B}_i^{(\mathcal{I}(\mathbf{u}_n,D))}(\mathbf{1}\{Y_i \leq \left(\mathbf{B}_i^{(\mathcal{I}'(\mathbf{u}_n,D))}\right)^\top\boldsymbol{\beta}_n(\tau)\} - \mathbf{1}\{Y_i \leq Q(X_i;\tau)\}),$$
$$= \frac{1}{n}\sum_{i=1}^n \mathbf{u}_n^\top \widetilde{J}_m^{-1}(\tau)\mathbf{B}_i^{(\mathcal{I}(\mathbf{u}_n,D))}(\mathbf{1}\{Y_i \leq \mathbf{B}_i^\top\boldsymbol{\beta}_n(\tau)^{(\mathcal{I}'(\mathbf{u}_n,D))}\} - \mathbf{1}\{Y_i \leq Q(X_i;\tau)\}),$$

where the third equality follows from the fact that $\mathbf{B}_i^{(\mathcal{I}(\mathbf{u}_n,D))} \neq \mathbf{0}$ can only happen for $i \in \{i : \mathbf{B}_i^{(\mathcal{I}'(\mathbf{u}_n,D)^c)} = \mathbf{0}\}$, because $\mathbf{B}$ can only be nonzero in $r$ *consecutive* entries by assumption (L), where $\mathcal{I}'(\mathbf{u}_n,D)^c = \{1,...,m\} - \mathcal{I}'(\mathbf{u}_n,D)$ is the complement of $\mathcal{I}'(\mathbf{u}_n,D)$ in $\{1,...,m\}$. By restricting ourselves on set $\{i : \mathbf{B}_i^{(\mathcal{I}'(\mathbf{u}_n,D)^c)} = \mathbf{0}\}$, it is enough to look at the coefficient $\boldsymbol{\beta}_n(\tau)^{(\mathcal{I}'(\mathbf{u}_n,D))}$ in the last equality in (A.9). Hence,

$$\sup_{\tau\in\mathcal{T}}|\widetilde{I}_2(\tau)| \leq \|\mathbb{P}_n - P\|_{\widetilde{\mathcal{G}}_5(\mathcal{I}(\mathbf{u}_n,D),\mathcal{I}'(\mathbf{u}_n,D))}$$

where for any two index sets $\mathcal{I}_1$ and $\mathcal{I}_1'$

$$\widetilde{\mathcal{G}}_5(\mathcal{I}_1, \mathcal{I}_1') = \left\{ (X, Y) \mapsto \mathbf{a}^\top \mathbf{B}(X)^{(\mathcal{I}_1)} \big(\mathbf{1}\{Y \leq \mathbf{B}(X)^\top \mathbf{b}^{(\mathcal{I}_1')}\} - \mathbf{1}\{Y \leq Q(X;\tau)\}\big) \big| \tau \in \mathcal{T}, \mathbf{b} \in \mathbb{R}^m, \mathbf{a} \in \mathcal{S}^{m-1} \right\}.$$

With the choice of $D = c \log n$, the cardinality of both $\mathcal{I}(\mathbf{u}_n, D)$ and $\mathcal{I}'(\mathbf{u}_n, D)$ is of order $O(\log n)$. Hence, the VC index of $\widetilde{\mathcal{G}}_5(\mathcal{I}(\mathbf{u}_n, D), \mathcal{I}'(\mathbf{u}_n, D))$ is bounded by $O(\log n)$. Note that for any $f \in \widetilde{\mathcal{G}}_5(\mathcal{I}(\mathbf{u}_n, D), \mathcal{I}'(\mathbf{u}_n, D))$, $|f| \lesssim \xi_m$ and $\|f\|_{L_2(P)} \lesssim \widetilde{c}_n$. Applying (S.2.3) yields

$$P\left( \sup_{\tau \in \mathcal{T}} |\widetilde{I}_2(\tau)| \leq C\left[ \left( \frac{\widetilde{c}_n (\log n)^2}{n} \right)^{1/2} + \frac{\xi_m (\log n)^2}{n} + \left( \frac{\widetilde{c}_n \kappa_n}{n} \right)^{1/2} + \frac{\kappa_n \xi_m}{n} \right] \right) \geq 1 - e^{\kappa_n}.$$

Taking $\kappa_n = C \log n$, $\widetilde{c}_n^2 = o(n^{-1/2})$ and $\xi_m^4 (\log n)^6 = o(n)$ in (B1) implies that $\sup_{\tau \in \mathcal{T}} |\widetilde{I}_2(\tau)| = o_P(n^{1/2})$.

**Step 2:** $\sup_{\tau \in \mathcal{T}} |\mathbf{u}_n^\top (\widetilde{U}_{1,n}(\tau) - U_n(\tau))| = o_P(n^{-1/2})$, for all $\mathbf{u}_n \in \mathcal{S}_\mathcal{I}^{m-1}$.
Observe that

$$\mathbf{u}_n^\top (\widetilde{U}_{1,n}(\tau) - U_n(\tau))$$
$$= \frac{1}{n} \Big( \{\mathbf{u}_n^\top (\widetilde{J}_m(\tau)^{-1} - J_m(\tau)^{-1})\} - \{\mathbf{u}_n^\top (\widetilde{J}_m(\tau)^{-1} - J_m(\tau)^{-1})\}^{(\mathcal{I}(\mathbf{u}_n, D))} \Big) \sum_{i=1}^n \mathbf{B}_i (\mathbf{1}\{Y_i \leq Q(X_i; \tau)\} - \tau)$$
$$+ \frac{1}{n} \{\mathbf{u}_n^\top (\widetilde{J}_m(\tau)^{-1} - J_m(\tau)^{-1})\}^{(\mathcal{I}(\mathbf{u}_n, D))} \sum_{i=1}^n \mathbf{B}_i (\mathbf{1}\{Y_i \leq Q(X_i; \tau)\} - \tau)$$
$$=: \widetilde{I}_3(\tau) + \widetilde{I}_4(\tau).$$

Applying (A.5) and (A.6) in Lemma A.1 with $D = c \log n$ where $c > 0$ is chosen sufficiently large, we have almost surely

$$\sup_{\tau \in \mathcal{T}} |\widetilde{I}_3(\tau)| \leq \left( \sup_{\tau \in \mathcal{T}} \|\mathbf{u}_n^\top \widetilde{J}_m^{-1}(\tau) - (\mathbf{u}_n^\top \widetilde{J}_m^{-1}(\tau))^{(\mathcal{I}(\mathbf{u}_n, D))}\| + \sup_{\tau \in \mathcal{T}} \|\mathbf{u}_n^\top J_m^{-1}(\tau) - (\mathbf{u}_n^\top J_m^{-1}(\tau))^{(\mathcal{I}(\mathbf{u}_n, D))}\| \right) \xi_m$$
$$\leq 2\|\mathbf{u}_n\|_\infty \|\mathbf{u}_n\|_0 n^{c \log \gamma} \xi_m = o(n^{-1/2}). \tag{A.10}$$

Now it is left to bound $\sup_{\tau \in \mathcal{T}} |\widetilde{I}_4(\tau)|$. We have

$$\widetilde{I}_4(\tau) = \frac{1}{n} \sum_{i=1}^n \{\mathbf{u}_n^\top (\widetilde{J}_m(\tau)^{-1} - J_m(\tau)^{-1})\}^{(\mathcal{I}(\mathbf{u}_n, D))} \mathbf{B}_i (\mathbf{1}\{Y_i \leq Q(X_i; \tau)\} - \tau)$$
$$= \frac{1}{n} \sum_{i=1}^n \mathbf{u}_n^\top (\widetilde{J}_m(\tau)^{-1} - J_m(\tau)^{-1}) \mathbf{B}_i^{(\mathcal{I}(\mathbf{u}_n, D))} (\mathbf{1}\{Y_i \leq Q(X_i; \tau)\} - \tau).$$

Hence,

$$\sup_{\tau \in \mathcal{T}} |\widetilde{I}_4(\tau)| \leq \sup_{\tau \in \mathcal{T}} \|\mathbf{u}_n^\top (\widetilde{J}_m(\tau)^{-1} - J_m(\tau)^{-1})\| \|\mathbb{P}_n - P\|_{\mathcal{G}_0(\mathcal{I}(\mathbf{u}_n, D)) \cdot \mathcal{G}_4}$$

where for any $\mathcal{I}$,

$$\mathcal{G}_0(\mathcal{I}) := \{(\mathbf{B}, Y) \mapsto \mathbf{a}^\top \mathbf{B}^{(\mathcal{I})} \mathbf{1}\{\|\mathbf{B}\| \leq \xi_m\} | \mathbf{a} \in \mathcal{S}^{m-1}\},$$
$$\mathcal{G}_4 := \{(X, Y) \mapsto \mathbf{1}\{Y_i \leq Q(X; \tau)\} - \tau | \tau \in \mathcal{T}\}.$$

The cardinality of the set $\mathcal{I}(\mathbf{u}_n, c \log n)$ is of order $O(\log n)$. Thus, the VC index for $\mathcal{G}_0(\mathcal{I}(\mathbf{u}_n, D))$ is of order $O(\log n)$. The VC index of $\mathcal{G}_4$ is 2 (see Lemma S.2.4). By Lemma S.2.2,

$$N(\mathcal{G}_0(\mathcal{I}(\mathbf{u}_n, D)) \cdot \mathcal{G}_4, L_2(\mathbb{P}_n); \varepsilon) \leq \left( \frac{A \|F\|_{L_2(\mathbb{P}_n)}}{\varepsilon} \right)^{v_0(n)},$$

where $v_0(n) = O(\log n)$. In addition, for any $f \in \mathcal{G}_0(\mathcal{I}(\mathbf{u}_n, D)) \cdot \mathcal{G}_4$, $|f| \lesssim \xi_m$ and $\|f\|_{L_2(P)} = O(1)$ by (A1). Furthermore, by assumptions (A1)-(A2) and the definition of $\widetilde{c}_n$,

$$\left\|\mathbf{u}_n^\top \left(\widetilde{J}_m(\tau)^{-1} - J_m(\tau)^{-1}\right)\right\| \leq \widetilde{c}_n \lambda_{\max}(\mathbb{E}[\mathbf{B}(X)\mathbf{B}(X)^\top]) \bar{f}' \lesssim \widetilde{c}_n. \tag{A.11}$$

By (S.2.3), we have for some constant $C > 0$,

$$P\left(\sup_{\tau \in \mathcal{T}} |\widetilde{I}_2(\tau)| \leq C\widetilde{c}_n \left[\left(\frac{(\log n)^2}{n}\right)^{1/2} + \frac{\xi_m (\log n)^2}{n} + \left(\frac{\kappa_n}{n}\right)^{1/2} + \frac{\kappa_n \xi_m}{n}\right]\right) \geq 1 - e^{\kappa_n}.$$

Taking $\kappa_n = C \log n$, an application of (B1) completes the proof. $\square$

**Proof for Example 2.3.** As $\widetilde{J}_m(\tau)$ is a band matrix, applying similar arguments as in the proof of Lemma 6.3 in Zhou et al. (1998) gives

$$\sup_{\tau, m} |(\widetilde{J}_m^{-1}(\tau))_{j,j'}| \leq C_1 \gamma^{|j-j'|}, \tag{A.12}$$

for some $\gamma \in (0,1)$ and $C_1 > 0$. Let $k_{\mathbf{B}(x)}$ be the index of the first nonzero element of the vector $\mathbf{B}(x)$. Then by (A.12), we have

$$\sup_{\tau,m} |(\mathbf{B}(x)^\top \widetilde{J}_m^{-1}(\tau))_j| \leq C_1 \|\mathbf{B}(x)\|_\infty \sum_{j'=k_{\mathbf{u}_n}}^{k_{\mathbf{B}(x)} + \|\mathbf{B}(x)\|_0} \gamma^{|j'-j|},$$

and also

$$\sup_{\tau,m} \mathbb{E}|\mathbf{B}(x)^\top \widetilde{J}_m^{-1}(\tau) \mathbf{B}(X)| \leq C_1 \|\mathbf{B}(x)\|_\infty \max_{l \leq m} \mathbb{E}|B_l(X)| \sum_{j=1}^{m} \sum_{j'=k_{\mathbf{B}(x)}}^{k_{\mathbf{B}(x)} + \|\mathbf{B}(x)\|_0} \gamma^{|j'-j|}. \tag{A.13}$$

Since $\|\mathbf{B}(x)\|_0$ is bounded by a constant, the sum in (A.13) is bounded uniformly. Moreover, in the present setting we have $\|\mathbf{B}(x)\|_\infty = O(m^{1/2})$ and $\max_{l \leq m} \mathbb{E}|B_l(X)| = O(m^{-1/2})$. Therefore, for each $m$ we have

$$\sup_{\tau \in \mathcal{T}, x \in \mathcal{X}} \mathbb{E}|\mathbf{B}(x)^\top \widetilde{J}_m^{-1}(\tau) \mathbf{B}(X)| = O(1).$$

$\square$

**Proof of Remark 2.6.** Consider the product $B_j(x) B_{j'}(x)$ of two B-spline functions. The fact that $B_j(x)$ is locally supported on $[t_j, t_{j+r}]$ implies that for all $j'$ satisfying $|j - j'| \geq r$, $B_j(x) B_{j'}(x) = 0$ for all $x$, where $r \in \mathbb{N}$ is the degree of spline. This also implies $J_m(\tau)$ and $\mathbb{E}[\mathbf{BB}^\top]$ are a band matrices with each column having at most $L_r := 2r + 1$ nonzero elements and each non-zero element is at most $r$ entries away from the main diagonal. Recall also the fact that $\max_{j \leq m} \sup_{t \in \mathbb{R}} |B_j(t)| \lesssim m^{1/2}$ (by the discussion following assumption (A1)).

Define $J_{m,D}(\tau) := D_m(\tau) \mathbb{E}[\mathbf{BB}^\top]$, where matrix $D_m(\tau) := diag(f_{Y|X}(Q(t_j; \tau)|t_j), j = 1, ..., m)$, and $R_m(\tau) := J_m(\tau) - J_{m,D}(\tau)$. Both $J_{m,D}(\tau)$ and $R_m(\tau)$ have the same band structure as $J_m(\tau)$. For arbitrary $j, j' = 1, ..., m$, $\tau \in \mathcal{T}$ and a universal constant $C > 0$,

$$\begin{aligned}
&|(R_m(\tau))_{j,j'}| \\
&= \left|\mathbb{E}\left[B_j(X) B_{j'}(X) \{f_{Y|X}(Q(X;\tau)|X) - f_{Y|X}(Q(t_j;\tau)|t_j)\}\right]\right| \\
&\leq 2 \max_{j \leq m} \sup_{t \in \mathbb{R}} |B_j(t)|^2 \int_0^1 \mathbf{1}\left\{|x - t_j| \leq C\frac{r}{m}\right\} |f_{Y|X}(Q(x;\tau)|x) - f_{Y|X}(Q(t_j;\tau)|t_j)| f_X(x) dx \\
&\leq 2Cm \int_0^1 \mathbf{1}\left\{|x - t_j| \leq C\frac{r}{m}\right\} |x - t_j| dx \\
&= O(m^{-1}),
\end{aligned} \tag{A.14}$$

where the second inequality is an application of the upper bound of $\max_{j\le m}\sup_{t\in\mathbb{R}}|B_j(t)| \lesssim m^{1/2}$ and the local support property of $B_j$; the third inequality follows by the assumption (2.10) and bounded $f_X(x)$. This shows that $\max_{j,j'=1,\ldots,m}\sup_{\tau\in\mathcal{T}}|(R_m(\tau))_{j,j'}| = O(m^{-1})$.

Now we show a stronger result that $\sup_{\tau\in\mathcal{T}}\|R_m(\tau)\| = O(m^{-1/2})$ for later use. Let $\mathbf{v} = (v_1,\ldots,v_m)$. Denote $k_j$ the index with the first nonzero entry in the $j$th column of $R_m(\tau)$. By the band structure of $R_m(\tau)$,

$$\sup_{\tau\in\mathcal{T}}\|R_m(\tau)\|^2 = \sup_{\tau\in\mathcal{T}}\sup_{\mathbf{v}\in\mathcal{S}^{m-1}}\|R_m(\tau)\mathbf{v}\|_2^2 = \sup_{\tau\in\mathcal{T}}\sup_{\mathbf{v}\in\mathcal{S}^{m-1}}\sum_{j=1}^{m}\left(\sum_{i=k_j}^{k_j+L_r-1} v_i(R_m(\tau))_{i,j}\right)^2$$
$$\lesssim \sup_{\tau\in\mathcal{T}}\max_{j,j'}|(R_m(\tau))_{j,j'}|^2 m = O(m^{-1}), \tag{A.15}$$

where the last equality follows by (A.14). Note that from assumptions (A1)-(A3) that

$$f_{\min}M^{-1} < \lambda_{\min}(J_{m,D}(\tau)) \le \lambda_{\max}(J_{m,D}(\tau)) < \bar{f}M \tag{A.16}$$

uniformly in $\tau \in \mathcal{T}$, where the constant $M > 0$ is defined as in Assumption (A1). Using (A.16), assumptions (A1)-(A3) and (A.15),

$$\|J_m^{-1}(\tau) - J_{m,D}^{-1}(\tau)\| \le \|J_{m,D}^{-1}(\tau)\|\|J_{m,D}(\tau) - J_m(\tau)\|\|J_m^{-1}(\tau)\| \lesssim \sup_{\tau\in\mathcal{T}}\|R_m(\tau)\| = O(m^{-1/2})$$

uniformly in $\tau \in \mathcal{T}$.

Without loss of generality, from now on we drop the term $\tau_1 \wedge \tau_2 - \tau_1\tau_2$ out of our discussion and focus on the matrix part in the covariance function $\widetilde{H}(\tau_1,\tau_2;\mathbf{u}_n)$ defined in (2.8). From (A1) we have $\|\mathbb{E}[\mathbf{B}\mathbf{B}^\top]\| < M$ for some constant $M > 0$ so for any $\tau_1, \tau_2 \in \mathcal{T}$,

$$\|\mathbf{u}_n\|^{-2}\left|\mathbf{u}_n^\top J_m^{-1}(\tau_1)\mathbb{E}[\mathbf{B}\mathbf{B}^\top]J_m^{-1}(\tau_2)\mathbf{u}_n - \mathbf{u}_n^\top J_{m,D}(\tau_1)^{-1}\mathbb{E}[\mathbf{B}\mathbf{B}^\top]J_{m,D}(\tau_2)^{-1}\mathbf{u}_n\right|$$
$$\lesssim \sup_{\tau\in\mathcal{T}}\|R_m(\tau)\|\sup_{\tau\in\mathcal{T}}\|\mathbb{E}[\mathbf{B}\mathbf{B}^\top]J_{m,D}(\tau)^{-1}\| + \sup_{\tau\in\mathcal{T}}\|R_m(\tau)\|\sup_{\tau\in\mathcal{T}}\|\mathbb{E}[\mathbf{B}\mathbf{B}^\top]J_m(\tau)^{-1}\|$$
$$= O(m^{-1/2}). \tag{A.17}$$

Moreover, note that

$$\mathbf{u}_n^\top J_{m,D}(\tau_1)^{-1}\mathbb{E}[\mathbf{B}\mathbf{B}^\top]J_{m,D}(\tau_2)^{-1}\mathbf{u}_n = \mathbf{u}_n^\top D_m(\tau_1)^{-1}\mathbb{E}[\mathbf{B}\mathbf{B}^\top]^{-1}D_m(\tau_2)^{-1}\mathbf{u}_n \tag{A.18}$$

If $\mathbf{u}_n = \mathbf{B}(x)$, observe that for $l = 1,\ldots,m$, as suggested by the local support property, we only need to focus on the index $l$ satisfying $|x - t_l| \le Cr/m$, for a universal constant $C > 0$. We have

$$(\mathbf{B}(x)^\top D_m(\tau)^{-1})_l = B_l(x)f_{Y|X}(Q(t_l;\tau)|t_l)^{-1} = B_l(x)f_{Y|X}(Q(x;\tau)|x)^{-1} + R'(t_l), \tag{A.19}$$

where by assumption (2.10), $|R'(t_l)| \le \max_{j\le m}\sup_{t\in\mathbb{R}}|B_j(t)|Cf_{\min}^{-2}|x - t_l| = O(m^{-1/2})$. Therefore, the sparse vector $\mathbf{B}(x)^\top D_m(\tau)^{-1} = f_{Y|X}(Q(x;\tau)|x)^{-1}\mathbf{B}(x)^\top + \mathbf{a}_{\mathbf{B}(x)}$, where $\mathbf{a}_{\mathbf{B}(x)} \in \mathbb{R}^m$ is a vector with the same support as $\mathbf{B}(x)$ (only $r < \infty$ nonzero components) and with nonzero components of order $O(m^{-1/2})$. Hence, $\|\mathbf{a}_{\mathbf{B}(x)}\| = O(m^{-1/2})$. Continued from (A.18), for any $x \in [0,1]$,

$$\|\mathbf{B}(x)\|^{-2}\left|\mathbf{B}(x)^\top D_m(\tau_1)^{-1}\mathbb{E}[\mathbf{B}\mathbf{B}^\top]^{-1}D_m(\tau_2)^{-1}\mathbf{B}(x) - \frac{\mathbf{B}(x)^\top\mathbb{E}[\mathbf{B}\mathbf{B}^\top]^{-1}\mathbf{B}(x)}{f_{Y|X}(Q(x;\tau_1)|x)f_{Y|X}(Q(x;\tau_2)|x)}\right|$$
$$\le \|\mathbf{B}(x)\|^{-1}\|\mathbf{a}_{\mathbf{B}(x)}\|\sup_{\tau\in\mathcal{T}}\|\mathbb{E}[\mathbf{B}\mathbf{B}^\top]^{-1}D_m(\tau)^{-1}\| + \|\mathbf{B}(x)\|^{-1}\|\mathbf{a}_{\mathbf{B}(x)}\|\|\mathbb{E}[\mathbf{B}\mathbf{B}^\top]^{-1}\|/f_{\min}$$
$$= O(m^{-1}).$$

We observe that $\mathbf{B}(x)^\top\mathbb{E}[\mathbf{B}\mathbf{B}^\top]^{-1}\mathbf{B}(x)$ does not depend on $\tau_1$ and $\tau_2$ and can be treated as a scaling factor and shifted out of the covariance function as (2.11). Therefore, we finish the proof.

□

### A.3. Proof of Theorem 3.1

Observe that
$$\widehat{\boldsymbol{\alpha}}_j(\cdot) - \boldsymbol{\alpha}_j(\cdot) = \mathbf{e}_j^\top (\widehat{\boldsymbol{\gamma}}^\dagger(\cdot) - \boldsymbol{\gamma}_n^\dagger(\cdot))$$
where $\mathbf{e}_j$ denotes the $j$-th unit vector in $\mathbb{R}^{m+k}$ for $j = 1, .., m+k$, and
$$\mathbf{w}_n^\top (\widehat{\boldsymbol{\beta}}^\dagger(\tau) - \boldsymbol{\beta}_n^\dagger(\tau)) = (\mathbf{0}_k^\top, \mathbf{w}_n^\top)(\widehat{\boldsymbol{\gamma}}^\dagger(\cdot) - \boldsymbol{\gamma}_n^\dagger(\cdot)).$$

Let $h_n^\dagger(w, \tau) = \widetilde{\mathbf{Z}}(w)^\top \boldsymbol{\beta}_n^\dagger(\tau)$. The following results will be established at the end of the proof.

$$\sup_{\tau \in \mathcal{T}, j=1,\ldots,k} \left| \mathbb{E}\left[ \mathbf{e}_j^\top J_m(\tau)^{-1} \mathbf{Z}(\mathbf{1}\{Y \le Q(X; \tau)\} - \mathbf{1}\{Y \le \boldsymbol{\alpha}(\tau)^\top V + h_n^\dagger(W, \tau)\}) \right] \right| = o(n^{-1/2}), \quad (A.20)$$

$$\sup_{\tau \in \mathcal{T}} \left| (\mathbf{0}_k^\top, \mathbf{w}_n^\top / \|\mathbf{w}_n\|) J_m(\tau)^{-1} \mathbb{E}\left[ \mathbf{Z}(\mathbf{1}\{Y \le Q(X; \tau)\} - \mathbf{1}\{Y \le \boldsymbol{\alpha}(\tau)^\top V + h_n^\dagger(W, \tau)\}) \right] \right| = o(n^{-1/2}). \quad (A.21)$$

From Theorem 5.4, we obtain that under Condition (B1')

$$\mathbf{e}_j^\top (\widehat{\boldsymbol{\gamma}}^\dagger(\tau) - \boldsymbol{\gamma}_n^\dagger(\tau)) = -n^{-1/2} \mathbf{e}_j^\top J_m(\tau)^{-1} \mathbb{G}_n(\psi(\cdot; \boldsymbol{\gamma}_n^\dagger(\tau), \tau)) + o_P(n^{-1/2}), \quad j = 1, \ldots, k; \quad (A.22)$$

$$(\mathbf{0}_k^\top, \mathbf{w}_n^\top)(\widehat{\boldsymbol{\gamma}}^\dagger(\tau) - \boldsymbol{\gamma}_n^\dagger(\tau)) = -n^{-1/2}(\mathbf{0}_k^\top, \mathbf{w}_n^\top) J_m(\tau)^{-1} \mathbb{G}_n(\psi(\cdot; \boldsymbol{\gamma}_n^\dagger(\tau), \tau)) + o_P(n^{-1/2}) \quad (A.23)$$

uniformly in $\tau \in \mathcal{T}$. Equation (A.20) implies that for $j = 1, \ldots, k$

$$\mathbf{e}_j^\top J_m(\tau)^{-1} \mathbb{E}[\psi(Y_i, \mathbf{Z}_i; \boldsymbol{\gamma}_n^\dagger(\tau), \tau)] = \mathbf{e}_j^\top J_m(\tau)^{-1} \mathbb{E}[\mathbf{1}\{Y_i \le Q(X_i; \tau)\} - \tau] + o(n^{-1/2}) = o(n^{-1/2}).$$

Following similar arguments as given in the proof of (2.2) in Section A.1.1, (A.20) and (A.22) imply that

$$\mathbf{e}_j^\top (\widehat{\boldsymbol{\gamma}}^\dagger(\tau) - \boldsymbol{\gamma}_n^\dagger(\tau)) = -n^{-1} \mathbf{e}_j^\top J_m(\tau)^{-1} \sum_{i=1}^n \mathbf{Z}_i(\mathbf{1}\{Y_i \le Q(X_i; \tau)\} - \tau) + o_P(n^{-1/2}), \quad j = 1, \ldots, k.$$

uniformly in $\tau \in \mathcal{T}$. Similarly, by (A.21) and (A.23) we have

$$\|\mathbf{w}_n\|^{-1} \mathbf{w}_n^\top (\widehat{\boldsymbol{\beta}}^\dagger(\tau) - \boldsymbol{\beta}_n^\dagger(\tau)) = -n^{-1}(\mathbf{0}_k^\top, \|\mathbf{w}_n\|^{-1} \mathbf{w}_n^\top) J_m(\tau)^{-1} \sum_{i=1}^n \mathbf{Z}_i(\mathbf{1}\{Y_i \le Q(X_i; \tau)\} - \tau) + o_P(n^{-1/2}).$$

Thus, the claim will follow once we prove

$$\mathbb{G}_n(\cdot) := (\mathbb{G}_{n,1}(\cdot), \ldots, \mathbb{G}_{n,k}(\cdot), \mathbb{G}_{n,h}(\cdot)) \rightsquigarrow \mathbb{G}(\cdot) \text{ in } (\ell^\infty(\mathcal{T}))^{k+1}$$

where

$$\mathbb{G}_{n,j}(\tau) := -n^{-1/2} \mathbf{e}_j^\top J_m(\tau)^{-1} \sum_{i=1}^n \mathbf{Z}_i(\mathbf{1}\{Y_i \le Q(X_i; \tau)\} - \tau), \quad j = 1, \ldots, k$$

and

$$\mathbb{G}_{n,h}(\tau) := -\|\mathbf{w}_n\|^{-1} n^{-1/2} (\mathbf{0}_k^\top, \mathbf{w}_n^\top) J_m(\tau)^{-1} \sum_{i=1}^n \mathbf{Z}_i(\mathbf{1}\{Y_i \le Q(X_i; \tau)\} - \tau).$$

We need to establish tightness and finite dimensional convergence. By Lemma 1.4.3 of van der Vaart and Wellner (1996), it is enough to show the tightness of $\mathbb{G}_{n,j}$'s and $\mathbb{G}_{n,h}$ individually. Tightness follows from asymptotic equicontinuity which can be proved by an application of Lemma A.3. More precisely, apply Lemma A.3 with $\mathbf{u}_n = \mathbf{e}_j$ to prove tightness of $\mathbb{G}_{n,j}(\cdot)$ for $j = 1, \ldots, k$, and Lemma A.3 with $\mathbf{u}_n^\top = (\mathbf{0}_k^\top, \mathbf{w}_n^\top)$ to prove tightness of $\mathbb{G}_{n,h}(w_0; \cdot)$. Continuity of the sample paths of $\mathbb{G}_{n,h}(w_0; \cdot)$ follows by the same arguments as given at the beginning of Section A.1.2.

Next, we prove finite-dimensional convergence. Observe the decomposition

$$\mathbb{G}_n(\tau) = -n^{-1/2} \sum_{i=1}^n \left\{ \begin{pmatrix} M_1(\tau)^{-1}(V_i - A(\tau)\widetilde{\mathbf{Z}}(W_i)) \\ \|\mathbf{w}_n\|^{-1}\mathbf{w}_n^\top M_2(\tau)^{-1}\widetilde{\mathbf{Z}}(W_i) \end{pmatrix} \left(\mathbf{1}\{Y_i \leq Q(X_i;\tau)\} - \tau\right) + \begin{pmatrix} 0 \\ \varphi_i(\tau) \end{pmatrix} \right\}$$

where

$$\varphi_i(\tau) := -\|\mathbf{w}_n\|^{-1}\mathbf{w}_n^\top A(\tau)^\top M_1(\tau)^{-1}(V_i - A(\tau)\widetilde{\mathbf{Z}}(W_i))\left(\mathbf{1}\{Y_i \leq Q(X_i;\tau)\} - \tau\right).$$

By definition, we have $\mathbb{E}[\varphi_i(\tau)] = 0$ and moreover

$$\mathbb{E}[\varphi_i^2(\tau)] \lesssim \|\mathbf{w}_n\|^{-2}\mathbf{w}_n^\top A(\tau)^\top M_1(\tau)^{-1}\mathbb{E}[(V_i - A(\tau)\widetilde{\mathbf{Z}}(W_i))(V_i - A(\tau)\widetilde{\mathbf{Z}}(W_i))^\top]M_1(\tau)^{-1}A(\tau)\mathbf{w}_n.$$

Since $f_{Y|X}(Q(X_i;\tau)|X)$ is bounded away from zero uniformly, it follows that

$$\|\mathbb{E}[(V_i - A(\tau)\widetilde{\mathbf{Z}}(W_i))(V_i - A(\tau)\widetilde{\mathbf{Z}}(W_i))^\top]\| \leq f_{\min}\lambda_{\max}(M_1(\tau)) < \infty,$$

by Remark 5.5. Moreover, by Lemma A.2 proven later, $\|A(\tau)\mathbf{w}_n\| = O(1)$ uniformly in $\tau$, and thus by $\|\mathbf{w}_n\| \to \infty$, $\sup_{\tau \in \mathcal{T}} \mathbb{E}[\varphi_i^2(\tau)] = o(1)$. This implies that $n^{-1/2}\sum_{i=1}^n \varphi_i(\tau) = o_P(1)$ for every fixed $\tau \in \mathcal{T}$. Hence it suffices to prove finite dimensional convergence of

$$-n^{-1/2}\sum_{i=1}^n \begin{pmatrix} M_1(\tau)^{-1}(V_i - A(\tau)\widetilde{\mathbf{Z}}(W_i)) \\ \|\mathbf{w}_n\|^{-1}\mathbf{w}_n^\top M_2(\tau)^{-1}\widetilde{\mathbf{Z}}(W_i) \end{pmatrix} \left(\mathbf{1}\{Y_i \leq Q(X_i;\tau)\} - \tau\right).$$

Observe that $\mathbb{E}[M_1(\tau)^{-1}(h_{VW}(W_i;\tau) - A(\tau)\widetilde{\mathbf{Z}}(W_i))(\mathbf{1}\{Y_i \leq Q(X_i;\tau)\} - \tau)] = 0$ and by assumptions (A1)-(A3), (C1)

$$\sup_{\tau \in \mathcal{T}} \mathbb{E}[\|M_1(\tau)^{-1}(h_{VW}(W;\tau) - A(\tau)\widetilde{\mathbf{Z}}(W))\|^2]$$

$$\leq \sup_{\tau \in \mathcal{T}} \frac{1}{f_{\min}\lambda_{\min}(M_1(\tau))}\mathbb{E}[\|f_{Y|X}(Q(X;\tau)|X)(h_{VW}(w;\tau) - A(\tau)\widetilde{\mathbf{Z}}(W))\|^2] = o(1).$$

Thus, $n^{-1/2}\sum_{i=1}^n M_1(\tau)^{-1}(h_{VW}(W_i;\tau) - A(\tau)\widetilde{\mathbf{Z}}(W_i))(\mathbf{1}\{Y_i \leq Q(X_i;\tau)\} - \tau) = o_P(1)$ for every fixed $\tau \in \mathcal{T}$. So, now we only need to consider finite dimensional convergence of

$$\sum_{i=1}^n \boldsymbol{\psi}_i(\tau) := -n^{-1/2}\sum_{i=1}^n \begin{pmatrix} M_1(\tau)^{-1}(V_i - h_{VW}(W_i;\tau)) \\ \|\mathbf{w}_n\|^{-1}\mathbf{w}_n^\top M_2(\tau)^{-1}\widetilde{\mathbf{Z}}(W_i) \end{pmatrix} \left(\mathbf{1}\{Y_i \leq Q(X_i;\tau)\} - \tau\right). \tag{A.24}$$

Note that

$$\mathbb{E}[\boldsymbol{\psi}_i(\tau_1)\boldsymbol{\psi}_i(\tau_2)^\top] = (\tau_1 \wedge \tau_2 - \tau_1\tau_2)\begin{pmatrix} \Gamma_{11}(\tau_1,\tau_2) + o(1) & \Gamma_{12}(\tau_1,\tau_2) \\ \Gamma_{12}(\tau_1,\tau_2)^\top & \Gamma_{22}(\tau_1,\tau_2) + o(1) \end{pmatrix},$$

where $\Gamma_{12}(\tau_1,\tau_2) := \|\mathbf{w}_n\|^{-1}\mathbb{E}[M_1(\tau_2)^{-1}(V_i - h_{VW}(W_i;\tau_1))\mathbf{w}_n^\top M_2(\tau_2)^{-1}\widetilde{\mathbf{Z}}(W_i)]$. We shall now show that $\Gamma_{12}(\tau_1,\tau_2) = o(1)$ uniformly in $\tau_1,\tau_2$. Note that from the definition of $h_{VW}(W;\tau)$ in (3.4), by standard argument we can write

$$h_{VW}(W;\tau) = \mathbb{E}[f(Q(X;\tau)|X)|W]^{-1}\mathbb{E}[Vf(Q(X;\tau)|X)|W]. \tag{A.25}$$

From (A3) we obtain $\mathbb{E}[f(Q(X;\tau)|X)|W] \geq f_{\min} > 0$, and from (C1), (A2) it follows that $|\mathbb{E}[V^{(j)}f(Q(X;\tau)|X)|W]| \leq C\bar{f} = O(1)$, i.e. the components of $h_{VW}(W;\tau)$ are bounded by a constant almost surely. Hence,

$$\|\Gamma_{12}(\tau_1,\tau_2)\| = \|\mathbf{w}_n\|^{-1}\|M_1(\tau_1)^{-1}\mathbb{E}[(V_i - h_{VW}(W_i;\tau))\mathbf{w}_n^\top M_2(\tau_2)^{-1}\widetilde{\mathbf{Z}}(W_i)]\|$$

$$\leq \|\mathbf{w}_n\|^{-1}\|M_1(\tau_1)^{-1}\|\|\mathbb{E}[(V_i - h_{VW}(W_i;\tau_1))\mathbf{w}_n^\top M_2(\tau_2)^{-1}\widetilde{\mathbf{Z}}(W_i)]\|$$

$$\lesssim \|\mathbf{w}_n\|^{-1}\mathbb{E}[\|V_i - h_{VW}(W_i;\tau_1)\|\|\mathbf{w}_n^\top M_2(\tau_2)^{-1}\widetilde{\mathbf{Z}}(W_i)\|]$$

$$\lesssim \|\mathbf{w}_n\|^{-1}\mathbb{E}[|\mathbf{w}_n^\top M_2(\tau_2)^{-1}\widetilde{\mathbf{Z}}(W_i)|]$$

$$= o(1),$$

where the third inequality applies the lower bound for $\inf_{\tau \in \mathcal{T}} \lambda_{\min}(M_1(\tau))$ in Remark 5.5; the fourth inequality follows from $\sup_{1 \leq j \leq k, \tau} |V^{(j)}| + |h_{VW}^{(j)}(W; \tau)| < \infty$ a.s., while the last equality follows by the assumptions of the Theorem.

Now we prove the finite dimensional convergence (A.24). Taking arbitrary collections $\{\tau_1, ..., \tau_J\} \subset \mathcal{T}$, $\mathbf{c}_1, .., \mathbf{c}_J \in \mathbb{R}^{k+1}$, we need to show that

$$\sum_{j=1}^{J} \mathbf{c}_j^\top \mathbb{G}_n(\tau_j) = \sum_{j=1}^{J} \sum_{i=1}^{n} \mathbf{c}_j^\top \boldsymbol{\psi}_i(\tau_j) \xrightarrow{d} \sum_{j=1}^{J} \mathbf{c}_j^\top \mathbb{G}(\tau_j).$$

Define $V_{i,J} = \sum_{j=1}^{J} \mathbf{c}_j^\top \boldsymbol{\psi}_i(\tau_j)$. Note that $\mathbb{E}[V_{i,J}] = 0$ and $|V_{i,J}| \lesssim n^{-1/2} \xi_m$. Using the results derived above, we have

$$\mathrm{var}(V_{i,J}) = o(n^{-1}) + n^{-1} \sum_{j,j'=1}^{J} (\tau_j \wedge \tau_{j'} - \tau_j \tau_{j'}) \mathbf{c}_j^\top \Gamma(\tau_j, \tau_{j'}) \mathbf{c}_{j'},$$

where $\Gamma(\tau_j, \tau_{j'})$ is defined as (3.11). If $\sum_{j,j'=1}^{J} (\tau_j \wedge \tau_{j'} - \tau_j \tau_{j'}) \mathbf{c}_j^\top \Gamma(\tau_j, \tau_{j'}) \mathbf{c}_{j'} = 0$, then the distribution of $\sum_{j=1}^{J} \mathbf{c}_j^\top \mathbb{G}(\tau_j)$ is a single point mass at 0, and $\sum_{j=1}^{J} \mathbf{c}_j^\top \mathbb{G}_n(\tau_j) = \sum_{j=1}^{J} \sum_{i=1}^{n} \mathbf{c}_j^\top \boldsymbol{\psi}_i(\tau_j)$ converges to 0 in probability by Markov's inequality.

If $n^{-1} \sum_{j,j'=1}^{J} (\tau_j \wedge \tau_{j'} - \tau_j \tau_{j'}) \mathbf{c}_j^\top \Gamma(\tau_j, \tau_{j'}) \mathbf{c}_{j'} > 0$ define $s_{n,J}^2 = \sum_{i=1}^{n} \mathrm{var}(V_{i,J})$. We will now verify that the triangular array of random variables $(V_{i,J})_{i=1,...,n}$ satisfies the Lindeberg condition. For any $v > 0$ and sufficiently large $n$, Markov's inequality gives

$$s_{n,J}^{-2} \sum_{i=1}^{n} \mathbb{E}\left[V_{i,J}^2 \mathbf{1}(V_{i,J} \geq v)\right] \lesssim \xi_m^2 s_{n,J}^{-2} \mathbb{E}\left[\mathbf{1}(V_{1,J} \geq v)\right] \lesssim \xi_m^2 s_{n,J}^{-2} v^{-2} n^{-1} s_{n,J}^2,$$

where $\xi_m^2 n^{-1} = o(1)$ by (B1'). Thus the Lindeberg condition holds and it follows that

$$\frac{\sum_{i=1}^{n} V_{i,J}}{s_{n,J}} \xrightarrow{d} \mathcal{N}(0, 1).$$

Finally, it remains to prove (A.20) and (A.21). Begin by observing that

$$\begin{aligned}
&\left\|\mathbb{E}\left[(V - A(\tau)\widetilde{\mathbf{Z}}(W))\left(F_{Y|X}(Q(X; \tau)|X) - F_{Y|X}(\boldsymbol{\alpha}(\tau)^\top V + h_n^\dagger(W, \tau)|X)\right)\right]\right\| \\
&\lesssim \left\|\mathbb{E}\left[(V - A(\tau)\widetilde{\mathbf{Z}}(W))f_{Y|X}(Q(X_i; \tau)|X)(h(W, \tau) - \boldsymbol{\beta}_n^\dagger(\tau)^\top \widetilde{\mathbf{Z}}(W))\right]\right\| \\
&\quad + \left\|\mathbb{E}\left[(V - A(\tau)\widetilde{\mathbf{Z}}(W))f'_{Y|X}(\bar{Q}^\dagger(X_i, \tau)|X)(h(W, \tau) - \boldsymbol{\beta}_n^\dagger(\tau)^\top \widetilde{\mathbf{Z}}(W))^2\right]\right\| \\
&\lesssim \left\|\mathbb{E}\left[(h_{VW}(W, \tau) - A(\tau)\widetilde{\mathbf{Z}}(W))f_{Y|X}(Q(X_i; \tau)|X)(h(W, \tau) - \boldsymbol{\beta}_n^\dagger(\tau)^\top \widetilde{\mathbf{Z}}(W))\right]\right\| + c_n^{\dagger 2} \\
&\lesssim c_n^\dagger \mathbb{E}\left[\left\|h_{VW}(W, \tau) - A(\tau)\widetilde{\mathbf{Z}}(W)\right\|f_{Y|X}(Q(X_i; \tau)|X)\right] + c_n^{\dagger 2} \\
&= O(c_n^\dagger \lambda_n + c_n^{\dagger 2}),
\end{aligned} \quad (A.26)$$

where the first inequality is an application of Taylor expansion, with $\bar{Q}^\dagger(X_i, \tau)$ lying on the line segment of $Q(X; \tau)$ and $\boldsymbol{\alpha}(\tau)^\top V_i + h_n^\dagger(W_i, \tau)$; the second inequality is the result of the orthogonality condition (3.5), Condition (A2), Conditions (C1)-(C2); the third inequality follows from Conditions (A2), (C1), and the last line follows from condition (C1) and the Hölder inequality. For a proof of (A.20) observe that by (5.10) $\mathbf{e}_j^\top J_m(\tau)^{-1} \mathbf{Z} = \mathbf{e}_j^\top M_1(\tau)^{-1}(V - A(\tau)\widetilde{\mathbf{Z}}(W))$ for $j = 1, ..., k$. Thus we obtain from Remark 5.5

$$\begin{aligned}
&\left|\mathbb{E}[\mathbf{e}_j^\top J_m(\tau)^{-1} \mathbf{Z}(\mathbf{1}\{Y \leq Q(X; \tau)\} - \mathbf{1}\{Y \leq \boldsymbol{\alpha}(\tau)^\top V + h_n^\dagger(W, \tau)\})]\right| \\
&= \left|\mathbf{e}_j^\top M_1(\tau)^{-1} \mathbb{E}\left[(V - A(\tau)\widetilde{\mathbf{Z}}(W))\left(F_{Y|X}(Q(X; \tau)|X) - F_{Y|X}(\boldsymbol{\alpha}(\tau)^\top V + h_n^\dagger(W, \tau)|X)\right)\right]\right| \\
&= O(c_n^\dagger \lambda_n + c_n^{\dagger 2}),
\end{aligned}$$

To prove (A.21), without loss of generality, let $\|\mathbf{w}_n\| = 1$. We note that by (5.10)

$$(\mathbf{0}_k^\top, \mathbf{w}_n^\top) J_m(\tau)^{-1} \mathbf{Z}_i = -\mathbf{w}_n^\top A(\tau)^\top M_1(\tau)^{-1}(V_i - A(\tau)\widetilde{\mathbf{Z}}(W_i)) + \mathbf{w}_n^\top M_2(\tau)^{-1} \widetilde{\mathbf{Z}}(W_i).$$

From (A.26), (5.11) in Remark 5.5 and Lemma A.2 we obtain

$$\left| \mathbb{E}\left[ \mathbf{w}_n^\top A(\tau)^\top M_1(\tau)^{-1}(V_i - A(\tau)\widetilde{\mathbf{Z}}(W_i)) \big( \mathbf{1}\{Y_i \leq Q(X_i; \tau)\} - \mathbf{1}\{Y_i \leq \boldsymbol{\alpha}(\tau)^\top V_i + h_n^\dagger(W_i, \tau)\} \big) \right] \right|$$
$$= O(c_n^\dagger \lambda_n + c_n^{\dagger 2}). \tag{A.27}$$

Moreover,

$$\left| \mathbb{E}\left[ \mathbf{w}_n^\top M_2(\tau)^{-1} \widetilde{\mathbf{Z}}(W_i) \big( \mathbf{1}\{Y_i \leq Q(X_i; \tau)\} - \mathbf{1}\{Y_i \leq \boldsymbol{\alpha}(\tau)^\top V_i + h_n^\dagger(W_i, \tau)\} \big) \right] \right|$$
$$\lesssim \left| \mathbf{w}_n^\top M_2(\tau)^{-1} \mathbb{E}\left[ \widetilde{\mathbf{Z}}(W_i) f_{Y|X}(Q(X_i;\tau)|X) \big( h(W, \tau) - \boldsymbol{\beta}_n^\dagger(\tau)^\top \widetilde{\mathbf{Z}}(W) \big) \right] \right|$$
$$\quad + \left| \mathbf{w}_n^\top M_2(\tau)^{-1} \mathbb{E}\left[ \widetilde{\mathbf{Z}}(W_i) f'_{Y|X}(\bar{Q}^\dagger(X_i,\tau)|X) \big( h(W, \tau) - \boldsymbol{\beta}_n^\dagger(\tau)^\top \widetilde{\mathbf{Z}}(W) \big)^2 \} \right] \right|$$
$$= \left| \mathbb{E}\left[ \mathbf{w}_n^\top M_2(\tau)^{-1} \widetilde{\mathbf{Z}}(W_i) f'_{Y|X}(\bar{Q}^\dagger(X_i,\tau)|X) \big( h(W, \tau) - \boldsymbol{\beta}_n^\dagger(\tau)^\top \widetilde{\mathbf{Z}}(W) \big)^2 \} \right] \right|$$
$$\lesssim \mathbb{E}\left[ |\mathbf{w}_n^\top M_2(\tau_2)^{-1} \widetilde{\mathbf{Z}}(W_i)| \right] c_n^{\dagger 2} = o(\|\mathbf{w}_n\| c_n^{\dagger 2}), \tag{A.28}$$

where the first inequality follows from the Taylor expansion with $\bar{Q}^\dagger(X_i, \tau)$ lying on the line segment of $Q(X; \tau)$ and $\boldsymbol{\alpha}(\tau)^\top V_i + h_n^\dagger(W_i, \tau)$; the second equality follows from the first order condition of (3.3), and the last line follows Conditions (A2), (C1) and the conditions of the theorem. Combining (A.27) and (A.28), and (B1') we obtain (A.21).

**Lemma A.2.** *Under Assumptions (A1)-(A3), (C2) and $\sup_{\tau \in \mathcal{T}} \mathbb{E}[|\widetilde{\mathbf{Z}}(W)^\top M_2(\tau)^{-1} \mathbf{w}_n|] = o(\|\mathbf{w}_n\|)$, we have*

$$\sup_\tau \|A(\tau) \mathbf{w}_n\|^2 = o(\|\mathbf{w}_n\|) \tag{A.29}$$

**Proof for Lemma A.2.** By the first order condition for obtaining $A(\tau)$,

$$A(\tau) = \mathbb{E}[h_{VW}(W; \tau) \widetilde{\mathbf{Z}}(W)^\top f_{Y|X}(Q(X;\tau)|X)] M_2(\tau)^{-1}.$$

By the orthogonality condition (3.5),

$$\|A(\tau) \mathbf{w}_n\| = \left\| \mathbb{E}[h_{VW}(W;\tau) \widetilde{\mathbf{Z}}(W)^\top f_{Y|X}(Q(X;\tau)|X)] M_2(\tau)^{-1} \mathbf{w}_n \right\|$$
$$= \left\| \mathbb{E}[(h_{VW}(W;\tau) - V + V) \widetilde{\mathbf{Z}}(W)^\top f_{Y|X}(Q(X;\tau)|X)] M_2(\tau)^{-1} \mathbf{w}_n \right\|$$
$$= \left\| \mathbb{E}[V \widetilde{\mathbf{Z}}(W)^\top f_{Y|X}(Q(X;\tau)|X)] M_2(\tau)^{-1} \mathbf{w}_n \right\|.$$

By the assumption that at fixed $j$, $|V_j| \leq C$, the uniform boundedness of the conditional density in (A2), and the hypothesis $\sup_{\tau \in \mathcal{T}} \mathbb{E}[|\widetilde{\mathbf{Z}}(W)^\top M_2(\tau)^{-1} \mathbf{w}_n|] = o(\|\mathbf{w}_n\|)$,

$$\sup_{\tau \in \mathcal{T}} \mathbb{E}[|V_j \widetilde{\mathbf{Z}}(W)^\top M_2(\tau)^{-1} \mathbf{w}_n | f_{Y|X}(Q(X;\tau)|X)] \leq \bar{f} C \sup_{\tau \in \mathcal{T}} \mathbb{E}[|\widetilde{\mathbf{Z}}(W)^\top M_2(\tau)^{-1} \mathbf{w}_n|] = o(\|\mathbf{w}_n\|). \tag{A.30}$$

This completes the proof of (A.29) by noting that $\sup_{\tau \in \mathcal{T}} \|M_2(\tau)^{-1}\| = O(1)$. □

### *A.4. Asymptotic Tightness of Quantile Process*

In this section we establish the asymptotic tightness of the process $n^{1/2} \mathbf{u}_n^\top \boldsymbol{U}_n(\tau)$ in $\ell^\infty(\mathcal{T})$ with $\mathbf{u}_n \in \mathbb{R}^m$ being an arbitrary vector, where

$$\boldsymbol{U}_n(\tau) = n^{-1} J_m^{-1}(\tau) \sum_{i=1}^n \mathbf{Z}_i \big( \mathbf{1}\{Y_i \leq Q(X_i; \tau)\} - \tau \big). \tag{A.31}$$

Note that the results obtained in this section, in particular Lemma A.3, apply to any series expansion $\mathbf{Z} = \mathbf{Z}(X_i)$ satisfying Assumptions (A1).

The following definition is only used in this subsection: For any non-decreasing, convex function $\Psi : \mathbb{R}^+ \to \mathbb{R}^+$ with $\Phi(0) = 0$, the *Orlicz norm* of a real-valued random variable $Z$ is defined as (see e.g. Chapter 2.2 of van der Vaart and Wellner (1996))

$$\|Z\|_\Psi = \inf \{C > 0 : \mathbb{E}\Phi(|Z|/C) \leq 1\}. \tag{A.32}$$

**Lemma A.3** (Asymptotic Equicontinuity of Quantile Process). *Under (A1)-(A3) and $\xi_m^2(\log n)^2 = o(n)$, we have for any $\varepsilon > 0$ and vector $\mathbf{u}_n \in \mathbb{R}^m$,*

$$\lim_{\delta \to 0} \limsup_{n \to \infty} P\Big(\|\mathbf{u}_n\|^{-1} n^{1/2} \sup_{\tau_1, \tau_2 \in \mathcal{T}, |\tau_1 - \tau_2| \leq \delta} \Big|\mathbf{u}_n^\top \boldsymbol{U}_n(\tau_1) - \mathbf{u}_n^\top \boldsymbol{U}_n(\tau_2)\Big| > \varepsilon\Big) = 0, \tag{A.33}$$

*where $\boldsymbol{U}_n(\tau)$ is defined in (A.31).*

**Proof of Lemma A.3.** Without loss of generality, we will assume that $\mathbf{u}_n$ is a sequence of vectors with $\|\mathbf{u}_n\| = 1$, which can always be achieved by rescaling. Define

$$\mathbb{G}_n(\tau) := n^{1/2} \mathbf{u}_n^\top \boldsymbol{U}_n(\tau).$$

Consider the decomposition

$$\mathbf{u}_n^\top \boldsymbol{U}_n(\tau_1) - \mathbf{u}_n^\top \boldsymbol{U}_n(\tau_2) = n^{-1} \mathbf{u}_n^\top (J_m^{-1}(\tau_1) - J_m^{-1}(\tau_2)) \sum_i \mathbf{Z}_i (\mathbf{1}\{Y_i \leq Q(X_i; \tau_1)\} - \tau_1)$$
$$+ n^{-1} \mathbf{u}_n^\top J_m^{-1}(\tau_2) \sum_i \mathbf{Z}_i \Delta_i(\tau_1, \tau_2),$$

where

$$\Delta_i(\tau_1, \tau_2) := \mathbf{1}\{Y_i \leq Q(X_i; \tau_1)\} - \tau_1 - (\mathbf{1}\{Y_i \leq Q(X_i; \tau_2)\} - \tau_2).$$

Note that for any $L \geq 2$,

$$\mathbb{E}\Big[\big|\mathbf{u}_n^\top J_m^{-1}(\tau_2) \mathbf{Z}_i \Delta_i(\tau_1, \tau_2)\big|^L\Big] \lesssim \xi_m^{L-2} \mathbb{E}\Big[\big|\mathbf{u}_n^\top J_m^{-1}(\tau_2) \mathbf{Z}_i \Delta_i(\tau_1, \tau_2)\big|^2\Big]$$
$$= \xi_m^{L-2} \mathbf{u}_n^\top J_m^{-1}(\tau_2) \mathbb{E}\big[\mathbf{Z}_i \mathbf{Z}_i^\top \Delta_i(\tau_1, \tau_2)^2\big] J_m^{-1}(\tau_2) \mathbf{u}_n$$
$$\lesssim \xi_m^{L-2} |\tau_1 - \tau_2|. \tag{A.34}$$

By the Lipschitz continuity of $\tau \mapsto J_m^{-1}(\tau)$ (cf. Lemma 13 of Belloni et al. (2016)) and positive definiteness of $J_m^{-1}(\tau)$, we have

$$\|J_m^{-1}(\tau_1) - J_m^{-1}(\tau_2)\| = \|J_m^{-1}(\tau_2)\{J_m(\tau_1) - J_m(\tau_2)\} J_m^{-1}(\tau_1)\|$$
$$\leq \frac{\bar{f}'}{f_{\min}} |\tau_1 - \tau_2| \Big(\inf_{\tau \in \mathcal{T}} \lambda_{\min}(J_m(\tau))\Big)^{-2} \lambda_{\max}(\mathbb{E}[\mathbf{Z}\mathbf{Z}^\top]),$$

where $\|\cdot\|$ denotes the operator norm of a matrix. Thus, we have for $L \geq 2$,

$$\mathbb{E}\Big[\big|\mathbf{u}_n^\top \{J_m^{-1}(\tau_1) - J_m^{-1}(\tau_2)\} \mathbf{Z}_i (\mathbf{1}(Y_i \leq Q(X_i; \tau)) - \tau)\big|^L\Big]$$
$$\lesssim \xi_m^{L-2} \mathbb{E}\Big[\big|\mathbf{u}_n^\top \{J_m^{-1}(\tau_1) - J_m^{-1}(\tau_2)\} \mathbf{Z}_i (\mathbf{1}(Y_i \leq Q(X_i; \tau)) - \tau)\big|^2\Big]$$
$$\leq \xi_m^{L-2} \mathbb{E}\Big[\big|\mathbf{u}_n^\top \{J_m^{-1}(\tau_1) - J_m^{-1}(\tau_2)\} \mathbf{Z}_i\big|^2\Big]$$
$$= \xi_m^{L-2} \mathbf{u}_n^\top \{J_m^{-1}(\tau_1) - J_m^{-1}(\tau_2)\} \mathbb{E}[\mathbf{Z}_i \mathbf{Z}_i^\top] \{J_m^{-1}(\tau_1) - J_m^{-1}(\tau_2)\} \mathbf{u}_n$$
$$\leq \xi_m^{L-2} \big\|\{J_m^{-1}(\tau_1) - J_m^{-1}(\tau_2)\} \mathbb{E}[\mathbf{Z}_i \mathbf{Z}_i^\top] \{J_m^{-1}(\tau_1) - J_m^{-1}(\tau_2)\}\big\|$$
$$\lesssim \xi_m^{L-2} |\tau_1 - \tau_2|^2. \tag{A.35}$$

To simplify notations, define
$$\widetilde{V}_{n,i}(\tau_1, \tau_2) := \mathbf{u}_n^\top \{J_m^{-1}(\tau_1) - J_m^{-1}(\tau_2)\}\mathbf{Z}_i(\mathbf{1}(Y_i \leq Q(X_i; \tau)) - \tau) + \mathbf{u}_n^\top J_m^{-1}(\tau_2)\mathbf{Z}_i\Delta_i(\tau_1, \tau_2).$$

Combining the bounds (A.34) and (A.35) yields

$$\mathbb{E}\big[|\widetilde{V}_{n,i}(\tau_1, \tau_2)|^{2L}\big]^{1/2L}$$
$$\leq \mathbb{E}\big[|\mathbf{u}_n^\top J_m^{-1}(\tau_2)\mathbf{Z}_i\Delta_i(\tau_1,\tau_2)|^{2L}\big]^{1/2L} + \mathbb{E}\big[|\mathbf{u}_n^\top\{J_m^{-1}(\tau_1) - J_m^{-1}(\tau_2)\}\mathbf{Z}_i(\mathbf{1}(Y_i \leq Q(X_i;\tau)) - \tau)|^{2L}\big]^{1/2L}$$
$$\lesssim \big(\xi_m^{2(L-1)}|\tau_1 - \tau_2|\big)^{1/2L}. \tag{A.36}$$

Note that (A.36) holds for all positive integers $L \geq 1$. By the fact that $\mathbb{G}_n(\tau_1) - \mathbb{G}_n(\tau_2) = n^{-1/2}\sum_{i=1}^n \widetilde{V}_{n,i}(\tau_1, \tau_2)$ and $\mathbb{E}\widetilde{V}_{n,i}(\tau_1, \tau_2) = 0$, we obtain from (A.36) that

$$\mathbb{E}[|\mathbb{G}_n(\tau_1) - \mathbb{G}_n(\tau_2)|^{2L}]$$
$$= n^{-L}\mathbb{E}\bigg[\bigg(\sum_{i=1}^n \widetilde{V}_{n,i}(\tau_1, \tau_2)\bigg)^{2L}\bigg]$$
$$= n^{-L}\bigg(\sum_{i=1}^n \mathbb{E}\big[\widetilde{V}_{n,i}(\tau_1,\tau_2)^{2L}\big] + \sum_{l=1}^{L-1}\sum_{\substack{1\leq i_1, i_2 \leq n \\ i_1 \neq i_2}} \mathbb{E}\big[\widetilde{V}_{n,i_1}(\tau_1,\tau_2)^{2L-2l}\big]\mathbb{E}\big[\widetilde{V}_{n,i_2}(\tau_1,\tau_2)^{2l}\big]$$
$$+ \sum_{\substack{l_1+l_2<L \\ l_1=1, l_2=1}}\sum_{\substack{1\leq i_1, i_2, i_3 \leq n \\ i_1\neq i_2 \neq i_3}} \mathbb{E}\big[\widetilde{V}_{n,i_1}(\tau_1,\tau_2)^{2L-2(l_1+l_2)}\big]\mathbb{E}\big[\widetilde{V}_{n,i_2}(\tau_1,\tau_2)^{2l_1}\big]\mathbb{E}\big[\widetilde{V}_{n,i_3}(\tau_1,\tau_2)^{2l_2}\big]$$
$$+ \ldots + \sum_{\substack{1\leq i_1,\ldots,i_L\leq n \\ i_1\neq\ldots\neq i_L}}\prod_{j=1}^L \mathbb{E}\big[\widetilde{V}_{n,i_j}(\tau_1,\tau_2)^2\big]\bigg)$$
$$\leq C_L n^{-L}\bigg(n\xi_m^{2(L-1)}|\tau_1-\tau_2| + \binom{n}{2}\xi_m^{2(L-1-1)}|\tau_1-\tau_2|^2 + \binom{n}{3}\xi_m^{2(L-1-2)}|\tau_1-\tau_2|^3 + \ldots + \binom{n}{L}|\tau_1-\tau_2|^L\bigg)$$
$$\lesssim \sum_{k=0}^{L-1}\frac{\xi_m^{2(L-k-1)}}{n^{(L-k-1)}}|\tau_1 - \tau_2|^{k+1}.$$

In particular we obtain for $|\tau_1 - \tau_2| \geq \xi_m^2/n$,
$$\mathbb{E}[|\mathbb{G}_n(\tau_1) - \mathbb{G}_n(\tau_2)|^{2L}] \lesssim |\tau_1 - \tau_2|^L. \tag{A.37}$$

For $\Psi(z) = z^{2L}$, the above equation implies that the Orlicz norm (defined in (A.32)) of $\mathbb{G}_n(\tau_1) - \mathbb{G}(\tau_2)$ satisfies
$$\|\mathbb{G}_n(\tau_1) - \mathbb{G}(\tau_2)\|_\Psi \lesssim |\tau_1 - \tau_2|^{1/2}.$$

Let $d(\tau, \tau') = \sqrt{|\tau - \tau'|}$, which is a metric on $\mathcal{T}$. The packing number $D(\epsilon, d)$ of $\mathcal{T}$ with respect to $d$ satisfies $D(\epsilon, d) \lesssim 1/\epsilon^2$. Let $\bar{\omega}_n = 2\xi_m/\sqrt{n} \to 0$ as $n \to \infty$. We have
$$\int_{\bar{\omega}_n/2}^\omega \Psi^{-1}(D(\epsilon, d))d\epsilon \lesssim \int_{\bar{\omega}_n/2}^\omega \epsilon^{-1/L}d\epsilon = \frac{\omega^{1-L^{-1}}}{1-L^{-1}} - \frac{(\bar{\omega}_n/2)^{1-L^{-1}}}{1-L^{-1}}. \tag{A.38}$$

For $\omega > 0$,
$$\Psi^{-1}(D^2(\omega, d)) \lesssim \Psi^{-1}\bigg(\frac{1}{\omega^4}\bigg) = \omega^{-2/L}. \tag{A.39}$$

Therefore, applying Lemma S.2.1 yields for any $\delta > 0$
$$\sup_{|\tau_1-\tau_2|\leq\delta}|\mathbb{G}_n(\tau_1) - \mathbb{G}_n(\tau_2)| = \sup_{|\tau_1-\tau_2|^{1/2}\leq\delta^{1/2}}|\mathbb{G}_n(\tau_1) - \mathbb{G}_n(\tau_2)|$$
$$\leq S_{1,n}(\delta) + 2\sup_{|\tau'-\tau|^{1/2}\leq\bar{\omega}_n, \tau\in\widetilde{\mathcal{T}}}|\mathbb{G}_n(\tau') - \mathbb{G}_n(\tau)|, \tag{A.40}$$

where $\widetilde{\mathcal{T}} \subset \mathcal{T}$ has at most $D(\bar{\omega}_n, d) \lesssim \bar{\omega}_n^{-2}$ points and $S_{1,n}(\delta)$ is a random variable that satisfies

$$P(|S_{1,n}(\delta)| > z) \leq \left( z \left[ 8K \left( \int_{\bar{\omega}_n/2}^{\omega} \Psi^{-1}(D(\epsilon, d)) d\epsilon + (\delta^{1/2} + 2\bar{\omega}_n) \Psi^{-1}(D^2(\omega, d)) \right) \right]^{-1} \right)^{-2L}$$

$$\lesssim \left( \frac{\frac{\omega^{1-L^{-1}}}{1-L^{-1}} - \frac{(\bar{\omega}_n/2)^{1-L^{-1}}}{1-L^{-1}} + (\delta^{1/2} + 2\bar{\omega}_n)\omega^{-2/L}}{z} \right)^{2L}. \quad (A.41)$$

for a constant $K > 0$. Let $\omega = \delta$ and $L = 6$. As $n \to \infty$, $\omega > \bar{\omega}_n$. We obtain $\lim_{\delta \to 0} \limsup_{n \to \infty} P(|S_{1,n}(\delta)| > z) = 0$ for any $z > 0$.

To bound the remaining term in (A.40), observe that

$$\sup_{d(\tau, \tau') \leq \bar{\omega}_n, \tau \in \widetilde{\mathcal{T}}} |\mathbb{G}_n(\tau') - \mathbb{G}_n(\tau)| = \sup_{|\tau - \tau'| \leq \bar{\omega}_n^2, \tau \in \widetilde{\mathcal{T}}} |\mathbb{G}_n(\tau') - \mathbb{G}_n(\tau)| \leq \sup_{|\tau - \tau'| \leq \bar{\omega}_n^2, \tau, \tau' \in \mathcal{T}} |\mathbb{G}_n(\tau') - \mathbb{G}_n(\tau)|.$$

Now by Lemma A.4 we have

$$P \left( \sup_{|\tau' - \tau| \leq \bar{\omega}_n^2, \tau, \tau' \in \mathcal{T}} |\mathbb{G}_n(\tau') - \mathbb{G}_n(\tau)| > r_n(\kappa_n) \right) < e^{-\kappa_n},$$

where

$$r_n(\kappa_n) := C \left[ \bar{\omega}_n \log^{1/2} \frac{\xi_m}{\bar{\omega}_n} + \frac{\xi_m}{\sqrt{n}} \log \frac{\xi_m}{\bar{\omega}_n} + \kappa_n^{1/2} \bar{\omega}_n + \frac{\xi_m}{\sqrt{n}} \kappa_n \right], \quad (A.42)$$

for a sufficiently large constant $C > 0$. Take $\kappa_n = \log n$. Since $\bar{\omega}_n (\log n)^{1/2} = 2\xi_m (\log n)^{1/2}/\sqrt{n} = o(1)$ and $\xi_m \log(n)/\sqrt{n} = o(1)$ by assumption, it follows that $r_n(\log n) \to 0$. Therefore, we conclude from Lemma A.4 that

$$\sup_{d(\tau, \tau') \leq \bar{\omega}_n, \tau \in \widetilde{\mathcal{T}}} |\mathbb{G}_n(\tau') - \mathbb{G}_n(\tau)| = o_P(1). \quad (A.43)$$

Applying bounds (A.41) and (A.43) to (A.40) verifies the asymptotic equicontinuity

$$\lim_{\delta \to 0} \limsup_{n \to \infty} P \left( \sup_{|\tau_1 - \tau_2| \leq \delta} |\mathbb{G}_n(\tau_1) - \mathbb{G}_n(\tau_2)| > z \right) = 0$$

for all $z > 0$. $\square$

The following result is applied in the proof of Lemma A.3.

**Lemma A.4.** *Under (A1)-(A3), we have for any $\kappa_n > 0$, $1/n \ll \delta < 1$,*

$$P \left( \sup_{0 \leq h \leq \delta} \sup_{\tau \in [\epsilon, 1-\epsilon-h]} |\mathbf{u}_n^\top \mathbf{U}_n(\tau + h) - \mathbf{u}_n^\top \mathbf{U}_n(\tau)| \geq C r_n(\delta, \kappa_n) \right) \leq 3 e^{-\kappa_n}, \quad (A.44)$$

*where $\kappa_n > 0$, $\mathbf{U}_n(\tau)$ is defined in (A.31) and $\mathbf{u}_n \in \mathbb{R}^m$ is arbitrary, and*

$$r_n(\delta, \kappa_n) = \|\mathbf{u}_n\| \left( \frac{\delta}{n} \log \frac{\xi_m}{\sqrt{\delta}} \right)^{1/2} + \frac{\|\mathbf{u}_n\| \xi_m}{n} \log \frac{\xi_m}{\sqrt{\delta}} + \|\mathbf{u}_n\| \left( \frac{\kappa_n \delta}{n} \right)^{1/2} + \frac{\|\mathbf{u}_n\| \xi_m}{n} \kappa_n.$$

To prove Lemma A.4, we need to establish some preliminary results. For any fixed vector $\mathbf{u} \in \mathbb{R}^m$ and $\delta > 0$, define the function classes

$$\mathcal{G}_3(\mathbf{u}) := \left\{ (\mathbf{Z}, Y) \mapsto \mathbf{u}^\top J_m(\tau)^{-1} \mathbf{Z} \mathbf{1}\{\|\mathbf{Z}\| \leq \xi_m\} \big| \tau \in \mathcal{T} \right\},$$

$$\mathcal{G}_4 := \left\{ (X, Y) \mapsto \mathbf{1}\{Y_i \leq Q(X; \tau)\} - \tau \big| \tau \in \mathcal{T} \right\},$$

$$\mathcal{G}_6(\mathbf{u}, \delta) := \left\{ (\mathbf{Z}, Y) \mapsto \mathbf{u}^\top \{J_m(\tau_1)^{-1} - J_m(\tau_2)^{-1}\} \mathbf{Z} \mathbf{1}\{\|\mathbf{Z}\| \leq \xi_m\} \big| \tau_1, \tau_2 \in \mathcal{T}, |\tau_1 - \tau_2| \leq \delta \right\},$$

$$\mathcal{G}_7(\delta) := \left\{ (X, Y) \mapsto \mathbf{1}\{Y_i \leq Q(X, \tau_1)\} - \mathbf{1}\{Y_i \leq Q(X, \tau_2)\} - (\tau_1 - \tau_2) \big| \tau_1, \tau_2 \in \mathcal{T}, |\tau_1 - \tau_2| \leq \delta \right\}.$$

Denote $G_3$, $G_6$ and $G_7$ as the envelope functions of $\mathcal{G}_3$, $\mathcal{G}_6$ and $\mathcal{C}_7$, respectively. The following covering number results will be shown in Section S.2.2: for any probability measure $Q$,

$$N(\epsilon \|G_3\|_{L_2(Q)}, \mathcal{G}_3(\mathbf{u}), L_2(Q)) \leq \frac{C_0}{\epsilon}, \tag{A.45}$$

$$N(\epsilon \|G_6\|_{L_2(Q)}, \mathcal{G}_6(\mathbf{u}, \delta), L_2(Q)) \leq 2\left(\frac{C_0}{\epsilon}\right)^2, \tag{A.46}$$

$$N(\epsilon \|G_7\|_{L_2(Q)}, \mathcal{G}_7(\delta), L_2(Q)) \leq \left(\frac{A_7}{\epsilon}\right)^2, \tag{A.47}$$

where $C_0 := \frac{\bar{f}'}{f_{\min}} \frac{\lambda_{\max}(\mathbb{E}[\mathbf{Z}\mathbf{Z}^\top])}{\inf_{\tau \in \mathcal{T}} \lambda_{\min}(J_m(\tau))} < \infty$ given Assumptions (A1)-(A3), and $A_7 > 0$ is a constant. Also, $\mathcal{G}_4$ has VC index 2 according to Lemma S.2.4.

**Proof of Lemma A.4.** Observe the decomposition

$$\mathbf{u}_n^\top \boldsymbol{U}_n(\tau_1) - \mathbf{u}_n^\top \boldsymbol{U}_n(\tau_2) = I_1(\tau_1, \tau_2) + I_2(\tau_1, \tau_2),$$

where

$$I_1(\tau_1, \tau_2) := n^{-1} \mathbf{u}_n^\top \{J_m(\tau_1)^{-1} - J_m(\tau_2)^{-1}\} \sum_{i=1}^n \mathbf{Z}_i \big(\mathbf{1}\{Y_i \leq Q(X_i, \tau_1)\} - \tau_1\big),$$

$$I_2(\tau_1, \tau_2) := n^{-1} \mathbf{u}_n^\top J_m(\tau_2)^{-1} \sum_{i=1}^n \mathbf{Z}_i \big(\mathbf{1}\{Y_i \leq Q(X_i, \tau_1)\} - \mathbf{1}\{Y_i \leq Q(X_i, \tau_2)\} - (\tau_1 - \tau_2)\big).$$

**Step 1: bounding $I_1(\tau_1, \tau_2)$.**

Note that $\sup_{\tau_1, \tau_2 \in \mathcal{T}, |\tau_1 - \tau_2| < \delta} |I_1(\tau_1, \tau_2)| \leq \|\mathbb{P}_n - P\|_{\mathcal{G}_6(\mathbf{u}_n, \delta) \cdot \mathcal{G}_4}$, where

$$\mathcal{G}_6(\mathbf{u}_n, \delta) \cdot \mathcal{G}_4 = \Big\{ \mathbf{u}_n^\top \{J_m(\tau_1)^{-1} - J_m(\tau_2)^{-1}\} \mathbf{Z}_i \big(\mathbf{1}\{Y_i \leq Q(X_i, \tau_3)\} - \tau_3\big) \Big| \tau_1, \tau_2, \tau_3 \in \mathcal{T}, |\tau_1 - \tau_2| \leq \delta \Big\}.$$

Theorem 2.6.7 of van der Vaart and Wellner (1996) and Part 1 of Lemma S.2.4 give

$$N(\epsilon \|G_4\|_{L_2(\mathbb{P}_n)}, \mathcal{G}_4, L_2(\mathbb{P}_n)) \leq \frac{A_4}{\epsilon}$$

where the envelope for $\mathcal{G}_4$ is $G_4 = 2$ and $A_4$ is a universal constant. Part 2 of Lemma S.2.4 and Part 2 of Lemma S.2.2 imply that

$$N(\epsilon \|G_6 G_4\|_{L_2(\mathbb{P}_n)}, \mathcal{G}_6(\mathbf{u}_n, \delta) \cdot \mathcal{G}_4, L_2(\mathbb{P}_n)) \leq \frac{2A_4}{\epsilon}\left(\frac{2C_0}{\epsilon}\right)^2 \leq \left(\frac{2A_4^{1/3} C_0^{2/3}}{\epsilon}\right)^3, \tag{A.48}$$

where $2A_4^{1/3} C_0^{2/3} \leq C$ for a large enough universal constant $C > 0$. To bound $\sup_{f \in \mathcal{G}_6(\mathbf{u}_n, \delta) \cdot \mathcal{G}_4} \|f\|_{L^2(P)}$, note that

$$\mathbb{E}\big(\mathbf{u}_n^\top \{J_m(\tau_1)^{-1} - J_m(\tau_2)^{-1}\} \mathbf{Z}_i \big(\mathbf{1}\{Y_i \leq Q(X_i, \tau_3)\} - \tau_3\big)\big)^2$$
$$\leq 4\|\mathbf{u}_n\|^2 \lambda_{\max}(\mathbb{E}[\mathbf{Z}\mathbf{Z}^\top]) [\inf_{\tau \in \mathcal{T}} \lambda_{\min}(J_m(\tau))]^{-2} C_0^2 \delta^2 \leq C\|\mathbf{u}_n\|^2 \delta^2,$$

for a large enough constant $C$. In addition, an upper bound for the functions in $\mathcal{G}_6(\mathbf{u}_n, \delta) \cdot \mathcal{G}_4$ is

$$2\xi_m \|\mathbf{u}_n\| [\inf_{\tau \in \mathcal{T}} \lambda_{\min}(J_m(\tau))]^{-1} C_0 \delta \leq C\delta \xi_m \|\mathbf{u}_n\|,$$

and we can take this upper bound as envelope.

Applying the bounds (S.2.2) and (S.2.3) and taking into account (A.48), for any $\mathbf{u}_n$ and $\delta > 0$,

$$\mathbb{E}\|\mathbb{P}_n - P\|_{\mathcal{G}_6(\mathbf{u}_n, \delta) \cdot \mathcal{G}_4} \leq c_1 \Big[\|\mathbf{u}_n\| \delta \Big(\frac{\log(\xi_m)}{n}\Big)^{1/2} + \frac{\delta \|\mathbf{u}_n\| \xi_m}{n} \log(\xi_m)\Big]. \tag{A.49}$$

Finally, for any $\kappa_n > 0$, let

$$r_{n,1}(\delta, \kappa_n) = \widetilde{C}\bigg[\|\mathbf{u}_n\|\delta\Big(\frac{1}{n}\log(\xi_m)\Big)^{1/2} + \frac{\delta\|\mathbf{u}_n\|\xi_m}{n}\log(\xi_m) + \Big(\frac{\kappa_n}{n}\Big)^{1/2}\|\mathbf{u}_n\|\delta + \frac{\delta\|\mathbf{u}_n\|\xi_m}{n}\kappa_n\bigg]$$

for a sufficiently large constant $\widetilde{C} > 0$. From this, we obtain

$$P\bigg\{\sup_{\tau_1,\tau_2 \in \mathcal{T}, |\tau_1-\tau_2|<\delta} |I_1(\tau_1,\tau_2)| \geq r_{n,1}(\delta, \kappa_n)\bigg\} \leq P\{\|\mathbb{P}_n - P\|_{\mathcal{G}_6(\mathbf{u}_n,\delta)\cdot\mathcal{G}_4} \geq r_{n,1}(\delta, \kappa_n)\} \leq e^{-\kappa_n}.$$

**Step 2: bounding $I_2(\tau_1, \tau_2)$.**

Note that $\sup_{\tau_1,\tau_2\in\mathcal{T},|\tau_1-\tau_2|<\delta} |I_2(\tau_1,\tau_2)| \leq \|\mathbb{P}_n - P\|_{\mathcal{G}_3(\mathbf{u}_n)\cdot\mathcal{G}_7(\delta)}$, where

$$\mathcal{G}_3(\mathbf{u}_n)\cdot\mathcal{G}_7(\delta) = \Big\{\mathbf{u}_n^\top J_m(\tau_3)^{-1}\mathbf{Z}_i\big(\mathbf{1}\{Y_i \leq Q(X_i,\tau_1)\} - \mathbf{1}\{Y_i \leq Q(X_i,\tau_2)\} - (\tau_1-\tau_2)\big)\Big|$$
$$\tau_1, \tau_2, \tau_3 \in \mathcal{T}, |\tau_1-\tau_2| \leq \delta\Big\}.$$

Lemma S.2.3, Part 2 of Lemma S.2.4 and Part 2 of Lemma S.2.2 imply that

$$N(\epsilon\|G_3 G_7\|_{L_2(\mathbb{P}_n)}, \mathcal{G}_3(\mathbf{u}_n)\cdot\mathcal{G}_7(\delta), L_2(\mathbb{P}_n)) \leq \frac{2C_0}{\epsilon}\Big(\frac{2A_7}{\epsilon}\Big)^2 \leq \Big(\frac{2C_0^{1/3} A_7^{2/3}}{\epsilon}\Big)^3, \tag{A.50}$$

where $2C_0^{1/3} A_7^{2/3} \leq C$ for a large enough constant $C > 0$. To bound $\sup_{f\in\mathcal{G}_3(\mathbf{u}_n)\cdot\mathcal{G}_7(\delta)} \|f\|_{L^2(P)}$ note that

$$\mathbb{E}\big(\mathbf{u}_n^\top J_m(\tau_3)^{-1}\mathbf{Z}_i\big(\mathbf{1}\{Y_i \leq Q(X_i,\tau_1)\} - \mathbf{1}\{Y_i \leq Q(X_i,\tau_2)\} - (\tau_1-\tau_2)\big)\big)^2$$
$$\leq 3\|\mathbf{u}_n\|^2 \lambda_{\max}(\mathbb{E}[\mathbf{Z}\mathbf{Z}^\top])[\inf_{\tau\in\mathcal{T}}\lambda_{\min}(J_m(\tau))]^{-2}\delta \leq C\|\mathbf{u}_n\|^2\delta.$$

Moreover

$$\sup_{f\in\mathcal{G}_3(\mathbf{u}_n)\cdot\mathcal{G}_7(\delta)} \sup \|f\|_\infty \leq 2\|\mathbf{u}_n\|\xi_m[\inf_{\tau\in\mathcal{T}}\lambda_{\min}(J_m(\tau))]^{-1} \leq C\|\mathbf{u}_n\|\xi_m$$

for some constant $C$. Applying the bounds (S.2.2) and (S.2.3) and taking into account (A.50)

$$\mathbb{E}\|\mathbb{P}_n - P\|_{\mathcal{G}_3(\mathbf{u}_n)\cdot\mathcal{G}_7(\delta)} \leq c_1\bigg[\|\mathbf{u}_n\|\Big(\frac{\delta}{n}\log\frac{\xi_m}{\sqrt{\delta}}\Big)^{1/2} + \frac{\|\mathbf{u}_n\|\xi_m}{n}\log\frac{\xi_m}{\sqrt{\delta}}\bigg]. \tag{A.51}$$

For any $\kappa_n > 0$, let

$$r_{n,2}(\delta, \kappa_n) = C\bigg[\|\mathbf{u}_n\|\Big(\frac{\delta}{n}\log\frac{\xi_m}{\sqrt{\delta}}\Big)^{1/2} + \frac{\|\mathbf{u}_n\|\xi_m}{n}\log\frac{\xi_m}{\sqrt{\delta}} + \|\mathbf{u}_n\|\Big(\frac{\kappa_n\delta}{n}\Big)^{1/2} + \frac{\|\mathbf{u}_n\|\xi_m}{n}\kappa_n\bigg]$$

for a constant $C > 0$ sufficiently large, we obtain

$$P\bigg\{\sup_{\tau_1,\tau_2\in\mathcal{T},|\tau_1-\tau_2|<\delta} |I_2(\tau_1,\tau_2)| \geq r_{n,2}(\delta,\kappa_n)\bigg\} \leq P\{\|\mathbb{P}_n - P\|_{\mathcal{G}_3(\mathbf{u}_n)\cdot\mathcal{G}_7(\delta)} \geq r_{n,2}(\delta,\kappa_n)\} \leq e^{-\kappa_n}.$$

Finally, $r_{n,1}(\delta, \kappa_n) \leq r_{n,2}(\delta, \kappa_n)$ when $\delta < 1$. Hence, we conclude (A.44). $\square$

### A.5. Proof of Corollary 4.1

As the argument $x_0$ in $Q(x_0; \tau)$ and $F_{Y|X}(y|x_0)$ is fixed, simplify notations by writing $Q(x_0; \tau) = Q(\tau), \widehat{Q}(x_0; \tau) = \widehat{Q}(\tau)$ and $F_{Y|X}(y|x_0) = F(y), \widehat{F}_{Y|X}(y|x_0) = \widehat{F}(y)$ as functions of the single arguments in $\tau$ and $y$, respectively. From Theorems 2.4, 3.1 or Corollary 2.2, we have

$$a_n\big(\widehat{Q}(\cdot) - Q(\cdot)\big) \rightsquigarrow \mathbb{G}(\cdot) \text{ in } \ell^\infty([\tau_L, \tau_U]), \tag{A.52}$$

where $a_n$ and $\mathbb{G}$ depend on the model for $Q(x;\tau)$ and $\mathbb{G}$ has continuous sample paths almost surely. Next, note that for $y \in \mathcal{Y}$
$$a_n\big(\widehat{F}(y) - F(y)\big) = a_n\big(\Phi(\widehat{Q})(y) - \Phi(Q)(y)\big).$$
Finally, observe that $\Phi(f)(y) = \tau_L + (\tau_U - \tau_L)(\Phi^* \circ R)(f)(y)$ where $\Phi^*(f)(y) := \int_0^1 \mathbf{1}\{f(u) < y\}du$ and $R(f)(y) := f(\tau_L + y(\tau_U - \tau_L))$. The map $R : \ell^\infty((\tau_L, \tau_U)) \to \ell^\infty((0,1))$ is linear and continuous, hence compactly differentiable with derivative $R$. The map $\Phi^*$ is compactly differentiable tangentially to $\mathcal{C}(0,1)$ at any strictly increasing, differentiable function $f_0$ and the derivative of $\Phi$ at $f_0$ is given by $d\Phi^*_{f_0}(h)(y) = -h(f_0^{-1}(y))/f_0'(f_0^{-1}(y))$ - see Corollary 1 in Chernozhukov et al. (2010). Hence the map $\Phi^* \circ R$ is compactly differentiable at any strictly increasing function $f_0 \in \ell^\infty((\tau_L, \tau_U))$ tangentially to $\mathcal{C}(\tau_L, \tau_U)$. Combining this with the representation $\Phi(f)(y) = \tau_L + (\tau_U - \tau_L)(\Phi^* \circ R)(f)(y)$ it follows that $\Phi$ is compactly differentiable at any strictly increasing function $f_0 \in \ell^\infty((\tau_L, \tau_U))$ with derivative $d\Phi_{f_0}(h)(y) = -h(f_0^{-1}(y))/f_0'(f_0^{-1}(y))$. Thus weak convergence of $a_n\big(\widehat{F}(y) - F(y)\big)$ follows from the functional delta method.

Next, observe that $\Psi(f) = \Theta \circ \Phi(f)$ where $\Theta(f)(\tau) = \inf\{y : f(y) \geq \tau\}$ denotes the generalized inverse. Compact differentiability of $\Theta$ at differentiable, strictly increasing functions $f_0$ tangentially to the space of contunuous functions is established in Lemma 3.9.23 of van der Vaart and Wellner (1996), and the derivative of $\Theta$ at $f_0$ is given by $d\Theta_{f_0}(h)(y) = -h(f_0^{-1}(y))/f_0'(f_0^{-1}(y))$. By the chain rule for Hadamard derivatives this implies compact differentiability of $\Psi$ tangentially to $\mathcal{C}(\tau_L, \tau_U)$. Thus the second weak convergence result again follows by the functional delta method. □

# APPENDIX B: Technical Remarks on Estimation Bias

**Remark B.1.** In this remark we show the bound $\widetilde{c}_n = o(m^{-\lfloor \eta \rfloor})$ for univariate spline models discussed in Example 2.3, as well as $c_n^\dagger = O(m^{-\lfloor \eta \rfloor/k'})$ for partial linear model in Section 3. We first show the latter.

Assume that $\mathcal{W} = [0,1]^{k'}$, that $h(\cdot;\cdot) \in \Lambda_c^\eta(\mathcal{W}, \mathcal{T})$ and that $\widetilde{\mathbf{Z}}$ corresponds to a tensor product $B$-spline basis of order $q$ on $\mathcal{W}$ with $m^{1/k'}$ equidistant knots in each coordinate. Moreover, assume that $(V, W)$ has a density $f_{V,W}$ such that $0 < \inf_{v,w} f_{V,W}(v,w) \leq \sup_{v,w} f_{V,W}(v,w) < \infty$. We shall show that in this case $c_n^\dagger = O(m^{-\lfloor \eta \rfloor/k'})$ where $c_n^\dagger$ is defined in Assumption (C1). Define

$$\boldsymbol{\beta}_{n,g}(\tau) := \underset{\boldsymbol{\beta} \in \mathbb{R}^m}{\operatorname{argmin}} \int \left(\widetilde{\mathbf{Z}}(w)^\top \boldsymbol{\beta} - h(w;\tau)\right)^2 \int f_{Y|X}(Q(v,w;\tau)|(v,w)) f_{V,W}(v,w) dv dw. \tag{B.1}$$

Note that $w \mapsto \widetilde{\mathbf{Z}}(w)^\top \boldsymbol{\beta}_{n,g}(\tau)$ can be viewed as a projection of a function $g : \mathcal{W} \to \mathbb{R}$ onto the spline space $\mathcal{B}_m(\mathcal{W}) := \{w \mapsto \widetilde{\mathbf{Z}}(w)^\top \mathbf{b} : \mathbf{b} \in \mathbb{R}^m\}$, with respect to the inner product $\langle g_1, g_2 \rangle = \int g_1(w) g_2(w) d\nu(w)$, where $d\nu(w) := \left(\int_v f_{Y|X}(Q(v,w;\tau)|v,w) f_{V,W}(v,w) dv\right) dw$.

We first apply Theorem A.1 on p.1630 of Huang (2003). To do so, we need to verify Condition A.1-A.3 of Huang (2003). Condition A.1 can be verified by invoking (A2)-(A3) in our paper and using the bounds on $f_W$. The choice of basis functions and knots ensures that Conditions A.2 and A.3 hold (see the discussion on p.1630 of Huang (2003)). Thus, Theorem A.1 on p.1630 of Huang (2003) implies there exists a constant $C$ independent of $n$ such that for any function on $\mathcal{W}$,

$$\sup_{w \in \mathcal{W}} \left|\widetilde{\mathbf{Z}}(w)^\top \boldsymbol{\beta}_{n,g}(\tau)\right| \leq C \sup_{w \in \mathcal{W}} |g(w)|.$$

Recall that $\mathcal{W}$ is a compact subset of $\mathbb{R}^d$ and $h(w;\tau) \in \Lambda_c^\eta(\mathcal{W}, \mathcal{T})$. Since $\mathcal{B}_m(\mathcal{W})$ is a finite dimensional vector space of functions, by a compactness argument there exists $g^*(\cdot;\tau) \in \mathcal{B}_m(\mathcal{W})$ such that $\sup_{w \in \mathcal{W}} |h(w;\tau) - g^*(w;\tau)| = \inf_{g \in \mathcal{B}_m(\mathcal{W})} \sup_{w \in \mathcal{W}} |h(w;\tau) - g(w)|$ for each fixed $\tau$. With $m > \eta$, the inequality in the proof for Theorem 12.8 in Schumaker (1981), with their "$m_i$" being our $\eta$ and $\overline{\Delta}_i \asymp m^{-1/k'}$ yields

$$\begin{aligned}
\widetilde{c}_n &= \sup_{w,\tau} \left|\widetilde{\mathbf{Z}}(w)^\top \boldsymbol{\beta}_{n,h(w;\tau)}(\tau) - h(w;\tau)\right| \\
&= \sup_{w,\tau} \left|\widetilde{\mathbf{Z}}(w)^\top \boldsymbol{\beta}_{n,h(w;\tau)}(\tau) - g^*(w;\tau) + g^*(w;\tau) - h(w;\tau)\right| \\
&\leq \sup_{w,\tau} \left|\widetilde{\mathbf{Z}}(w)^\top \boldsymbol{\beta}_{n,h(w;\tau)-g^*(w;\tau)}(\tau)\right| + \sup_{w,\tau} |g^*(w;\tau) - h(w;\tau)| \\
&\leq (C+1) \sup_{\tau \in \mathcal{T}} \inf_{g \in \mathcal{B}(\mathcal{W})} \sup_{x} |h(w;\tau) - g(w)| \\
&\lesssim m^{-\lfloor \eta \rfloor/k'} k' \max_{|\boldsymbol{j}| \leq \eta} \sup_{\tau \in \mathcal{T}} \sup_{w \in \mathcal{W}} |D^{\boldsymbol{j}} h(w;\tau)|,
\end{aligned}$$

where $\lfloor \eta \rfloor$ is the greatest integer less than $\eta$, $\max_{|\boldsymbol{j}| \leq \eta} \sup_{\tau \in \mathcal{T}} \sup_{x \in \mathcal{X}} |D^{\boldsymbol{j}} h(w;\tau)| = O(1)$ by the assumption that $h(w;\tau) \in \Lambda_c^\eta(\mathcal{W}, \mathcal{T})$ and fixed $k'$. An extension of Theorem 12.8 of Schumaker (1981) to Besov spaces (see Example 6.29 of Schumaker (1981)) in similar manner as Theorem 6.31 of Schumaker (1981) could refine the rate to $\widetilde{c}_n \lesssim m^{-\eta/k'}$, but we do not pursue this direction here.

Next we show the bound $\widetilde{c}_n = o(m^{-\lfloor \eta \rfloor})$ in the setting of Example 2.3. Assume the density $f_X(x)$ of $X$ exists and $0 < \inf_{x \in \mathcal{X}} f_X(x) \leq \sup_{x \in \mathcal{X}} f_X(x) < \infty$. Define the measure $\nu(u)$ by $d\nu(u) = f(Q(u;\tau)|u) f_X(u) du$. Thus, $x \mapsto \mathbf{B}(x)^\top \boldsymbol{\beta}_{n,g}(\tau)$ with $\boldsymbol{\beta}_{n,g}$ defined similarly to (B.1) is now viewed as a projection of a function $g : \mathcal{X} \to \mathbb{R}$ onto the space $\mathcal{B}(\mathcal{X})$ with respect to the inner product $\langle g_1, g_2 \rangle = \int g_1(u) g_2(u) d\nu(u)$. The remaining proof is similar to the partial linear model, with $h(w;\tau)$ being replaced by $Q(x;\tau)$, and we omit the details. □

# SUPPLEMENTARY MATERIAL: QUANTILE PROCESS FOR SEMI AND NONPARAMETRIC REGRESSION

In this supplemental material, we provide the auxiliary proofs needed in the appendices. Section S.1 develops the technicalities for Bahadur representations. Section S.2 presents some empirical process results and computes the covering number of some function classes encountered in the proofs of asymptotic tightness of quantile process.

## S.1: Proofs for Bahadur Representations

### S.1.1. Proof of Theorem 5.1

Some rearranging of terms yields

$$\begin{aligned}
\mathbb{P}_n \psi(\cdot; \widehat{\boldsymbol{\gamma}}(\tau), \tau) &= n^{-1/2}\mathbb{G}_n(\psi(\cdot; \widehat{\boldsymbol{\gamma}}(\tau), \tau)) - n^{-1/2}\mathbb{G}_n(\psi(\cdot; \boldsymbol{\gamma}_n(\tau), \tau)) \\
&\quad + \widetilde{J}_m(\tau)(\widehat{\boldsymbol{\gamma}}(\tau) - \boldsymbol{\gamma}_n(\tau)) + n^{-1/2}\mathbb{G}_n(\psi(\cdot; \boldsymbol{\gamma}_n(\tau), \tau)) \\
&\quad + \mu(\boldsymbol{\gamma}_n(\tau), \tau) + \Big(\mu(\widehat{\boldsymbol{\gamma}}(\tau), \tau) - \mu(\boldsymbol{\gamma}_n(\tau), \tau) - \widetilde{J}_m(\tau)(\widehat{\boldsymbol{\gamma}}(\tau) - \boldsymbol{\gamma}_n(\tau))\Big).
\end{aligned}$$

In other words

$$\widehat{\boldsymbol{\gamma}}(\tau) - \boldsymbol{\gamma}_n(\tau) = -n^{-1/2} J_m(\tau)^{-1} \mathbb{G}_n(\psi(\cdot; \boldsymbol{\gamma}_n(\tau), \tau)) + r_{n,1}(\tau) + r_{n,2}(\tau) + r_{n,3}(\tau) + r_{n,4}(\tau) \quad \text{(S.1.1)}$$

where

$$\begin{aligned}
r_{n,1}(\tau) &:= \widetilde{J}_m(\tau)^{-1} \mathbb{P}_n \psi(\cdot; \widehat{\boldsymbol{\gamma}}(\tau), \tau), \\
r_{n,2}(\tau) &:= -\widetilde{J}_m(\tau)^{-1}\Big(\mu(\widehat{\boldsymbol{\gamma}}(\tau), \tau) - \mu(\boldsymbol{\gamma}_n(\tau), \tau) - \widetilde{J}_m(\tau)(\widehat{\boldsymbol{\gamma}}(\tau) - \boldsymbol{\gamma}_n(\tau))\Big), \\
r_{n,3}(\tau) &:= -n^{-1/2} \widetilde{J}_m(\tau)^{-1}\Big(\mathbb{G}_n(\psi(\cdot; \widehat{\boldsymbol{\gamma}}(\tau), \tau)) - \mathbb{G}_n(\psi(\cdot; \boldsymbol{\gamma}_n(\tau), \tau))\Big), \\
r_{n,4}(\tau) &:= -n^{-1/2}(J_m(\tau)^{-1} - \widetilde{J}_m(\tau)^{-1})\mathbb{G}_n(\psi(\cdot; \boldsymbol{\gamma}_n(\tau), \tau)) - \widetilde{J}_m(\tau)^{-1}\mu(\boldsymbol{\gamma}_n(\tau), \tau).
\end{aligned}$$

The remaining proof consists in bounding the individual remainder terms.

The bound on $r_{n,1}$ follows from results on duality theory for convex optimization, see Lemma 26 on page 66 in Belloni et al. (2016) for a proof.

To bound $r_{n,2}$ and $r_{n,3}$, define the class of functions

$$\mathcal{G}_1 := \big\{(\mathbf{Z}, Y) \mapsto \mathbf{a}^\top \mathbf{Z}(\mathbf{1}\{Y \leq \mathbf{Z}^\top \mathbf{b}\} - \tau)\mathbf{1}\{\|\mathbf{Z}\| \leq \xi_m\} \big| \tau \in \mathcal{T}, \mathbf{b} \in \mathbb{R}^m, \mathbf{a} \in \mathcal{S}^{m-1}\big\}. \quad \text{(S.1.2)}$$

Moreover, let

$$s_{n,1} := \|\mathbb{P}_n - P\|_{\mathcal{G}_1}.$$

Observe that by Lemma S.1.2 with $t = 2$ we have

$$\Omega_{1,n} := \Big\{\sup_{\tau \in \mathcal{T}} \|\widehat{\boldsymbol{\gamma}}(\tau) - \boldsymbol{\gamma}_n(\tau)\| \leq \frac{4(s_{n,1} + g_n)}{\inf_{\tau \in \mathcal{T}} \lambda_{\min}(\widetilde{J}_m(\tau))}\Big\} \supseteq \Big\{s_{n,1} + g_n < \frac{\inf_{\tau \in \mathcal{T}} \lambda_{\min}^2(\widetilde{J}_m(\tau))}{8\xi_m \overline{f'} \lambda_{\max}(\mathbb{E}[\mathbf{Z}\mathbf{Z}^\top])}\Big\} =: \Omega_{2,n}.$$

Define the event

$$\Omega_{3,n} := \Big\{s_{n,1} \leq C\Big[\Big(\frac{m}{n}\log n\Big)^{1/2} + \frac{m\xi_m}{n}\log n + \Big(\frac{\kappa_n}{n}\Big)^{1/2} + \frac{\xi_m \kappa_n}{n}\Big]\Big\}.$$

Now it follows from Lemma S.1.3 that $P(\Omega_{3,n}) \geq 1 - e^{-\kappa_n}$ [note that $\xi_m = O(n^b)$ yields $\log \xi_m = O(\log n)$]. Moreover, the assumption $m\xi_m^2 \log n = o(n), \xi_m = O(n^b), \xi_m g_n = o(1)$ implies that for $\kappa_n \ll n/\xi_m^2$ and large enough $n$,

$$C\Big[\Big(\frac{m}{n}\log n\Big)^{1/2} + \frac{m\xi_m}{n}\log n + \Big(\frac{\kappa_n}{n}\Big)^{1/2} + \frac{\xi_m \kappa_n}{n}\Big] + g_n \leq \frac{\inf_{\tau \in \mathcal{T}} \lambda_{\min}^2(\widetilde{J}_m(\tau))}{8\xi_m \overline{f'}\lambda_{\max}(\mathbb{E}[\mathbf{Z}\mathbf{Z}^\top])} \quad (S.1.3)$$

for $n$ large enough. Thus, for all $n$ for which (S.1.3) holds, $\Omega_{3,n} \subseteq \Omega_{2,n} \subseteq \Omega_{1,n}$. From this we obtain that on $\Omega_{3,n}$, for a constant $C_2$ which is independent of $n$, we have

$$\sup_{\tau \in \mathcal{T}} \|\widehat{\boldsymbol{\gamma}}(\tau) - \boldsymbol{\gamma}_n(\tau)\| \leq C_2\Big[\Big(\frac{m}{n}\log n\Big)^{1/2} + \frac{m\xi_m}{n}\log n + \Big(\frac{\kappa_n}{n}\Big)^{1/2} + \frac{\xi_m \kappa_n}{n} + g_n\Big].$$

In particular, for all $n$ for which (S.1.3) holds,

$$P\Big(\sup_{\tau \in \mathcal{T}} \|\widehat{\boldsymbol{\gamma}}(\tau) - \boldsymbol{\gamma}_n(\tau)\| \leq C_2\Big[\Big(\frac{m}{n}\log n\Big)^{1/2} + \frac{m\xi_m}{n}\log n + \Big(\frac{\kappa_n}{n}\Big)^{1/2} + \frac{\xi_m \kappa_n}{n} + g_n\Big]\Big) \geq 1 - e^{-\kappa_n}. \quad (S.1.4)$$

The bound on $r_{n,2}$ is now a direct consequence of Lemma S.1.1 and the fact that $n^{-1}m\xi_m \log n = o((n^{-1}m\log n)^{1/2})$ and $\xi_m \kappa_n n^{-1} = o(\kappa_n^{1/2} n^{-1/2})$. The bound on $r_{n,4}$ follows once we observe that for any $\mathbf{a} \in \mathcal{S}^{m-1}$

$$|\mathbf{a}^\top (J_m(\tau) - \widetilde{J}_m(\tau))\mathbf{a}| \leq \overline{f'} c_n \lambda_{\max}(\mathbb{E}[\mathbf{Z}\mathbf{Z}^\top]). \quad (S.1.5)$$

Together with the identity $A^{-1} - B^{-1} = B^{-1}(B-A)A^{-1}$ this implies that for sufficiently large $n$ we have $\sup_{\tau \in \mathcal{T}} \|r_{n,4}(\tau)\| \leq C_1(c_n s_{n,1} + g_n)$ for a constant $C_1$ which does not depend on $n$.

Thus, it remains to bound $r_{n,3}$. Observe that on the set $\{\sup_{\tau \in \mathcal{T}} \|\widehat{\boldsymbol{\gamma}}(\tau) - \boldsymbol{\gamma}_n(\tau)\| \leq \delta\}$ we have the bound

$$\sup_{\tau \in \mathcal{T}} \|r_{n,3}(\tau)\| \leq \frac{1}{\inf_{\tau \in \mathcal{T}} \lambda_{\min}(\widetilde{J}_m(\tau))} \|\mathbb{P}_n - P\|_{\mathcal{G}_2(\delta)},$$

where the class of functions $\mathcal{G}_2(\delta)$ is defined as follows

$$\mathcal{G}_2(\delta) := \big\{(\mathbf{Z}, Y) \mapsto \mathbf{a}^\top \mathbf{Z}(\mathbf{1}\{Y \leq \mathbf{Z}^\top \mathbf{b}_1\} - \mathbf{1}\{Y \leq \mathbf{Z}^\top \mathbf{b}_2\})\mathbf{1}\{\|\mathbf{Z}\| \leq \xi_m\}\big|$$
$$\mathbf{b}_1, \mathbf{b}_2 \in \mathbb{R}^m, \|\mathbf{b}_1 - \mathbf{b}_2\| \leq \delta, \mathbf{a} \in \mathcal{S}^{m-1}\big\}. \quad (S.1.6)$$

It thus follows that for any $\delta, \alpha > 0$

$$P\Big(\sup_{\tau \in \mathcal{T}} \|r_{n,3}(\tau)\| \geq \alpha\Big) \leq P\Big(\sup_{\tau \in \mathcal{T}} \|\widehat{\boldsymbol{\gamma}}(\tau) - \boldsymbol{\gamma}_n(\tau)\| \geq \delta\Big) + P\Big(\frac{\|\mathbb{P}_n - P\|_{\mathcal{G}_2(\delta)}}{\inf_{\tau \in \mathcal{T}} \lambda_{\min}(\widetilde{J}_m(\tau))} \geq \alpha\Big).$$

Letting $\delta := C((n^{-1}m\log n)^{1/2} + (\kappa_n/n)^{1/2} + g_n)$ and

$$\alpha := C\zeta_n(\delta_n, \kappa_n)$$
$$= C\Big\{\Big(\Big(\frac{m\log n}{n}\Big)^{1/2} + \Big(\frac{\kappa_n}{n}\Big)^{1/2} + g_n\Big)^{1/2}\Big(\Big(\frac{m\xi_m \log n}{n}\Big)^{1/2} + \Big(\frac{\xi_m \kappa_n}{n}\Big)^{1/2}\Big) + \frac{m\xi_m \log n}{n} + \frac{\xi_m \kappa_n}{n}\Big\}, \quad (S.1.7)$$

where $\zeta_n$ is defined in (S.1.12) in the statement of Lemma S.1.3, with a suitable constant $C$. Observe that the assumption $m\xi_m^2 \log n = o(n)$ implies $(mn^{-1}\log n)^{1/4} > (m\xi_m n^{-1}\log n)^{1/2}$ and $(\kappa_n/n)^{1/4} > (\xi_m \kappa_n/n)^{1/2}$ for sufficiently large $n$, so the last two terms are less than the first term in (S.1.7). Hence, for some large enough constant $C > 0$,

$$\alpha \leq C\Big(\Big(\frac{m\log n}{n}\Big)^{1/2} + \Big(\frac{\kappa_n}{n}\Big)^{1/2} + g_n\Big)^{1/2}\Big(\Big(\frac{m\xi_m \log n}{n}\Big)^{1/2} + \Big(\frac{\xi_m \kappa_n}{n}\Big)^{1/2}\Big)$$

Finally, the bounds in Lemma S.1.3 and (S.1.4) yield the desired bound. □

*S.1.1.1. Technical details for the proof of Theorem 5.1*

**Lemma S.1.1.** *Under assumptions (A1)-(A3) we have for any $\delta > 0$,*

$$\sup_{\tau \in \mathcal{T}} \sup_{\|\mathbf{b} - \boldsymbol{\gamma}_n(\tau)\| \leq \delta} \|\mu(\mathbf{b}, \tau) - \mu(\boldsymbol{\gamma}_n(\tau), \tau) - \widetilde{J}_m(\tau)(\mathbf{b} - \boldsymbol{\gamma}_n(\tau))\| \leq \lambda_{\max}(\mathbb{E}[\mathbf{Z}\mathbf{Z}^\top])\overline{f'}\delta^2 \xi_m,$$

**Proof of Lemma S.1.1.** Note that $\mu'(\boldsymbol{\gamma}_n(\tau), \tau) = \mathbb{E}[\mathbf{Z}\mathbf{Z}^\top f_{Y|X}(\mathbf{Z}^\top \boldsymbol{\gamma}_n(\tau)|X)] = \widetilde{J}_m(\tau)$ where we use the notation $\mu'(\mathbf{b}, \tau) := \partial_{\mathbf{b}} \mu(\mathbf{b}, \tau)$. Additionally, we have

$$\mu(\mathbf{b}, \tau) = \mu(\boldsymbol{\gamma}_n(\tau), \tau) + \mu'(\bar{\boldsymbol{\gamma}}_n, \tau)(\mathbf{b} - \boldsymbol{\gamma}_n(\tau)),$$

where $\bar{\boldsymbol{\gamma}}_n = \mathbf{b} + \lambda_{\mathbf{b},\tau}(\boldsymbol{\gamma}_n(\tau) - \mathbf{b})$ for some $\lambda_{\mathbf{b},\tau} \in [0,1]$. Moreover, for any $\mathbf{a} \in \mathbb{R}^m$,

$$\mathbf{a}^\top [\mu(\mathbf{b},\tau) - \mu(\boldsymbol{\gamma}_n(\tau),\tau) - \mu'(\boldsymbol{\gamma}_n(\tau),\tau)(\mathbf{b}-\boldsymbol{\gamma}_n(\tau))] = \mathbf{a}^\top[(\mu'(\bar{\boldsymbol{\gamma}}_n,\tau) - \mu'(\boldsymbol{\gamma}_n(\tau),\tau))(\mathbf{b}-\boldsymbol{\gamma}_n(\tau))]$$

and thus we have for any $\|\mathbf{b} - \boldsymbol{\gamma}_n(\tau)\| \leq \delta$

$$\begin{aligned}
&\|\mu(\mathbf{b},\tau) - \mu(\boldsymbol{\gamma}_n(\tau),\tau) - \mu'(\boldsymbol{\gamma}_n(\tau),\tau)(\mathbf{b}-\boldsymbol{\gamma}_n(\tau))\| \\
&\leq \sup_{\|\mathbf{a}\|=1} \left|\mathbb{E}\left[(\mathbf{a}^\top \mathbf{Z})\mathbf{Z}^\top(\mathbf{b}-\boldsymbol{\gamma}_n(\tau))\left(f_{Y|X}(\mathbf{Z}^\top \bar{\boldsymbol{\gamma}}_n|X) - f_{Y|X}(\mathbf{Z}^\top \boldsymbol{\gamma}_n(\tau)|X)\right)\right]\right| \\
&\leq \overline{f'} \sup_{\|\mathbf{a}\|=1} \mathbb{E}\left[|\mathbf{a}^\top \mathbf{Z}||\mathbf{Z}^\top(\bar{\boldsymbol{\gamma}}_n - \boldsymbol{\gamma}_n(\tau))||\mathbf{Z}^\top(\mathbf{b}-\boldsymbol{\gamma}_n(\tau))|\right] \\
&\leq \overline{f'}\xi_m \mathbb{E}\left[|\mathbf{Z}^\top(\bar{\boldsymbol{\gamma}}_n - \boldsymbol{\gamma}_n(\tau))||\mathbf{Z}^\top(\mathbf{b}-\boldsymbol{\gamma}_n(\tau))|\right] \\
&\leq \xi_m \delta^2 \overline{f'} \sup_{\|\mathbf{a}\|=1} \mathbb{E}[|\mathbf{a}^\top \mathbf{Z}|^2],
\end{aligned}$$

here the last inequality follows by Chauchy-Schwarz. Since the last line does not depend on $\tau, \mathbf{b}$, this completes the proof. $\square$

**Lemma S.1.2.** *Let assumptions (A1)-(A3) hold. Then, for any $t > 1$*

$$\left\{\sup_{\tau \in \mathcal{T}} \|\widehat{\boldsymbol{\gamma}}(\tau) - \boldsymbol{\gamma}_n(\tau)\| \leq \frac{2t(s_{n,1}+g_n)}{\inf_{\tau \in \mathcal{T}} \lambda_{\min}(\widetilde{J}_m(\tau))}\right\} \supseteq \left\{(s_{n,1}+g_n) < \frac{\inf_{\tau \in \mathcal{T}} \lambda^2_{\min}(\widetilde{J}_m(\tau))}{4t\xi_m \overline{f'}\lambda_{\max}(\mathbb{E}[\mathbf{Z}\mathbf{Z}^\top])}\right\},$$

*where $s_{n,1} := \|\mathbb{P}_n - P\|_{\mathcal{G}_1}$ and $\mathcal{G}_1$ is defined in (S.1.2).*

**Proof of Lemma S.1.2.** Observe that $f: \mathbf{b} \mapsto \mathbb{P}_n \rho_\tau(Y_i - \mathbf{Z}_i^\top \mathbf{b})$ is convex, and the vector $\mathbb{P}_n \psi(\cdot; \mathbf{b}, \tau)$ is a subgradient of $f$ at the point $\mathbf{b}$. Recalling that $\widehat{\boldsymbol{\gamma}}(\tau)$ is a minimizer of $\mathbb{P}_n \rho_\tau(Y_i - \mathbf{Z}_i^\top \mathbf{b})$, it follows that for any $a > 0$,

$$\{\sup_{\tau \in \mathcal{T}} \|\widehat{\boldsymbol{\gamma}}(\tau) - \boldsymbol{\gamma}_n(\tau)\| \leq a(s_{n,1}+g_n)\} \supseteq \{\inf_\tau \inf_{\|\boldsymbol{\delta}\|=1} \boldsymbol{\delta}^\top \mathbb{P}_n \psi(\cdot; \boldsymbol{\gamma}_n(\tau) + a(s_{n,1}+g_n)\boldsymbol{\delta}, \tau) > 0\}. \tag{S.1.8}$$

To see this, define $\boldsymbol{\delta} := (\widehat{\boldsymbol{\gamma}}(\tau) - \boldsymbol{\gamma}_n(\tau))/\|\widehat{\boldsymbol{\gamma}}(\tau) - \boldsymbol{\gamma}_n(\tau)\|$ and note that by definition of the subgradient we have for any $\widetilde{\zeta}_n > 0$

$$\mathbb{P}_n \rho_\tau(Y_i - \mathbf{Z}_i^\top \widehat{\boldsymbol{\gamma}}(\tau)) \geq \mathbb{P}_n \rho_\tau(Y_i - \mathbf{Z}_i^\top(\boldsymbol{\gamma}_n(\tau) + \widetilde{\zeta}_n \boldsymbol{\delta})) + (\|\widehat{\boldsymbol{\gamma}}(\tau) - \boldsymbol{\gamma}_n(\tau)\| - \widetilde{\zeta}_n)\boldsymbol{\delta}^\top \mathbb{P}_n \psi(\cdot; \boldsymbol{\gamma}_n(\tau) + \widetilde{\zeta}_n \boldsymbol{\delta}, \tau).$$

Set $\widetilde{\zeta}_n = a(s_{n,1}+g_n)$. By the definition of $\widehat{\boldsymbol{\gamma}}(\tau)$ as minimizer, the inequality above can only be true if $(\|\widehat{\boldsymbol{\gamma}}(\tau) - \boldsymbol{\gamma}_n(\tau)\| - \widetilde{\zeta}_n)\boldsymbol{\delta}^\top \mathbb{P}_n \psi(\cdot; \boldsymbol{\gamma}_n(\tau) + \widetilde{\zeta}_n \boldsymbol{\delta}, \tau) \leq 0$, which yields (S.1.8).

The proof is finished once we minorize the empirical score $\boldsymbol{\delta}^\top \mathbb{P}_n \psi(\cdot; \boldsymbol{\gamma}_n(\tau) + a(s_{n,1}+g_n)\boldsymbol{\delta}, \tau)$ in (S.1.8) in terms of $s_{n,1} + g_n$. To proceed, observe that under assumptions (A1)-(A3) we have by Lemma S.1.1

$$\sup_{\|\boldsymbol{\delta}\|=1} \left|\mathbb{E}[\boldsymbol{\delta}^\top \{\psi(Y, \mathbf{Z}; \mathbf{b}, \tau) - \psi(Y, \mathbf{Z}; \boldsymbol{\gamma}_n(\tau), \tau) - \mathbf{Z}\mathbf{Z}^\top f_{Y|X}(\boldsymbol{\gamma}_n(\tau)^\top \mathbf{Z}|X)(\mathbf{b} - \boldsymbol{\gamma}_n(\tau))\}]\right|$$
$$\leq \xi_m \overline{f'} \lambda_{\max}(\mathbb{E}[\mathbf{Z}\mathbf{Z}^\top]) \|\mathbf{b} - \boldsymbol{\gamma}_n(\tau)\|^2. \tag{S.1.9}$$

Therefore, we have for arbitrary $\|\boldsymbol{\delta}\| = 1, \tau \in \mathcal{T}$ that

$$\boldsymbol{\delta}^\top \mathbb{P}_n \psi(\cdot; \boldsymbol{\gamma}_n(\tau) + a(s_{n,1} + g_n)\boldsymbol{\delta}, \tau)$$
$$\geq -s_{n,1} - g_n + \boldsymbol{\delta}^\top \left( \mathbb{E}[\psi(Y, \mathbf{Z}; \boldsymbol{\gamma}_n(\tau) + a(s_{n,1} + g_n)\boldsymbol{\delta}, \tau)] - \mathbb{E}[\psi(Y, \mathbf{Z}; \boldsymbol{\gamma}_n(\tau), \tau)] \right)$$
$$\geq a(s_{n,1} + g_n) \inf_{\tau \in \mathcal{T}} \lambda_{\min}(\widetilde{J}_m(\tau)) - s_{n,1} - g_n - \xi_m \overline{f'} \lambda_{\max}(\mathbb{E}[\mathbf{Z}\mathbf{Z}^\top]) a^2 (s_{n,1} + g_n)^2,$$

where for the first inequality we recall the definition $g_n = \sup_{\tau \in \mathcal{T}} \|\mathbb{E}[\psi(Y, \mathbf{Z}; \boldsymbol{\gamma}_n(\tau), \tau)]\|$ and $s_{n,1} = \|\mathbb{P}_n - P\|_{\mathcal{G}_1}$; the second inequality follows by (S.1.9). Setting $a = 2t/\inf_{\tau \in \mathcal{T}} \lambda_{\min}(\widetilde{J}_m(\tau))$ in $\zeta_n$, we see that the right-hand side of the display above is positive when

$$s_{n,1} + g_n < \frac{(2t-1)\inf_{\tau \in \mathcal{T}} \lambda_{\min}^2(\widetilde{J}_m(\tau))}{4t^2 \xi_m \overline{f'} \lambda_{\max}(\mathbb{E}[\mathbf{Z}\mathbf{Z}^\top])}.$$

Observing that for $t > 1$ we have $(2t-1)/t^2 \geq 1/t$ and plugging $a = 2t/\inf_{\tau \in \mathcal{T}} \lambda_{\min}(\widetilde{J}_m(\tau))$ in equation (S.1.8) completes the proof. $\square$

**Lemma S.1.3.** *Consider the classes of functions $\mathcal{G}_1, \mathcal{G}_2(\delta)$ defined in (S.1.2) and (S.1.6), respectively. Under assumptions (A1)-(A3) we have for some constant $C$ independent of $n$ and all $\kappa_n > 0$ provided that $\xi_m = O(n^b)$ for some fixed $b$*

$$P\left( \|\mathbb{P}_n - P\|_{\mathcal{G}_1} \geq C\left[ \left(\frac{m}{n} \log \xi_m\right)^{1/2} + \frac{m\xi_m}{n} \log \xi_m + \left(\frac{\kappa_n}{n}\right)^{1/2} + \frac{\xi_m \kappa_n}{n} \right] \right) \leq e^{-\kappa_n}. \quad (S.1.10)$$

*For any $\delta_n$ satisfying $\xi_m \delta_n \gg n^{-1}$, we have for sufficiently large $n$ and arbitrary $\kappa_n > 0$*

$$P\left( \|\mathbb{P}_n - P\|_{\mathcal{G}_2(\delta_n)} \geq C\zeta_n(\delta_n, \kappa_n) \right) \leq e^{-\kappa_n}, \quad (S.1.11)$$

*where*

$$\zeta_n(t, \kappa_n) := t^{1/2} \left( \frac{m\xi_m}{n} \log(\xi_m \vee n) \right)^{1/2} + \frac{m\xi_m}{n} \log(\xi_m \vee n) + t^{1/2} \left( \frac{\xi_m \kappa_n}{n} \right)^{1/2} + \frac{\xi_m \kappa_n}{n}. \quad (S.1.12)$$

**Proof of Lemma S.1.3.** Observe that for each $f \in \mathcal{G}_1$ we have $|f(x,y)| \leq \xi_m$, and the same holds for $\mathcal{G}_2(\delta)$ for any value of $\delta$. Additionally, similar arguments as those in the proof of Lemma 18 in Belloni et al. (2016) imply together with Theorem 2.6.7 in van der Vaart and Wellner (1996) that, almost surely,

$$N(\mathcal{G}_2(\delta), L_2(\mathbb{P}_n); \varepsilon) \leq \left( \frac{A\|F\|_{L^2(\mathbb{P}_n)}}{\varepsilon} \right)^{v_1(m)}, \quad N(\mathcal{G}_1, L_2(\mathbb{P}_n); \varepsilon) \leq \left( \frac{A\|F\|_{L^2(\mathbb{P}_n)}}{\varepsilon} \right)^{v_2(m)},$$

where $A$ is some constant and $v_1(m) = O(m), v_2(m) = O(m)$. Finally, for each $f \in \mathcal{G}_1$ we have

$$\mathbb{E}[f^2] \leq \sup_{\|\mathbf{a}\|=1} \mathbf{a}^\top \mathbb{E}[\mathbf{Z}\mathbf{Z}^\top] \mathbf{a} = \lambda_{\max}(\mathbb{E}[\mathbf{Z}\mathbf{Z}^\top]).$$

On the other hand, each $f \in \mathcal{G}_2(\delta_n)$ satisfies

$$\mathbb{E}[f^2] \leq \sup_{\|\mathbf{a}\|=1} \sup_{\|\mathbf{b}_1 - \mathbf{b}_2\| \leq \delta_n} \mathbb{E}\left[ (\mathbf{a}^\top \mathbf{Z})^2 \mathbf{1}\{|Y - \mathbf{Z}^\top \mathbf{b}_1| \leq |\mathbf{Z}^\top (\mathbf{b}_1 - \mathbf{b}_2)|\} \right]$$
$$\leq \sup_{\mathbf{b} \in \mathbb{R}^m} \sup_{\|\mathbf{a}\|=1} \mathbb{E}\left[ (\mathbf{a}^\top \mathbf{Z})^2 \mathbf{1}\{|Y - \mathbf{Z}^\top \mathbf{b}| \leq \xi_m \delta_n\} \right]$$
$$\leq 2\overline{f} \xi_m \delta_n \lambda_{\max}(\mathbb{E}[\mathbf{Z}\mathbf{Z}^\top]).$$

Note that under assumptions (A1)-(A3) the right-hand side is bounded by $c\xi_m \delta_n$ where $c$ is a constant that does not depend on $n$. Thus the bound in (S.2.2) implies that for $\xi_m \delta_n \gg n^{-1}$ we have for some constant $C$ which is independent of $n$,

$$\mathbb{E}\|\mathbb{P}_n - P\|_{\mathcal{G}_1} \leq C\left[ \left(\frac{m}{n} \log(\xi_m \vee n)\right)^{1/2} + \frac{m\xi_m}{n} \log(\xi_m \vee n) \right], \quad (S.1.13)$$

$$\mathbb{E}\|\mathbb{P}_n - P\|_{\mathcal{G}_2(\delta_n)} \leq C\left[ \xi_m^{1/2} \delta_n^{1/2} \left( \frac{m}{n} \log(\xi_m \vee n) \right)^{1/2} + \frac{m\xi_m}{n} \log(\xi_m \vee n) \right]. \quad (S.1.14)$$

Thus (S.1.10) and (S.1.11) follow from the bound in (S.2.3) by setting $t = \kappa_n$. $\square$

### *S.1.2. Proof of Theorem 5.2*

We begin with the following useful decomposition

$$\widehat{\boldsymbol{\beta}}(\tau) - \boldsymbol{\beta}_n(\tau) = -n^{-1/2}\widetilde{J}_m^{-1}(\tau)\mathbb{G}_n(\psi(\cdot;\boldsymbol{\beta}_n(\tau),\tau)) + \widetilde{J}_m^{-1}(\tau)\sum_{k=1}^{4} R_{n,k}(\tau) \quad (\text{S.1.15})$$

where

$$\begin{aligned} R_{n,1}(\tau) &:= \mathbb{P}_n \psi(\cdot;\widehat{\boldsymbol{\beta}}(\tau),\tau), \\ R_{n,2}(\tau) &:= -\Big(\mu(\widehat{\boldsymbol{\beta}}(\tau),\tau) - \mu(\boldsymbol{\beta}_n(\tau),\tau) - \widetilde{J}_m(\tau)(\widehat{\boldsymbol{\beta}}(\tau) - \boldsymbol{\beta}_n(\tau))\Big), \\ R_{n,3}(\tau) &:= -n^{-1/2}\Big(\mathbb{G}_n(\psi(\cdot;\widehat{\boldsymbol{\beta}}(\tau),\tau)) - \mathbb{G}_n(\psi(\cdot;\boldsymbol{\beta}_n(\tau),\tau))\Big), \\ R_{n,4}(\tau) &:= -\mu(\boldsymbol{\beta}_n(\tau),\tau). \end{aligned}$$

Define $r_{n,2}(\tau,\mathbf{u}_n) := \mathbf{u}_n^\top \widetilde{J}_m^{-1}(\tau) R_{n,2}(\tau)$, $r_{n,k}(\tau,\mathbf{u}_n) := (\mathbf{u}_n^\top \widetilde{J}_m^{-1}(\tau))^{(\mathcal{I}(\mathbf{u}_n,D))} R_{n,k}(\tau)$ for $k=1,3$ and

$$r_{n,4}(\tau,\mathbf{u}_n) := \mathbf{u}_n^\top \widetilde{J}_m^{-1}(\tau) R_{n,4}(\tau) + \Big(\mathbf{u}_n^\top \widetilde{J}_m^{-1}(\tau) - (\mathbf{u}_n^\top \widetilde{J}_m^{-1}(\tau))^{(\mathcal{I}(\mathbf{u}_n,D))}\Big)\Big(R_{n,1}(\tau) + R_{n,3}(\tau)\Big).$$

With those definitions we obtain

$$\mathbf{u}_n^\top\big(\widehat{\boldsymbol{\beta}}(\tau) - \boldsymbol{\beta}_n(\tau)\big) = -n^{-1/2}\mathbf{u}_n^\top \widetilde{J}_m^{-1}(\tau)\mathbb{G}_n(\psi(\cdot;\boldsymbol{\beta}_n(\tau),\tau)) + \sum_{k=1}^{4} r_{n,k}(\tau,\mathbf{u}_n).$$

We will now show that the terms $r_{n,k}(\tau,\mathbf{u}_n)$ defined above satisfy the bounds given in the statement of Theorem 5.2 if we let $D = c\log n$ for a sufficiently large constant $c$.

The bound on $r_{n,1}$ follows from Lemma S.1.7. To bound $r_{n,2}$ apply Lemma S.1.4 and Lemma S.1.6. To bound $r_{n,3}$ observe that by Lemma S.1.4 the probability of the event

$$\Omega_1 := \Big\{ \sup_{\tau,x} |\mathbf{B}(x)^\top \widehat{\boldsymbol{\beta}}(\tau) - \mathbf{B}(x)^\top \boldsymbol{\beta}_n(\tau)| \leq C\Big(\widetilde{c}_n^2 + \frac{\xi_m(\log n + \kappa_n^{1/2})}{n^{1/2}}\Big)\Big\}.$$

is at least $1 - (m+1)e^{-\kappa_n}$. Letting $\delta_n := C\Big(\widetilde{c}_n^2 + \frac{\xi_m(\log n \vee \kappa_n^{1/2})}{n^{1/2}}\Big)$ we find that on $\Omega_1$

$$\sup_{\mathbf{u}_n \in \mathcal{S}_\mathcal{I}^{m-1}} \sup_\tau |r_{n,3}(\tau,\mathbf{u}_n)| \lesssim \|\mathbb{P}_n - P\|_{\mathcal{G}_2(\delta_n,\mathcal{I}(\mathbf{u}_n,D),\mathcal{I}'(\mathbf{u}_n,D))},$$

this follows since

$$n^{-1/2}(\mathbf{u}_n^\top \widetilde{J}_m^{-1}(\tau))^{(\mathcal{I}(\mathbf{u}_n,D))}\mathbb{G}_n(\psi(\cdot;\widehat{\boldsymbol{\beta}}(\tau),\tau)) = n^{-1/2}(\mathbf{u}_n^\top \widetilde{J}_m^{-1}(\tau))^{(\mathcal{I}(\mathbf{u}_n,D))}\mathbb{G}_n(\psi^{(\mathcal{I}'(\mathbf{u}_n,D))}(\cdot;\widehat{\boldsymbol{\beta}}(\tau),\tau))$$

and a similar identity holds with $\boldsymbol{\beta}_n$ instead of $\widehat{\boldsymbol{\beta}}$. Note that $\max\{|\mathcal{I}(\mathbf{u}_n,D)|, |\mathcal{I}'(\mathbf{u}_n,D)|\} \lesssim L + c\log n$. Hence, for any $\omega_n > 0$,

$$P\Big(\sup_{\mathbf{u}_n \in \mathcal{S}_\mathcal{I}^{m-1}} \sup_\tau |r_{n,3}(\tau,\mathbf{u}_n)| > \omega_n\Big)$$
$$\leq P\big(\sup_{\tau,x} |\mathbf{B}(x)^\top \widehat{\boldsymbol{\beta}}(\tau) - \mathbf{B}(x)^\top \boldsymbol{\beta}_n(\tau)| > \delta_n\big) + P\big(\|\mathbb{P}_n - P\|_{\mathcal{G}_2(\delta_n,\mathcal{I}(\mathbf{u}_n,D),\mathcal{I}'(\mathbf{u}_n,D))} > \omega_n\big).$$

The bound on $r_{n,3}$ now follows form the bound for the event $\Omega_1^c$, Lemma S.1.5 under $\omega_n := C\zeta(\delta_n,\mathcal{I}(\mathbf{u}_n,D),\mathcal{I}'(\mathbf{u}_n,D),\kappa_n)$ and the observation from the assumption $m\xi_m^2(\log n)^2 = o(n)$ that $\xi_m n^{-1}(\log n)^2 < \xi_m^{1/2} n^{-3/4}(\log n)^{3/2}$ when $n$ is sufficiently large. To bound the first part of $r_{n,4}(\tau,\mathbf{u}_n)$, we proceed as in the proof of equation (S.1.19) to obtain

$$\Big|(\mathbf{u}_n^\top \widetilde{J}_m^{-1}(\tau))\mu(\boldsymbol{\beta}_n(\tau),\tau)\Big| \leq \frac{1}{2}\overline{f'}\widetilde{c}_n^2 \widetilde{\mathcal{E}}(\mathbf{u}_n,\mathbf{B})$$

where the last line follows after a Taylor expansion. To bound the second part of $r_{n,4}(\tau,\mathbf{u}_n)$ note that $\sup_\tau \|R_{n,1}(\tau)\| + \|R_{n,3}(\tau)\| \leq 3\xi_m$ almost surely and thus choosing $D = c\log n$ with $c$ sufficiently large yields $\|(\mathbf{u}_n^\top \widetilde{J}_m^{-1}(\tau))^{(\mathcal{I}(\mathbf{u}_n,D))} - \mathbf{u}_n^\top \widetilde{J}_m^{-1}(\tau)\| \leq n^{-1} 3^{-1}\xi_m^{-1}$ where we used (A.5). This completes the proof. □

*S.1.2.1. Technical details for the proof of Theorem 5.2*

**Lemma S.1.4.** *Under the assumptions of Theorem 5.2 we have for sufficiently large $n$ and any $\kappa_n \ll n/\xi_m^2$,*

$$P\Big(\sup_{\tau,x}|\mathbf{B}(x)^\top(\widehat{\boldsymbol{\beta}}(\tau) - \boldsymbol{\beta}_n(\tau))| \geq C\Big(\frac{\xi_m \log n + \xi_m \kappa_n^{1/2}}{n^{1/2}} + \widetilde{c}_n^2\Big)\Big) \leq (m+1)e^{-\kappa_n}.$$

*where the constant $C$ does not depend on $n$.*

**Proof of Lemma S.1.4.** Apply (A.5) with $\mathbf{a} = \mathbf{B}(x)$ to obtain

$$\|\mathbf{B}(x)^\top \widetilde{J}_m^{-1}(\tau) - (\mathbf{B}(x)^\top \widetilde{J}_m^{-1}(\tau))^{(\mathcal{I}(\mathbf{B}(x),D))}\| \lesssim m\xi_m \gamma^D,$$

where $\mathcal{I}(\mathbf{B}(x), D)$ is defined as (A.2), and $\gamma \in (0,1)$ is a constant independent of $n$. Next observe the decomposition

$$\mathbf{B}(x)^\top(\widehat{\boldsymbol{\beta}}(\tau) - \boldsymbol{\beta}_n(\tau)) = -n^{-1/2}(\mathbf{B}(x)^\top \widetilde{J}_m^{-1}(\tau))^{(\mathcal{I}(\mathbf{B}(x),D))}\mathbb{G}_n(\psi(\cdot;\widehat{\boldsymbol{\beta}}(\tau),\tau)) + \sum_{k=1}^{4} r_{n,k}(\tau,x)$$

where

$$\begin{aligned}
r_{n,1}(\tau,x) &:= (\mathbf{B}(x)^\top \widetilde{J}_m^{-1}(\tau))^{(\mathcal{I}(\mathbf{B}(x),D))} \mathbb{P}_n \psi(\cdot;\widehat{\boldsymbol{\beta}}(\tau),\tau), \\
r_{n,2}(\tau,x) &:= -(\mathbf{B}(x)^\top \widetilde{J}_m^{-1}(\tau))\Big(\mu(\widehat{\boldsymbol{\beta}}(\tau),\tau) - \mu(\boldsymbol{\beta}_n(\tau),\tau) - \widetilde{J}_m(\tau)(\widehat{\boldsymbol{\beta}}(\tau) - \boldsymbol{\beta}_n(\tau))\Big), \\
r_{n,3}(\tau,x) &:= -(\mathbf{B}(x)^\top \widetilde{J}_m^{-1}(\tau))\mu(\boldsymbol{\beta}_n(\tau),\tau)).
\end{aligned}$$

and

$$r_{n,4}(\tau,x) := \Big(\mathbf{B}(x)^\top \widetilde{J}_m^{-1}(\tau) - (\mathbf{B}(x)^\top \widetilde{J}_m^{-1}(\tau))^{(\mathcal{I}(\mathbf{B}(x),D))}\Big)\Big(\mathbb{P}_n \psi(\cdot;\widehat{\boldsymbol{\beta}}(\tau),\tau) - n^{-1/2}\mathbb{G}_n(\psi(\cdot;\widehat{\boldsymbol{\beta}}(\tau),\tau))\Big).$$

Letting $D = c \log n$ for a sufficiently large constant $c$, (S.1.15) yields $\sup_{\tau,x}|r_{n,4}(\tau,x)| \leq n^{-1}$ almost surely. Lemma S.1.7 yields the bound

$$\sup_{x,\tau}|r_{n,1}(\tau,x)| \lesssim \frac{\xi_m^2 \log n}{n} \quad a.s. \tag{S.1.16}$$

Let $\delta_n := \sup_{x,\tau}|\mathbf{B}(x)^\top \widehat{\boldsymbol{\beta}}(\tau) - \mathbf{B}(x)^\top \boldsymbol{\beta}_n(\tau)|$. From Lemma S.1.6 we obtain under (L)

$$\sup_{x,\tau}|r_{n,2}(\tau,x)| \leq \delta_n^2 \sup_{x,\tau}\mathbb{E}\big[|\mathbf{B}(x)^\top \widetilde{J}_m^{-1}(\tau)\mathbf{B}|\big] \asymp \delta_n^2, \tag{S.1.17}$$

where $\sup_{x,\tau}\mathbb{E}\big[|\mathbf{B}(x)^\top \widetilde{J}_m^{-1}(\tau)\mathbf{B}|\big] = O(1)$ by assumption (L). Finally, note that

$$n^{-1/2}(\mathbf{B}(x)^\top \widetilde{J}_m^{-1}(\tau))^{(\mathcal{I}(\mathbf{B}(x),D))}\mathbb{G}_n(\psi(\cdot;\widehat{\boldsymbol{\beta}}(\tau),\tau))$$
$$= n^{-1/2}(\mathbf{B}(x)^\top \widetilde{J}_m^{-1}(\tau))^{(\mathcal{I}(\mathbf{B}(x),D))}\mathbb{G}_n(\psi^{(\mathcal{I}'(\mathbf{B}(x),D))}(\cdot;\widehat{\boldsymbol{\beta}}(\tau),\tau))$$

where $\mathcal{I}'(\mathbf{B}(x), D)$ is defined as (A.3). This yields

$$\sup_{\tau,x}\Big|n^{-1/2}(\mathbf{B}(x)^\top \widetilde{J}_m^{-1}(\tau))^{(\mathcal{I}(\mathbf{B}(x),D))}\mathbb{G}_n(\psi(\cdot;\widehat{\boldsymbol{\beta}}(\tau),\tau))\Big| \lesssim \xi_m \sup_{x\in\mathcal{X}}\|\mathbb{P}_n - P\|_{\mathcal{G}_1(\mathcal{I}(\mathbf{B}(x),D),\mathcal{I}'(\mathbf{B}(x),D))}$$

By the definition of $\mathcal{I}(\mathbf{B}(x), D), \mathcal{I}'(\mathbf{B}(x), D)$, the supremum above ranges over at most $m$ distinct terms. Additionally, $\sup_{x\in\mathcal{X}}|\mathcal{I}(\mathbf{B}(x), c\log n)| + |\mathcal{I}'(\mathbf{B}(x), c\log n)| \lesssim \log n$. Thus Lemma S.1.5 yields

$$P\Big(\sup_{\tau,x}\Big|n^{-1/2}(\mathbf{B}(x)^\top \widetilde{J}_m^{-1}(\tau))^{(\mathcal{I}(\mathbf{B}(x),c\log n))}\mathbb{G}_n(\psi(\cdot;\widehat{\boldsymbol{\beta}}(\tau),\tau))\Big| \geq C\Big(\frac{\xi_m^2(\log n)^2}{n}\Big)^{1/2} + C\frac{\xi_m \kappa_n^{1/2}}{n^{1/2}}\Big) \leq me^{-\kappa_n}. \tag{S.1.18}$$

Finally, from the definition of $\boldsymbol{\beta}_n$ as minimizer we obtain

$$\begin{aligned} |r_{n,3}(\tau,x)| &= \left|(\mathbf{B}(x)^\top \widetilde{J}_m^{-1}(\tau))\mu(\boldsymbol{\beta}_n(\tau),\tau))\right| \\ &= \left|(\mathbf{B}(x)^\top \widetilde{J}_m^{-1}(\tau))\mathbb{E}[\mathbf{B}(\mathbf{1}\{Y \leq \boldsymbol{\beta}_n(\tau)^\top \mathbf{B}\} - \tau)]\right| \\ &= \left|(\mathbf{B}(x)^\top \widetilde{J}_m^{-1}(\tau))\mathbb{E}[\mathbf{B}(F_{Y|X}(\boldsymbol{\beta}_n(\tau)^\top \mathbf{B}|X) - F_{Y|X}(Q(X;\tau)|X))]\right| \\ &\leq \frac{1}{2}\overline{f'}\widetilde{c}_n^2 O(1) \end{aligned} \qquad (\text{S.1.19})$$

where the last line follows after a Taylor expansion and the fact that $\mathbb{E}[\mathbf{B} f_{Y|X}(Q(X;\tau)|X)(\boldsymbol{\beta}_n^\top \mathbf{B} - Q(X;\tau))] = 0$ from the definition of $\boldsymbol{\beta}_n$ and making use of (L). Combining this with (S.1.16) - (S.1.18) yields

$$\delta_n \leq C\Big[\Big(\frac{\xi_m^2(\log n)^2}{n}\Big)^{1/2} + \frac{\xi_m \kappa_n^{1/2}}{n^{1/2}} + \frac{\xi_m^2 \log n}{n} + \delta_n^2 + \widetilde{c}_n^2\Big]$$

with probability at least $me^{-\kappa_n}$. By Lemma S.1.2 we have $P(\delta_n \geq 1/(2C)) \leq e^{-\kappa_n}$ for any $\kappa_n$ satisfying $\xi_m \kappa_n \gg n^{-1}$ This yields the assertion. $\square$

**Lemma S.1.5.** *Let $\mathcal{Z} := \{\mathbf{B}(x)|x \in \mathcal{X}\}$ where $\mathcal{X}$ is the support of $X$. For $\mathcal{I}_1, \mathcal{I}_1' \subset \{1,...,m\}$, define the classes of functions*

$$\widetilde{\mathcal{G}}_1(\mathcal{I}_1, \mathcal{I}_1') := \big\{(\mathbf{Z},Y) \mapsto \mathbf{a}^\top \mathbf{Z}^{(\mathcal{I}_1)}(\mathbf{1}\{Y \leq \mathbf{Z}^\top \mathbf{b}^{(\mathcal{I}_1')}\} - \tau)\mathbf{1}\{\|\mathbf{Z}\| \leq \xi_m\}\big|\tau \in \mathcal{T}, \mathbf{b} \in \mathbb{R}^m, \mathbf{a} \in \mathcal{S}^{m-1}\big\}, \quad (\text{S.1.20})$$

$$\widetilde{\mathcal{G}}_2(\delta, \mathcal{I}_1, \mathcal{I}_1') := \big\{(\mathbf{Z},Y) \mapsto \mathbf{a}^\top \mathbf{Z}^{(\mathcal{I}_1)}(\mathbf{1}\{Y \leq \mathbf{Z}^\top \mathbf{b}_1^{(\mathcal{I}_1')}\} - \mathbf{1}\{Y \leq \mathbf{Z}^\top \mathbf{b}_2^{(\mathcal{I}_1')}\})\mathbf{1}\{\mathbf{Z} \in \mathcal{Z}\}\big|$$
$$\mathbf{b}_1, \mathbf{b}_2 \in \mathbb{R}^m, \sup_{\mathbf{v} \in \mathcal{Z}} \|\mathbf{v}^\top \mathbf{b}_1 - \mathbf{v}^\top \mathbf{b}_2\| \leq \delta, \mathbf{a} \in \mathcal{S}^{m-1}\big\}. \quad (\text{S.1.21})$$

*Under assumptions (A1)-(A3) we have*

$$P\Big(\|\mathbb{P}_n - P\|_{\widetilde{\mathcal{G}}_1(\mathcal{I}_1, \mathcal{I}_1')} \geq C\Big[\Big(\frac{\max(|\mathcal{I}_1|,|\mathcal{I}_1'|)}{n}\log \xi_m\Big)^{1/2} + \frac{\max(|\mathcal{I}_1|,|\mathcal{I}_1'|)\xi_m}{n}\log \xi_m + \Big(\frac{\kappa_n}{n}\Big)^{1/2} + \frac{\xi_m \kappa_n}{n}\Big]\Big) \leq e^{-\kappa_n} \qquad (\text{S.1.22})$$

*and for any $\delta_n$ satisfying $\xi_m \delta_n \gg n^{-1}$ we have for sufficiently large $n$ and arbitrary $\kappa_n > 0$*

$$P\Big(\|\mathbb{P}_n - P\|_{\widetilde{\mathcal{G}}_2(\delta, \mathcal{I}_1, \mathcal{I}_1')} \geq C\zeta_n(\delta, \mathcal{I}_1, \mathcal{I}_1', \kappa_n)\Big) \leq e^{-\kappa_n}, \qquad (\text{S.1.23})$$

*where*

$$\zeta_n(t, \mathcal{I}_1, \mathcal{I}_1', \kappa_n) := t^{1/2}\Big(\frac{\max(|\mathcal{I}_1|,|\mathcal{I}_1'|)}{n}\log(\xi_m \vee n)\Big)^{1/2} + \frac{\max(|\mathcal{I}_1|,|\mathcal{I}_1'|)\xi_m}{n}\log(\xi_m \vee n)$$
$$+ n^{-1/2}(t\kappa_n)^{1/2} + n^{-1}\xi_m \kappa_n.$$

**Proof of Lemma S.1.5.** We begin by observing that

$$N\big(\widetilde{\mathcal{G}}_2(\delta, \mathcal{I}_1, \mathcal{I}_1'), L_2(\mathbb{P}_n); \varepsilon\big) \leq \Big(\frac{A\|F\|_{L^2(\mathbb{P}_n)}}{\varepsilon}\Big)^{v_1(m)}, \quad N\big(\widetilde{\mathcal{G}}_1(\mathcal{I}_1, \mathcal{I}_1'), L_2(\mathbb{P}_n); \varepsilon\big) \leq \Big(\frac{A\|F\|_{L^2(\mathbb{P}_n)}}{\varepsilon}\Big)^{v_2(m)},$$

where $v_1(m) = O(\max(|\mathcal{I}_1|,|\mathcal{I}_1'|)), v_2(m) = O(\max(|\mathcal{I}_1|,|\mathcal{I}_1'|))$. The proof of the bound for $\widetilde{\mathcal{G}}_1(\mathcal{I}_1, \mathcal{I}_1')$ now follows by similar arguments as the proof of Lemma S.1.3. For a proof of the second part, note that for $f \in \widetilde{\mathcal{G}}_2$ we have

$$\begin{aligned} \mathbb{E}[f^2] &\leq \sup_{\|\mathbf{a}\|=1} \sup_{\mathbf{b}_1, \mathbf{b}_2 \in R(\delta_n)} \mathbb{E}\Big[(\mathbf{a}^\top \mathbf{B}^{(\mathcal{I}_1)})^2 \mathbf{1}\Big\{|Y - \mathbf{B}^\top \mathbf{b}_1^{(\mathcal{I}_1')}| \leq |\mathbf{B}^\top(\mathbf{b}_1^{(\mathcal{I}_1')} - \mathbf{b}_2^{(\mathcal{I}_1')})|\Big\}\Big] \\ &\leq \sup_{\mathbf{b} \in \mathbb{R}^m} \sup_{\|\mathbf{a}\|=1} \mathbb{E}\Big[\big((\mathbf{a}^{(\mathcal{I}_1)})^\top \mathbf{B}\mathbf{B}^\top \mathbf{a}^{(\mathcal{I}_1)}\big)^2 \mathbf{1}\{|Y - \mathbf{B}^\top \mathbf{b}^{(\mathcal{I}_1')}| \leq \delta_n\}\Big] \\ &\leq 2\overline{f}\delta_n \lambda_{\max}(\mathbb{E}[\mathbf{B}\mathbf{B}^\top]) \end{aligned}$$

where we defined $R(\delta) := \left\{ \mathbf{b}_1, \mathbf{b}_2 \in \mathbb{R}^m, \sup_{\mathbf{v} \in \mathcal{Z}} \|\mathbf{v}^\top \mathbf{b}_1 - \mathbf{v}^\top \mathbf{b}_2\| \leq \delta \right\}$. The rest of the proof follows by similar arguments as the proof of Lemma S.1.3. $\square$

**Lemma S.1.6.** *Under assumptions (A1)-(A3) we have for any* $\mathbf{a}, \mathbf{b} \in \mathbb{R}^m$

$$\left| \mathbf{a}^\top \widetilde{J}_m(\tau)^{-1} \mu(\mathbf{b}, \tau) - \mathbf{a}^\top \widetilde{J}_m(\tau)^{-1} \mu(\boldsymbol{\beta}_n(\tau), \tau) - \mathbf{a}^\top (\mathbf{b} - \boldsymbol{\beta}_n(\tau)) \right|$$
$$\leq \overline{f'} \sup_x |\mathbf{B}(x)^\top \mathbf{b} - \mathbf{B}(x)^\top \boldsymbol{\beta}_n(\tau)|^2 \mathbb{E}[|\mathbf{a}^\top \widetilde{J}_m(\tau)^{-1} \mathbf{B}|].$$

**Proof of Lemma S.1.6.** Note that $\mu'(\boldsymbol{\beta}_n(\tau), \tau) = \mathbb{E}[\mathbf{B}\mathbf{B}^\top f_{Y|X}(\mathbf{B}^\top \boldsymbol{\beta}_n(\tau)|X)] = \widetilde{J}_m(\tau)$. Additionally, we have

$$\mu(\mathbf{b}, \tau) = \mu(\boldsymbol{\beta}_n(\tau), \tau) + \mu'(\bar{\mathbf{b}}, \tau)(\mathbf{b} - \boldsymbol{\beta}_n(\tau)),$$

where $\bar{\mathbf{b}} = \mathbf{b} + \lambda_{\mathbf{b},\tau}(\boldsymbol{\beta}_n(\tau) - \mathbf{b})$ for some $\lambda_{\mathbf{b},\tau} \in [0,1]$. Moreover,

$$\mathbf{a}^\top [\widetilde{J}_m(\tau)^{-1} \mu(\mathbf{b}, \tau) - \widetilde{J}_m(\tau)^{-1} \mu(\boldsymbol{\beta}_n(\tau), \tau) - (\mathbf{b} - \boldsymbol{\beta}_n(\tau))] = \mathbf{a}^\top \widetilde{J}_m(\tau)^{-1} [\mu'(\bar{\mathbf{b}}, \tau) - \widetilde{J}_m(\tau)](\mathbf{b} - \boldsymbol{\beta}_n(\tau))$$

and thus

$$\left| \mathbf{a}^\top \widetilde{J}_m(\tau)^{-1} \mu(\mathbf{b}, \tau) - \mathbf{a}^\top \widetilde{J}_m(\tau)^{-1} \mu(\boldsymbol{\beta}_n(\tau), \tau) - \mathbf{a}^\top (\mathbf{b} - \boldsymbol{\beta}_n(\tau)) \right|$$
$$= \left| \mathbb{E}\left[ \mathbf{a}^\top \widetilde{J}_m(\tau)^{-1} \mathbf{B}\mathbf{B}^\top \left( f_{Y|X}(\mathbf{B}^\top \bar{\mathbf{b}}|X) - f_{Y|X}(\mathbf{B}^\top \boldsymbol{\beta}_n(\tau)|X) \right) (\mathbf{b} - \boldsymbol{\beta}_n(\tau)) \right] \right|$$
$$\leq \overline{f'} \mathbb{E}\left[ \left| \mathbf{a}^\top \widetilde{J}_m(\tau)^{-1} \mathbf{B} \right| \left\{ \mathbf{B}^\top (\mathbf{b} - \boldsymbol{\beta}_n(\tau)) \right\}^2 \right]$$
$$\leq \overline{f'} \sup_x |\mathbf{B}(x)^\top \mathbf{b} - \mathbf{B}(x)^\top \boldsymbol{\beta}_n(\tau)|^2 \mathbb{E}[|\mathbf{a}^\top \widetilde{J}_m(\tau)^{-1} \mathbf{B}|].$$

$\square$

**Lemma S.1.7.** *Under assumptions (A1)-(A3) and (L) we have for any* $\mathbf{a} \in \mathbb{R}^m$ *having zero entries everywhere except at L consecutive positions:*

$$\left| \mathbf{a}^\top \mathbb{P}_n \psi(\cdot; \widehat{\boldsymbol{\beta}}(\tau), \tau) \right| \leq \frac{(L+2r) \|\mathbf{a}\| \xi_m}{n}.$$

**Proof of Lemma S.1.7.** From standard arguments of the optimization condition of quantile regression (p.35 of Koenker (2005), also see equation (2.2) on p.224 of Knight (2008)), we know that for any $\tau \in \mathcal{T}$,

$$\mathbb{P}_n \psi(\cdot; \widehat{\boldsymbol{\beta}}(\tau), \tau) = \frac{1}{n} \sum_{i=1}^n \mathbf{B}_i \left( \mathbf{1}\{Y_i \leq \mathbf{B}_i^\top \widehat{\boldsymbol{\beta}}(\tau)\} - \tau \right) = \frac{1}{n} \sum_{i \in H_\tau} v_i \mathbf{B}_i$$

where $v_i \in [-1, 1]$ and $H_\tau = \{i : Y_i = \mathbf{B}_i^\top \widehat{\boldsymbol{\beta}}(\tau)\}$. Since $\mathbf{a}$ has at most $L$ non-zero entries, the dimension of the subspace spanned by $\{\mathbf{B}_i : \mathbf{a}^\top \mathbf{B}_i \neq 0\}$ is at most $L + 2r$ [each vector $\mathbf{B}_i$ by construction has at most $r$ nonzero entries and all of those entries are consecutive]. Since the conditional distribution of $Y$ given covariates has a density, the data are in general position almost surely, i.e. no more than $k$ of the points $(\mathbf{B}_i, Y_i)$ lie in any $k$-dimensional linear space, it follows that the cardinality of the set $H_\tau \cap \{i : \mathbf{a}^\top \mathbf{B}_i \neq 0\}$ is bounded by $L + 2r$. The assertion follows after an elementary calculation. $\square$

### S.1.3. Proof of Theorem 5.4

The statement follows from Theorem 5.1 if we prove that the vector $\boldsymbol{\gamma}_n^\dagger(\tau)$ satisfies

$$\sup_{\tau \in \mathcal{T}} \left\| \mu(\boldsymbol{\gamma}_n^\dagger(\tau); \tau) \right\| = O(\xi_m c_n^{\dagger 2}) \tag{S.1.24}$$

as $c_n^\dagger = o(\xi_m^{-1})$ in Condition (C1), and establish the identity in (5.10). For the identity (5.10), we first observe the representation

$$J_m(\tau) = \begin{pmatrix} M_1(\tau) + A(\tau)M_2(\tau)A(\tau)^\top & A(\tau)M_2(\tau) \\ M_2(\tau)A(\tau)^\top & M_2(\tau) \end{pmatrix}, \quad (S.1.25)$$

which follows from (3.5) and

$$\mathbb{E}[(V - A(\tau)\widetilde{\mathbf{Z}}(W))\widetilde{\mathbf{Z}}(W)^\top f_{Y|X}(Q(X;\tau)|X)] = 0, \text{ for all } \tau \in \mathcal{T}.$$

To simplify the notations, we suppress the argument in $\tau$ in the following matrix calculations. Recall the following identity for the inverse of $2 \times 2$ block matrix (see equation (6.0.8) on p.165 of Puntanen and Styan (2005))

$$\begin{pmatrix} A & B \\ C & D \end{pmatrix}^{-1} = \begin{pmatrix} (A - BD^{-1}C)^{-1} & -(A - BD^{-1}C)^{-1}BD^{-1} \\ -D^{-1}C(A - BD^{-1}C)^{-1} & D^{-1} + D^{-1}C(A - BD^{-1}C)^{-1}BD^{-1} \end{pmatrix}.$$

Identifying the blocks in the representation (S.1.25) with the blocks in te above representation yields the result after some simple calculations. For a proof of (S.1.24) observe that

$$\mu(\boldsymbol{\gamma}_n^\dagger(\tau);\tau) = \mathbb{E}[(V^\top, \widetilde{\mathbf{Z}}(W)^\top)^\top (F_{Y|X}(\boldsymbol{\gamma}_n^\dagger(\tau)^\top \mathbf{Z}|X) - \tau)]$$

Now on one hand we have, uniformly in $\tau \in \mathcal{T}$,

$$\left\|\mathbb{E}[V(F_{Y|X}(\boldsymbol{\alpha}(\tau)^\top V + \boldsymbol{\beta}_n^\dagger(\tau)^\top \widetilde{\mathbf{Z}}(W)|X) - \tau)]\right\|$$
$$= \left\|\mathbb{E}[V f_{Y|X}(Q(X;\tau)|X)(Q(X;\tau) - \boldsymbol{\alpha}(\tau)^\top V - \boldsymbol{\beta}_n^\dagger(\tau)^\top \widetilde{\mathbf{Z}}(W))]\right\| + O(c_n^{\dagger 2})$$
$$= \left\|\mathbb{E}[V f_{Y|X}(Q(X;\tau)|X)(h(W;\tau) - \boldsymbol{\beta}_n^\dagger(\tau)^\top \widetilde{\mathbf{Z}}(W))]\right\| + O(c_n^{\dagger 2})$$
$$= \left\|\mathbb{E}[(V - h_{VW}(W;\tau) + h_{VW}(W;\tau))f_{Y|X}(Q(X;\tau)|X)(h(W;\tau) - \boldsymbol{\beta}_n^\dagger(\tau)^\top \widetilde{\mathbf{Z}}(W))]\right\| + O(c_n^{\dagger 2})$$
$$= \left\|\mathbb{E}[h_{VW}(W;\tau)f_{Y|X}(Q(X;\tau)|X)(h(W;\tau) - \boldsymbol{\beta}_n^\dagger(\tau)^\top \widetilde{\mathbf{Z}}(W))]\right\| + O(c_n^{\dagger 2})$$
$$= \left\|\mathbb{E}[(h_{VW}(W;\tau) - A(\tau)\widetilde{\mathbf{Z}}(W) + A(\tau)\widetilde{\mathbf{Z}}(W))f_{Y|X}(Q(X;\tau)|X)(h(W;\tau) - \boldsymbol{\beta}_n^\dagger(\tau)^\top \widetilde{\mathbf{Z}}(W))]\right\| + O(c_n^{\dagger 2})$$
$$= \left\|\mathbb{E}[(h_{VW}(W;\tau) - A(\tau)\widetilde{\mathbf{Z}}(W))f_{Y|X}(Q(X;\tau)|X)(h(W;\tau) - \boldsymbol{\beta}_n^\dagger(\tau)^\top \widetilde{\mathbf{Z}}(W))]\right\| + O(c_n^{\dagger 2})$$
$$= O(c_n^{\dagger 2} + \lambda_n c_n^\dagger).$$

Here, the first equation follows after a Taylor expansion taking into account that, by the definition of the conditional quantile function, $F_{Y|X}(Q(X;\tau)|X) \equiv \tau$. The fourth equality is a consequence of (3.5), the sixth equality follows since

$$\mathbb{E}[\widetilde{\mathbf{Z}}(W)f_{Y|X}(Q(X;\tau)|X)(h(W;\tau) - \boldsymbol{\beta}_n^\dagger(\tau)^\top \widetilde{\mathbf{Z}}(W))] = 0 \quad (S.1.26)$$

by the definition of $\boldsymbol{\beta}_n^\dagger(\tau)$ as minimizer. The last line follows by the Cauchy-Schwarz inequality. On the other hand

$$\mathbb{E}[\widetilde{\mathbf{Z}}(W)(F_{Y|X}(\boldsymbol{\alpha}(\tau)^\top V + \boldsymbol{\beta}_n^\dagger(\tau)^\top \widetilde{\mathbf{Z}}(W)|X) - \tau)]$$
$$= \mathbb{E}[\widetilde{\mathbf{Z}}(W)f_{Y|X}(Q(X;\tau)|X)(h(W;\tau) - \boldsymbol{\beta}_n^\dagger(\tau)^\top \widetilde{\mathbf{Z}}(W))]$$
$$+ \frac{1}{2}\mathbb{E}[\widetilde{\mathbf{Z}}(W)f'_{Y|X}(\zeta(X;\tau)|X)(h(W;\tau) - \boldsymbol{\beta}_n^\dagger(\tau)^\top \widetilde{\mathbf{Z}}(W))^2].$$

By (S.1.26), the first term in the representation above is zero, and the norm of the second term is of the order $O(\xi_m c_n^{\dagger 2})$. This completes the proof. $\square$

## APPENDIX S.2: Auxiliary Results

### S.2.1. Results on empirical process theory

In this section, we collect some basic results from empirical process theory needed in our proofs. Denote by $\mathcal{G}$ a class of functions that satisfies $|f(x)| \leq F(x) \leq U$ for every $f \in \mathcal{G}$ and let $\sigma^2 \geq \sup_{f \in \mathcal{G}} Pf^2$. Additionally, let for some $A > 0, V > 0$ and all $\varepsilon > 0$,

$$N(\varepsilon, \mathcal{G}, L_2(\mathbb{P}_n)) \leq \Big(\frac{A\|F\|_{L^2(\mathbb{P}_n)}}{\varepsilon}\Big)^V. \tag{S.2.1}$$

Note that if $\mathcal{G}$ is a VC-class, then $V$ is the VC-index of the set of subgraphs of functions in $\mathcal{G}$. In that case, the symmetrization inequality and inequality (2.2) from Koltchinskii (2006) yield

$$\mathbb{E}\|\mathbb{P}_n - P\|_{\mathcal{G}} \leq c_0 \Big[\sigma\Big(\frac{V}{n}\log\frac{A\|F\|_{L^2(P)}}{\sigma}\Big)^{1/2} + \frac{VU}{n}\log\frac{A\|F\|_{L^2(P)}}{\sigma}\Big] \tag{S.2.2}$$

for a universal constant $c_0 > 0$ provided that $1 \geq \sigma^2 > \text{const} \times n^{-1}$ [in fact, the inequality in Koltchinskii (2006) is for $\sigma^2 = \sup_{f \in \mathcal{G}} Pf^2$. However, this is not a problem since we can replace $\mathcal{G}$ by $\mathcal{G}\sigma/(\sup_{f \in \mathcal{G}} Pf^2)^{1/2}$]. The second inequality (a refined version of Talagrand's concentration inequality) states that for any countable class of measurable functions $\mathcal{F}$ with elements mapping into $[-M, M]$

$$P\Big\{\|\mathbb{P}_n - P\|_{\mathcal{F}} \geq 2\mathbb{E}\|\mathbb{P}_n - P\|_{\mathcal{F}} + c_1 n^{-1/2}\Big(\sup_{f \in \mathcal{F}} Pf^2\Big)^{1/2}\sqrt{t} + n^{-1}c_2 Mt\Big\} \leq e^{-t}, \tag{S.2.3}$$

for all $t > 0$ and universal constants $c_1, c_2 > 0$. This is a special case of Theorem 3 in Massart (2000) [in the notation of that paper, set $\varepsilon = 1$].

**Lemma S.2.1** (Lemma 7.1 of Kley et al. (2016)). *Let $\{\mathbb{G}_t : t \in T\}$ be a separable stochastic process with $\|\mathbb{G}_s - \mathbb{G}_t\|_{\Psi} \leq Cd(s,t)$ ($\|\cdot\|_{\Psi}$ is defined in (A.32)) for all $s, t$ satisfying $d(s,t) \geq \bar{\omega}/2 \geq 0$. Denote by $D(\epsilon, d)$ the packing number of the metric space $(T, d)$. Then, for any $\delta > 0$, $\omega \geq \bar{\omega}$, there exists a random variable $S_1$ and a constant $K < \infty$ such that*

$$\sup_{d(s,t) \leq \delta} |\mathbb{G}_s - \mathbb{G}_t| \leq S_1 + 2\sup_{d(s,t) \leq \bar{\omega}, t \in \widetilde{T}} |\mathbb{G}_s - \mathbb{G}_t|, \tag{S.2.4}$$

*where the set $\widetilde{T}$ contains at most $D(\bar{\omega}, d)$ points, and $S_1$ satisfies*

$$\|S_1\|_{\Psi} \leq K\Big[\int_{\bar{\omega}/2}^{\omega} \Psi^{-1}\big(D(\epsilon, d)\big)d\epsilon + (\delta + 2\bar{\omega})\Psi^{-1}\big(D^2(\omega, d)\big)\Big] \tag{S.2.5}$$

$$P(|S_1| > x) \leq \Big(\Psi\Big\{x\Big[8K\Big(\int_{\bar{\omega}/2}^{\omega} \Psi^{-1}\big(D(\epsilon, d)\big)d\epsilon + (\delta + 2\bar{\omega})\Psi^{-1}\big(D^2(\omega, d)\big)\Big)\Big]^{-1}\Big\}\Big)^{-1}. \tag{S.2.6}$$

### S.2.2. Covering number calculation

A few useful lemmas on covering number are given in this section.

**Lemma S.2.2.** *Suppose $\mathcal{F}$ and $\mathcal{G}$ are two function classes with envelopes $F$ and $G$.*

1. *The class $\mathcal{F} - \mathcal{G} := \{f - g | f \in \mathcal{F}, g \in \mathcal{G}\}$ with envelope $F + G$ and*

$$\sup_Q N\big(\epsilon\|F+G\|_{Q,2}, \mathcal{F}-\mathcal{G}, L_2(Q)\big) \leq \sup_Q N\Big(\epsilon\frac{\|F\|_{Q,2}}{\sqrt{2}}, \mathcal{F}, L_2(Q)\Big)\sup_Q N\Big(\epsilon\frac{\|G\|_{Q,2}}{\sqrt{2}}, \mathcal{G}, L_2(Q)\Big), \tag{S.2.7}$$

2. The class $\mathcal{F} \cdot \mathcal{G} := \{fg : f \in \mathcal{F}, g \in \mathcal{G}\}$ with envelope $FG$ and

$$\sup_Q N\big(\epsilon \|FG\|_{Q,2}, \mathcal{F} \cdot \mathcal{G}, L_2(Q)\big) \leq \sup_Q N\left(\frac{\epsilon \|F\|_{Q,2}}{2}, \mathcal{F}, L_2(Q)\right) \sup_Q N\left(\frac{\epsilon \|G\|_{Q,2}}{2}, \mathcal{G}, L_2(Q)\right), \quad \text{(S.2.8)}$$

where the suprema are taken over the appropriate subsets of all finitely discrete probability measures $Q$.

**Proof of Lemma S.2.2.**

1. It is obvious that $|f - g| \leq |f| + |g| \leq F + G$ for any $f \in \mathcal{F}$ and $g \in \mathcal{G}$. Hence, $F + G$ is an envelop for the function class $\mathcal{F} - \mathcal{G}$. Suppose that $A = \{f_1, ... f_J\}$ and $B = \{g_1, ..., g_K\}$ are the centers of $\frac{\epsilon \|F\|_{Q,2}}{\sqrt{2}}$-net for $\mathcal{F}$ and $\frac{\epsilon \|G\|_{Q,2}}{\sqrt{2}}$-net for $\mathcal{F}$ and $\mathcal{G}$ respectively. For any $f - g$, there exists $f_j$ and $g_k$ such that

$$\|(f_j - g_k) - (f - g)\|_{Q,2}^2 = \|(f_j - f) - (g_k - g)\|_{Q,2}^2 \leq 2\big(\|f_j - f\|_{Q,2}^2 + \|g_k - g\|_{Q,2}^2\big)$$
$$\leq \epsilon^2(\|F\|_{Q,2}^2 + \|G\|_{Q,2}^2) \leq \epsilon^2 \|F + G\|_{Q,2}^2,$$

where the last inequality follows from the fact that both $F$ and $G$ are nonnegative. Hence, $\{f_j + g_k : 1 \leq j \leq J, 1 \leq k \leq K\}$ forms an $\epsilon \|F + G\|_{Q,2}$-net for the class $\mathcal{F} - \mathcal{G}$, with cardinality $JK$.

2. See Lemma 6 of Belloni et al. (2016).

□

For any fixed vector $\mathbf{u} \in \mathbb{R}^m$ and $\delta > 0$, recall the function classes

$$\mathcal{G}_3(\mathbf{u}) = \big\{(\mathbf{Z}, Y) \mapsto \mathbf{u}^\top J_m(\tau)^{-1}\mathbf{Z} \big| \tau \in \mathcal{T}\big\},$$
$$\mathcal{G}_4 = \big\{(X, Y) \mapsto \mathbf{1}\{Y_i \leq Q(X; \tau)\} - \tau \big| \tau \in \mathcal{T}\big\},$$
$$\mathcal{G}_6(\mathbf{u}, \delta) = \big\{(\mathbf{Z}, Y) \mapsto \mathbf{u}^\top \{J_m(\tau_1)^{-1} - J_m(\tau_2)^{-1}\}\mathbf{Z} \big| \tau_1, \tau_2 \in \mathcal{T}, |\tau_1 - \tau_2| \leq \delta\big\},$$
$$\mathcal{G}_7(\delta) = \big\{(X, Y) \mapsto \mathbf{1}\{Y_i \leq Q(X, \tau_1)\} - \mathbf{1}\{Y_i \leq Q(X, \tau_2)\} - (\tau_1 - \tau_2) \big| \tau_1, \tau_2 \in \mathcal{T}, |\tau_1 - \tau_2| \leq \delta\big\}.$$

Recall the following Lipschitz continuity property of $J_m^{-1}(\tau)$ by Lemma 13 of Belloni et al. (2016): for $\tau_1, \tau_2 \in \mathcal{T}$,

$$\|J_m^{-1}(\tau_1) - J_m^{-1}(\tau_2)\| \leq \frac{\bar{f}'}{f_{\min}} |\tau_1 - \tau_2| \Big(\inf_{\tau \in \mathcal{T}} \lambda_{\min}(J_m(\tau))\Big)^{-2} \lambda_{\max}(\mathbb{E}[\mathbf{Z}\mathbf{Z}^\top])$$
$$:= C_0 \inf_{\tau \in \mathcal{T}} \lambda_{\min}(J_m(\tau))^{-1} |\tau_1 - \tau_2|. \quad \text{(S.2.9)}$$

where $C_0 = \frac{\bar{f}'}{f_{\min}} \frac{\lambda_{\max}(\mathbb{E}[\mathbf{Z}\mathbf{Z}^\top])}{\inf_{\tau \in \mathcal{T}} \lambda_{\min}(J_m(\tau))}$.

**Lemma S.2.3.** $\mathcal{G}_3(\mathbf{u})$ has an envelope $G_3(\mathbf{Z}) = \|\mathbf{u}\| \xi_m [\inf_{\tau \in \mathcal{T}} \lambda_{\min}(J_m(\tau))]^{-1}$ and

$$N(\epsilon \|G_3\|_{L_2(Q)}, \mathcal{G}_3(\mathbf{u}), L_2(Q)) \leq \frac{C_0}{\epsilon},$$

where $C_0 = \frac{\bar{f}'}{f_{\min}} \frac{\lambda_{\max}(\mathbb{E}[\mathbf{Z}\mathbf{Z}^\top])}{\inf_{\tau \in \mathcal{T}} \lambda_{\min}(J_m(\tau))} < \infty$, for any probability measure $Q$ and $\mathbf{u}$.

**Proof of Lemma S.2.3.** By (S.2.9), for any $\tau_1, \tau_2 \in \mathcal{T}$ and $\mathbf{u}$,

$$\big|\mathbf{u}^\top J_m(\tau_1)^{-1}\mathbf{Z} - \mathbf{u}^\top J_m(\tau_2)^{-1}\mathbf{Z}\big| \leq \|\mathbf{u}\|\xi_m \frac{\bar{f}'}{f_{\min}} \lambda_{\max}(\mathbb{E}[\mathbf{Z}\mathbf{Z}^\top])[\inf_{\tau \in \mathcal{T}} \lambda_{\min}(J_m(\tau))]^{-2} |\tau_1 - \tau_2|$$
$$= C_0 \|G_3\|_{L_2(Q)} |\tau_1 - \tau_2|.$$

Applying the relation of the covering and bracketing number on p.84 and Theorem 2.7.11 of van der Vaart and Wellner (1996) yields for each $\mathbf{u}$ and any probability measure $Q$,

$$N\big(\epsilon \|G_3\|_{L_2(Q)}, \mathcal{G}_3(\mathbf{u}), L_2(Q)\big) \leq N_{[\,]}\big(2\epsilon \|G_3\|_{L_2(Q)}, \mathcal{G}_3(\mathbf{u}), L_2(Q)\big) \leq N\big(\epsilon, \mathcal{T}, |\cdot|\big) \leq \frac{C_0}{\epsilon}.$$

□

**Lemma S.2.4.** *We have the following results:*

1. $\mathcal{G}_4$ *is a VC-class with VC index 2.*
2. *The envelopes for $\mathcal{G}_6(\mathbf{u}, \delta)$ and $\mathcal{G}_7$ are $G_6 = \xi_m^2 [\inf_{\tau \in \mathcal{T}} \lambda_{\min}(J_m(\tau))]^{-1} C_0 \delta$ and $G_7 = 2$. Furthermore, it holds for any fixed $x$ and $\delta \leq |\mathcal{T}|$ that*

$$N(\epsilon \|G_6\|_{L_2(Q)}, \mathcal{G}_6(\mathbf{u}, \delta), L_2(Q)) \leq 2\left(\frac{C_0}{\epsilon}\right)^2, \quad (S.2.10)$$

$$N(\epsilon \|G_7\|_{L_2(Q)}, \mathcal{G}_7(\delta), L_2(Q)) \leq \left(\frac{A_7}{\epsilon}\right)^4, \quad (S.2.11)$$

*where $A_7$ is a universal constant and $Q$ is an arbitrary probability measure.*

**Proof of Lemma S.2.4.**

1. Due to the fact that $Q(X; \tau)$ is monotone in $\tau$, it can be argued with basic VC subgraph argument that $\mathcal{G}_4$ has VC index 2, under the definition given in p.135 of van der Vaart and Wellner (1996).
2. By (S.2.9), the envelope for $\mathcal{G}_6(\mathbf{u}_n, \delta)$ is $\|\mathbf{u}\|\xi_m \lambda_{\min}(J_m(\tau))^{-1} C_0 \delta$. The envelope for $\mathcal{G}_7(\delta)$ is obvious. By the fact that $\mathcal{G}_7(\delta) \subset \mathcal{G}_4 - \mathcal{G}_4$ and the covering number of $\mathcal{G}_4$ (implied by Theorem 2.6.7 of van der Vaart and Wellner (1996)), (S.2.11) thus follows by (S.2.7) of Lemma S.2.2. As for (S.2.10), we note that $\mathcal{G}_6(\mathbf{u}, \delta) \subset \mathcal{G}_3(\mathbf{u}) - \mathcal{G}_3(\mathbf{u})$. Then, (S.2.7) of Lemma S.2.2 and Lemma S.2.3 imply

$$N(\epsilon\|G_6\|_{L_2(Q)}, \mathcal{G}_6(\mathbf{u}, \delta), L_2(Q)) \leq N\left(\frac{\epsilon}{\sqrt{2}}\|G_3\|_{L_2(Q)}, \mathcal{G}_3(\mathbf{u}), L_2(Q)\right)^2 \leq 2\left(\frac{C_0}{\epsilon}\right)^2,$$

where $Q$ is an arbitrary probability measure.

□